\documentclass[review,onefignum,onetabnum,dvipsnames]{siamart171218}
\usepackage{chngcntr}
\counterwithin{figure}{subsection}

\counterwithin{table}{subsection}
\counterwithin{equation}{subsection}

\counterwithin{algorithm}{subsection}

\usepackage{titlesec}

\def \vu{ \boldsymbol{u} }

\def \ve{ \mathbf{e} }

\def \vu{ \boldsymbol{u} }

\def \ve{ \boldsymbol{e} }

\def \vu{ \boldsymbol{u} }

\newcommand{\N}{\mathbbm{N}}
\newcommand{\R}{\mathbbm{R}}

\DeclareMathOperator*{\argmin}{argmin}
\newcommand{\bs}[1]{\boldsymbol{#1}}
\renewcommand{\hat}[1]{\widehat{#1}}

\usepackage{amsfonts, amssymb,bbm }
\usepackage{graphicx}
\usepackage{mathrsfs,mathtools}
\usepackage{epstopdf}
\usepackage{algorithm}
\usepackage{algorithmic}
\usepackage{cleveref}
\usepackage{enumitem}

\setlist[enumerate,1]{label=(\arabic*),ref=(\arabic*)}

\ifpdf
  \DeclareGraphicsExtensions{.eps,.pdf,.png,.jpg}
\else
  \DeclareGraphicsExtensions{.eps}
\fi

\theorembodyfont{\upshape}

\newsiamremark{remark}{Remark}
\newsiamremark{hypothesis}{Hypothesis}
\crefname{hypothesis}{Hypothesis}{Hypotheses}
\newsiamthm{claim}{Claim}

\headers{Structure-preserving Nonlinear Filtering: CG \& DG}{Vidhi Zala, Robert M. Kirby, and Akil Narayan}

\title{Convex Optimization-Based Structure-Preserving Filter for Multidimensional Finite Element Simulations]{Structure-preserving convex optimization for multidimensional finite element simulations}\thanks{Accepted by editors (JCP) 07/08/2023.
}}

\author{Vidhi Zala\thanks{Scientific Computing and Imaging Institute and School of Computing, University of Utah, Salt Lake City, UT 84112
  (\email{vidhi.zala@utah.edu}\url{}).}
  \and Robert M. Kirby\thanks{
Scientific Computing and Imaging Institute and School of Computing, University of Utah, Salt Lake City, UT 84112
  (\email{kirby@cs.utah.edu}\url{ }).}
\and Akil Narayan\thanks{
Scientific Computing and Imaging Institute and Department of Mathematics, University of
Utah, Salt Lake City, UT 84112
  (\email{akil@sci.utah.edu}\url{}).}
}

\usepackage{amsopn}

\usepackage[export]{adjustbox} 
\usepackage{aliascnt}
\usepackage[]{caption}
\usepackage[]{float}

\usepackage[]{appendix}

\begin{document}
\nolinenumbers
\maketitle
\begin{abstract}
In simulation sciences, capturing the real-world problem features as accurately as possible is desirable. Methods popular for scientific simulations such as the finite element method (FEM) and finite volume method (FVM) use {piecewise} polynomials to  {approximate various characteristics of a problem, such as the concentration profile and the temperature distribution across the domain}. Polynomials are prone to creating artifacts such as Gibbs oscillations while capturing a complex profile. An {efficient and accurate} approach must be applied to deal with such inconsistencies to obtain accurate simulations. This often entails dealing with negative values for the concentration of chemicals, exceeding a percentage value over 100, and other such problems. We consider these inconsistencies in the context of partial differential equations (PDEs). {We propose an innovative filter based on convex optimization to deal with the inconsistencies observed in polynomial-based simulations. In two or three spatial dimensions, additional complexities are involved }in solving the problems related to structure preservation. We present the construction and application of a structure-preserving filter with a focus on multidimensional PDEs. Methods used such as the Barycentric interpolation for polynomial evaluation at arbitrary points in the domain and an optimized root-finder to identify points of interest, improve the filter efficiency, usability, and robustness. Lastly, we present numerical experiments in 2D and 3D using discontinuous Galerkin formulation and demonstrate the filter's efficacy to preserve the desired structure. As a real-world application, implementation of the {mathematical biology} model involving platelet aggregation and blood coagulation has been reviewed and the issues around FEM implementation of the model are resolved by applying the proposed structure-preserving filter. 
\end{abstract}

\section{Introduction}

{A widely used application of mathematical modeling is creating simulations used to predict and analyze different processes. Simulations obtained from numerical solutions provide physically meaningful results if the values of simulation variables at each timestep follow the \textit{structure} of the exact solution.} Structure in this context refers to properties such as non-negative values for a chemical concentration or simulating a quantity requiring an upper bound, e.g., [0,1], and other such aspects of the numerical solution. If the solution fails to comply with these structural restrictions,  the resulting solution is not meaningful. Therefore, designing an approach {to regulating the solution} is crucial for simulations in many applications. Examples include thermodynamics, numerical weather predictions, bio-chemical processes such as platelet aggregation and blood coagulation, and fluid flow problems. \\ 

{Computing solutions to mathematical models described in multiple dimensions} further adds to the complexity of solving a structure-preservation problem.  To implement multidimensional models using numerical methods,  we require approximations of the quantities of interest in the domain.  Popular choices for numerical methods such as the finite element method (FEM) and the finite volume method (FVM) employ {piecewise} polynomial-based approximations.  These approximations are prone to artifacts such as the Gibbs phenomenon, leading to a violation of structure. We focus on the FEM implementation for the {advection-diffusion-reaction (ADR)} problem and design an algorithm that computes a solution at each timestep {adhering to a particular structure }(e.g., positivity and monotonicity). The nonphysical values from such simulations may result in a propagation of nonphysical values, which often cause a blow up, resulting in the failure of the entire simulation. {For example, when the concentration values are expected to be bound within $[0,1]$ but the numerical solution produces a value greater than one, even a small change in the next timestep causes a rapid compounding of the solution.}\\

The problem of structure preservation has been well studied in different fields for applications to various partial differential equations (PDEs) and using different numerical approaches.  The next section provides a brief overview of some existing methods proposed to solve the structure-preservation problem.

\subsection{Review}
Many sophisticated approaches have been considered in the literature to tackle different aspects of the structure-preservation problem.  We make a broad classification into `intrusive' and `nonintrusive' classes of solutions. Many popular approaches generally involve changing either the solution or the PDE and require the imposition of stringent conditions on solver parameters, such as the step-size of the numerical method or the domain shape and granularity. Such methods can be classified as intrusive. On the other hand, some methods constrain the solution obtained from the solver with minimal, often superficial, change to the numerical scheme, i.e., in a nonintrusive way.  \\

Prominent among the intrusive class of approaches is the transformation of structure preservation into a PDE-constrained optimization problem, methods that modify the spatial discretization \cite{zahr2018optimization}, and limiters derived from Karush-Kahan-Tucker (KKT) optimality conditions \cite{KKT}. Another example of the intrusive approach \cite{jcp1} proposes an operator-splitting positivity-preservation method based on the energy dissipation law. \\

The strategies for constraint satisfaction in \cite{a1} consider a version of the problem that is a {specialization of the one solved by the solution proposed in this manuscript, which falls under the nonintrusive class}.  In \cite{a2,allen_bounds-constrained_2021} and extensions \cite{a1,l_1}, the authors explore approaches for positivity preservation on the basis functions derived from Bernstein polynomials such that they are non-negative and possess a partition of unity property.  Therefore, the interpolated solution respects the original bounds at any point in the domain{, nonintrusively}. Other nonintrusive methods prescribe solutions such as limiters or truncation and {modification of the domain (e.g., curvilinear mesh adaptation)} to adhere to the desired solution structures.  One such method is established in \cite{doi:10.1137/17M1116453} where the authors present a new approach for multi-material arbitrary Lagrangian--Eulerian (ALE) hydrodynamics simulations based on high-order finite elements posed on high-order curvilinear meshes in which conservative fields are remapped and additional synchronization steps are introduced for structural preservation. Additionally, many nonintrusive methods perform pre- or postprocessing of intermediate solutions {without modifying the underlying problem to incorporate the constraints or making changes to the domain.} Some nonintrusive approaches to structure preservation can be found in  \cite{TVD,A,B,C,E}.  In \cite{jcp2}, the authors present a positivity-preserving limiter by maintaining cell averages within a certain tolerance.  The work in \cite{laiu_positive_2016} imposes positivity on a discrete grid for kinetic equations via a nonlinear filtering procedure. {Satisfying constraints up to a numerical tolerance is sufficient for many applications. However some applications require the solution to adhere to a strict non-negative structure. For example, in \cite{DELAHAIJE2020116405}, the authors enforce strict, global non-negativity of the diffusion propagator by formulating constraints specific to the propagator model. The algorithm in \cite{van2004fast} employs alternating least squares procedures to implement general linear equality and inequality constraints by reducing it to a non-negativity-constrained least squares (NNLS) problem. The non-negativity constraints based on square root representations have also been of interest in many other works \cite{goh2011nonparametric,cheng2014non}.} We {propose} a nonintrusive method based on a general filtering approach in \cite{p0,p1}. \\

Many of the existing techniques prescribe approaches to preserve the structure that work only for low-order problems. Additionally, an increase in the dimensionality of the problem further exacerbates the issue of accuracy and convergence.  Many approaches struggle to maintain high-order accuracy and convergence robustly, including the filter's ability to preserve structure without {any additional dependence or restrictions imposed on a per-case basis.  This includes problem-specific modifications as well as varying solutions based on other aspects of the problem.  For example,  \cite{repairorig, repair} presents redistribution-based local and global repair methods depending on a sweep through cells to locate and rectify the violating structures.  The final solution depends on the sweep order,  producing different solutions for different sweep orders.  In addition,  this repair procedure preserves element-averaged bound/positivity constraints at particular finite difference nodes. It does not preserve them pointwise, e.g., at every spatial location. The proposed method tries to preserve the constraints at every spatial location.} An ideal approach would achieve a balance between robust structure preservation, maintaining accuracy as well as convergence,  and providing usability by deterministic termination within an acceptable time.  \\



\subsection{Contribution}

We pose the problem of preserving the structure in polynomial-based numerical methods as a convex optimization problem.  Our previous works discuss a filtering algorithm to solve this problem and demonstrate the solution for function approximations \cite{p0} and PDE solutions in 1D \cite{p1}.  The robust design of the problem enables the preservation of different properties of the solution simultaneously (e.g., positivity,  monotonicity,  and boundedness). It transforms the problem into a composite constraint satisfaction problem that can be solved by convex minimization. To address this problem, we propose a filter that can be applied as a postprocessing step to deal with structural inconsistencies. {The novel idea in the design of the filter is the efficiency achieved by applying the filter only to the violating regions in a subdomain and structure-preservation continuously \textit{throughout} the domain, and not just a subset of points within the domain.} {This guarantee of structure-preservation depends on the granularity of the global minimization procedure in 2D and 3D, which is expensive and less reliable compared to its 1D counterpart.  In 1D, the minimum can be found accurately using the confederate matrix approach to find zeros of the derivative of a polynomial. For finding the critical points in higher dimensions, an efficient gradient-based method is required which solves the nonconvex problem of finding the global minimum. In the case of the proposed filter, the global minimum represents a structure violation.  Gradient-based methods are prone to get stuck at a local minimum,  causing it to miss the point of larger structural inconsistency, thus failing to provide strong guarantees of preservation. Another important aspect that affects the outcome of the gradient-based method is the choice of a starting point.  To address this aspect, we propose an investigation of structural inconsistencies on a lattice of points as the initial step in the gradient-based method. The point with the most inconsistency is chosen as a starting point. The accuracy of the minimum found depends on the pattern and granularity of the initial lattice used, which is detailed in \Cref{sec:lattice}}.\\

To streamline the minima-finding stage and reduce the time and space complexity,  we adopt the Barycentric interpolation approach from \cite{Trefethen04},  expanded upon by \cite{laughton2021fast}. Using these techniques, a considerable speed-up in time taken by minimization and overall filter computations is achieved.\\

This paper focuses on designing a structure-preserving filter using a convex optimization approach for 2D and 3D problems on different element types. We will present the formulation derived from \cite{p0} in \Cref{sec:method} and extend the concepts to more complex problems tackled in this paper. The remainder of the paper is organized as follows: \Cref{sec:method} discusses the setup for applications of the filter to time-dependent PDEs in 2D and 3D. \Cref{sec:formalism} introduces the notations, details the design of different building blocks of the filter, and summarizes the 1D problem formulation from \cite{p0,p1}. A procedure for the filtering process developed for multidimensional applications is presented in \Cref{sec:algor}. \Cref{sec:results} describes the numerical results to demonstrate the filter's efficacy in preserving the desired structures in different application scenarios. This section is divided into subsections describing the process of choosing the parameters to run the experiments, the 2D and 3D canonical experiments, and an advection problem on different homogeneous meshes using the discontinuous Galerkin (dG) formulation.  We conclude with a demonstration of the use of the proposed filter on a real-world application: the model of platelet aggregation and blood coagulation problem.

\section{Setup}
\label{sec:method}

In this section, we discuss the setup for filter applications to function approximations and time-dependent PDEs to solve an advection problem using the FEM with method-of-lines discretizations with a focus on the dG formulation. The choice of advection problem is for illustration purposes only. The setup remains the same for any time-dependent PDE solution using a polynomial-based method. The filter behaves as a postprocessing step that {preserves} the desired structure at each timestep {within a defined tolerance and is agnostic to the numerical method used to obtain the solutions. In all our experiments, we consider the numerical zero to be $10^{-7}$. The idea behind this choice is that in all the positivity preserving experiments the gradient-based method stagnates when the difference in minimum found by two consecutive iterations dip below $-10^{-7}$. In such a situation, there are no significant changes to the interpolation coefficients between iterations of the filter, and therefore, the convergence is considerably slow. Certain variations of the filter as investigated in \cite[Section 5.1.2]{p0} can improve the rate of convergence below a particular tolerance. } {\Cref{sec:results} shows numerical examples in higher dimensions using the same setup on element types such as triangle, quadrilateral, hexahedron, tetrahedron, prism, and meshes comprised of these element types.}

\subsection{Detecting and resolving structural inconsistency in 2D and 3D}\label{sec:lattice}

If the polynomial projections lose structural conformity, they are likely to do so at or between the quadrature points. We hope to capture the structural inconsistencies by checking the values on a lattice of quadrature points and the centroids formed by adjacent quadrature points.  {Here, the arbitrary choice of the centroids is one of the many possible options to choose the lattice.} A quadrilateral is used as a sample element to describe the setup; however, the same applies to any canonical element type in 2D and 3D. Let us call such a lattice of points $(X)$. For a quadrilateral, let $Q$ represent the number of quadrature points in one dimension. $X$ is defined by a combination of the quadrature grid on the quadrilateral (total $Q^2$ points and $(Q-1)^2$ centroids formed by the midpoints of the quadrature grid) as defined in \Cref{eq:lattice}.
\begin{equation}\label{eq:lattice}
X = \{Q^2  \textrm{ points in the quadrature grid }\} \cup \{(Q-1)^2 \textrm{ points in the staggered quadrature grid}\}, 
\end{equation}

\begin{figure}[h]
\begin{center}
\includegraphics[width=0.3\textwidth]{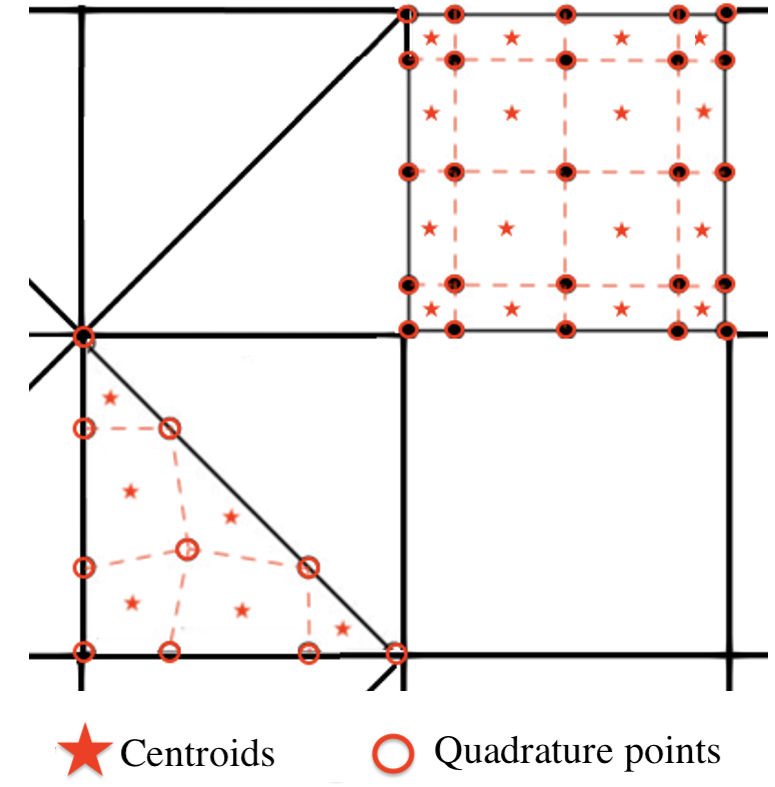}
\caption{Example lattice used for detecting structural nonconformity in a mesh.}
\label{fig:latticevisual}
\end{center}
\end{figure}

For a triangle, we can follow the same idea of combining the quadrature points and the centroids to construct the lattice. An example lattice on a mesh with quadrilaterals and triangles is shown in \Cref{fig:latticevisual}.\\

Consider a function in 2D defined by \Cref{eq:func31}. 

\begin{equation} 
\label{eq:func31}
f(x,y) = \nu \sin(2\pi x)\sin(2\pi y- 0.85\pi),
\end{equation}
where
\[   
\nu(x,y)= 
     \begin{cases}
       \text{1,} &\quad\text{if } x\geq0 \text{ and  }{x\leq0.5} \text{ and  }{y\geq0.4}    \text{ and  }{y\leq0.85} \\
        \text{0,} &\quad\text{ otherwise.}\\
        \end{cases}
\]

\begin{figure}[h]
\begin{center}
\includegraphics[width=0.9\textwidth]{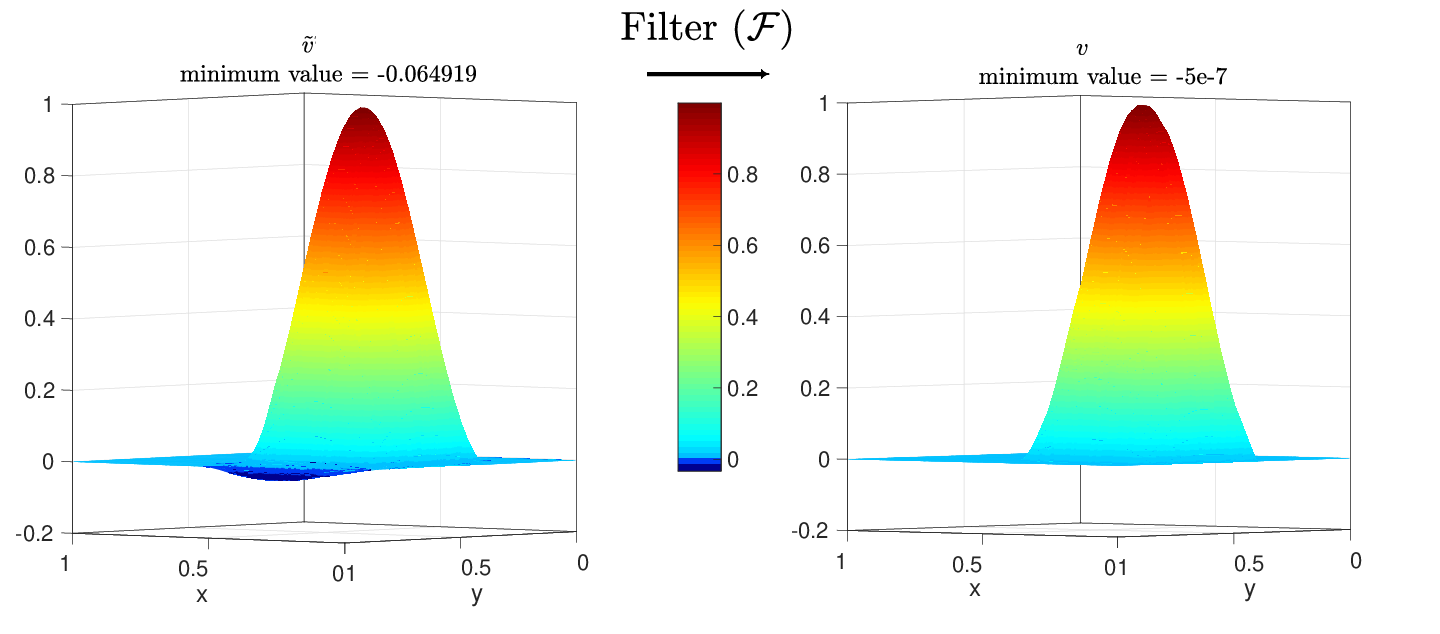}
\caption{\textit{Left:} Projection of \Cref{eq:func31} using polynomial order $= 5$. \textit{Right:} Filtered version after application of the structure-preserving filter. The filter converges to a solution that preserves positivity within a tolerance (set to $10^{-6}$ here).}
\label{fig:init2d}
\end{center}
\end{figure}

Projecting \Cref{eq:func31} using polynomial order $N = 5$ on a quadrilateral in the domain $\Omega = [0,1]\times[0,1]$, we get the coefficients $\bs{\tilde{v}}=\{\tilde{v}_j\}_{j = 0}^{P-1}$. The projection is shown in the left subfigure of \Cref{fig:init2d}.  In this case, $P= (N+1)^d$, where $d$ is dimension. Since {$d = 2$}, $P = 36$.  The original function is non-negative,  so the projection should preserve positivity. To this end, we employ the optimization outlined in \Cref{sec:algor} as a postprocessing filter. The process behaves as a spectral filter based on the empirical evidence presented in \cite[Section 5.2]{p0} and further theoretical explanation provided in \cite[Proof 5.1]{p0}. We will hereafter refer to this optimization as a nonlinear filter ($\mathcal{F}$). Let $\bs{{v}}$ be the filtered version of $\bs{\tilde{v}}$. The filtered projection is shown in the right subfigure of \Cref{fig:init2d}. \\

\subsection{Solution to advection problem in dG }
The advection equation for a quantity described by a scalar field $u$ is expressed mathematically by a continuity equation \Cref{eq:adv} in the domain $\Omega \subset {\R}^d$.
\begin{equation}\label{eq:adv}
{u}_t + \bs{a} \cdot \nabla {u} = 0,
\end{equation}
 where $\bs{a}$ represents the velocity of advection. Assume proper initial and boundary conditions are specified. \\

Galerkin-type methods assume an ansatz for $u$ as a time-varying element of a fixed {$P$}-dimensional linear subspace $V$ that contains the polynomial functions for a fixed degree {$P$}, where frequently $V \subset L^2(\Omega)$. 
\begin{align}\label{eq:u-ansatz}
  u(x,t) \approx u_{P}(x,t) &\coloneqq \sum_{i=0}^{P-1} \hat{v}_i(t) \phi_i(x), & V &= \mathrm{span} \{ \phi_0 \ldots, \phi_{P-1} \},
\end{align} 
where $x \in \Omega$, and {$\{\phi_j\}_{j=0}^{P-1}$} represents the traditional FEM basis, which is not necessarily orthogonal.{ Let $\{\psi_j\}_{j=0}^{P-1}$ represent a collection of orthonormal basis functions We can transform the vector of coefficients $\bs{\hat{v}}$ into its orthonormal form $\bs{\tilde{v}}$. }\\

Consider a partition of $\Omega$ consisting of $E$ {subdomains} defined by $\tau(\Omega) = \{e_0, e_1,\cdots,e_{E-1}\}$.  The discontinuous Galerkin formulation assumes that $V$ is comprised of functions that are polynomials of a fixed degree $P$ on each {element} on $\Omega$; where discontinuities in derivatives are allowed only at partition boundaries. The semidiscrete form for \Cref{eq:adv} is derived in the standard Galerkin way by using the ansatz \Cref{eq:u-ansatz} and forcing the residual to be $L^2$-orthogonal to $V$.  Usually, integration by parts is performed in the residual orthogonalization step, and often depending on the equation and spatial discretization, numerical flux and/or stabilization terms are included in the weak formulation.\\

The result is a system of ordinary differential equations prescribing time-evolution of the discrete degrees of freedom represented by vector {$\bs{\tilde{v}} =  \{\tilde{v}_0,\ldots,\tilde{v}_{(P-1)}\}. $}

\begin{align}\label{eq:disadv}
  \bs{M} \frac{\partial}{\partial t}{\bs{\tilde{v}}} + \bs{A}{\bs{\tilde{v}}} = \bs{F(\tilde{v})}
\end{align}
where $\bs{M} and \bs{A}$ are the $P \times P$ mass and advection matrices, respectively,  defined as
\begin{align*}
  (M_{i,j}) &= \left\langle \phi_i, \phi_j \right\rangle, & 
   (A_{i,j}) &= \left\langle \phi_i, a \frac{\partial}{\partial x} \phi_j(x) \right\rangle,
\end{align*}
and $\textbf{F}$ is a general term to compensate for any numerical fluxes or stabilization terms. Finally, \Cref{eq:disadv} is transformed into a fully discrete system that can be solved using an appropriate numerical integration scheme.\\

Consider the case when the solution defined by $\hat{v}^n$ on a particular partition $e$ at a timestep $n$ does not satisfy the desired structural properties.  To preserve the structure, we employ the optimization defined in \Cref{eq:V-opt} as a postprocessing filter. Let $\bs{v}$ represent the filtered solution that satisfies the structure. The process can be illustrated in \Cref{filter}. The proposed filter works by simply augmenting this scheme nonintrusively. 
\begin{equation}\label{filter}
  \bs{{v}}^n \xrightarrow{\textrm{Timestepper for \Cref{eq:adv}}}  \bs{\hat{v}} \xrightarrow{ (\mathcal{F})} \bs{{v}},
\end{equation} 

Since DG discretizations have degrees of freedom that are decoupled across the subdomains, the total degrees of freedom in the optimization step need to be the only ones over the subdomain. In addition, the optimization needs to be applied only in the {subdomains} that have structural violations discovered at the lattice points, which leads to a parallelizable set of independent optimizations each in {$P$ dimensions}.  We can filter the selected {subdomains} simultaneously, which further adds to the procedure's efficiency. \\


For a mesh with the number of elements $E > 1$, the process of enforcing positivity per element leads to changes in solution properties; in particular, elementwise boundary values and the mean of the discrete solution are not preserved. Assume that numerical fluxes are computed as explicit functions of element boundary values. A shift in element boundary values by the filter causes changes in the corresponding fluxes, thus adding errors to the simulation. We can resolve this issue by imposing additional equality constraints (i.e., function values at element boundaries) in the filter.  {However, for the flux conservation, we cannot make an assertion that the resulting filtered solutions should satisfy the original flux/jump conditions or some updated conditions, based on the new filtered basis coefficients. In the case of mass (integral) conservation, additional equality constraints can be imposed.  {The values of fluxes are often used to ensure the properties of numerical schemes therefore by preserving fluxes one can attest that the properties of the original scheme have been retained.} Preserving fluxes in this context corresponds to corrections in one element need not have an impact on neighboring elements, and one way to promote isolation of the optimization effect to the element under consideration is to ensure that interelement communication through fluxes is undisturbed.  The formulation for flux conservation in \cite[Section 3.2 and 3.3]{p1} is unique to one-dimensional scenarios. {Generalization to higher dimensions is not immediate. Additional analysis is needed to incorporate different elements and, therefore, beyond the scope of this paper. } The robust design of the filter allows for the incorporation of an arbitrary number of structural constraints into one optimization problem. However, increasing the number of constraints reduces the degrees of freedom for the filter accordingly. In particular, if the filter has $P$ degrees of freedom and we have {$q$} linear equality constraints, then the dimension of the optimization problem can be {reduced to {$(P - q)$}} using the procedure
described in \cite[Section 3.2]{p1}.

\section{Formulation of structure-preservation problem and solution design}\label{sec:formalism}

This section establishes the notations and summarizes formulation, filter design, and implementation ideas from \cite{p0,p1}.

\subsection{The problem}

Let $u(x) \in V$ be a function,  where $V$ is a finite-dimensional {subspace of a Hilbert space} of real-valued functions on $\Omega \in \R^d$.  {Let V contain a collection of orthonormal basis functions $\{\psi_n\}_{n=1}^N$ for some {$N \in \N$}.} \\

{Therefore,  we have}
{{
\begin{align*}
  V &= \mathrm{span}\left\{ \psi_1 \ldots, \psi_{N} \right\}, & \left\langle \psi_i, \psi_j \right\rangle &= \delta_{ij}, & i, j &= 1 \ldots, N,
\end{align*}
}}

{where $\left\langle\cdot,\cdot\right\rangle$ is the inner product on $V$, and {$\delta_{ij}$}, the Kronecker delta function.}\\

{The value of $d$ is set to 1 for simplicity. The formulation follows the same steps for higher dimensions. } A comprehensive constraint-satisfaction problem can be constructed as detailed in this section. \\

Let $K$ be the collection of $u$ that satisfies the constraints. Therefore,  $K$ is the feasible subset of $V$.  An example of such a family of linear constraints is \{positivity,  boundedness\}}. {Consider a family of linear constraints such as the one in \Cref{eq:constraints}, where ${L}$ is a linear operator bounded on $V$,  and $\ell$ is a function on the domain $\Omega$. }

\begin{align}\label{eq:constraints}
  {L}(u) &\leq \ell(x), & x &\in \Omega.
\end{align}

{Note that the framework in \cite{p0} allows a finite number of such constraints to be considered simultaneously.}\\

{According to the Riesz representation theorem, if $V^*$ is a dual of $V$,  then a functional $L \in V^*$ can be associated with a unique V-representor $\ell \in V$ satisfying}

{
\begin{align*}
L(u) &= \left\langle u, \ell \right\rangle, & \forall u \in V.
\end{align*}}

{Furthermore, this $L \leftrightarrow \ell$ identification is an isometry. We will use these facts in what follows. Given $L$ that identifies $\ell$, we consider the coordinates $\widehat{\ell}_j$ of $\ell$ in a $V$-orthonormal basis,}

{
\begin{align*}
  \ell(x) &= \sum_{j=1}^N \widehat{\ell}_j \psi_j(x), & \widehat{\ell}_j &= \left\langle \ell, \psi_j \right\rangle = L(\psi_j).
\end{align*}}

{Then we have the following relations:}
{
\begin{align*}
 \left\| L \right\|_{V^\ast} = \left\| \ell\right\|_V &= \| \bs{\widehat{\ell}} \|_2, & \bs{\widehat{\ell}} &= \left( \widehat{\ell}_1, \; \widehat{\ell}_2, \; \ldots, \; \widehat{\ell_N} \right)^T,
\end{align*}}

{where $\|\cdot\|_2$ is the Euclidean norm on vectors in $\R^N$. }\\

The problem is essentially that of developing a map $\mathcal{F}$ from a function $u \in V$,  to a unique function $\mathcal{F}(u) \in K$.  Initially, {$u$ may} not satisfy the structural constraints but its mapped version $\mathcal{F}(u)$ does.  In \cite{p0} we introduce a strategy to solve the following well-posed convex feasibility problem.

\begin{align}\label{eq:V-opt}
  \mathcal{F}(u) \coloneqq \argmin_{f \in K} \| u - f\|_V,
\end{align}

Under certain conditions, such as if $0 \in K$,  $\mathcal{F}$ is norm-contractive{; therefore, the operation can be called a filter}. For a brief discussion of the {norm-contractive} nature of \Cref{eq:V-opt}, see \cite[Proposition 5.1]{p0}.

\subsection{Toward construction of a convex optimization-based solution to \Cref{eq:V-opt}} \label{sec:towards}


We first transform the continuous problem \Cref{eq:V-opt} to a discrete version to construct a feasible solution. Let  $\{\psi_j\}_{j=0}^{P-1}$ be a collection of orthonormal basis functions in $V${, different from the standard FEM basis ($\phi$), which are not orthonormal}. \\


Consider $C$ as the affine conic region in $\R^P$ corresponding to $K\subset V$. From previously established notations, $\bs{\tilde{v}}$ is a vector of expansion coefficients of $u$ that does not satisfy that desired constraints, and $\bs{v}$ is a vector of filtered coefficients of $u$ that satisfies the desired constraints.  Representing the Euclidean 2-norm on vectors as $\|\cdot\|_2$,  and the fact that $C \subset V$,  we get the equivalent of the optimization problem \Cref{eq:V-opt} as follows:

\begin{align}\label{eq:v-optimization}
  \argmin_{\bs{{v}} \in C} \| \tilde{\bs{v}} - \bs{{v}}\|_2.
\end{align}

\begin{figure}[htbp]
  \begin{center}
        \includegraphics[width=0.49\textwidth]{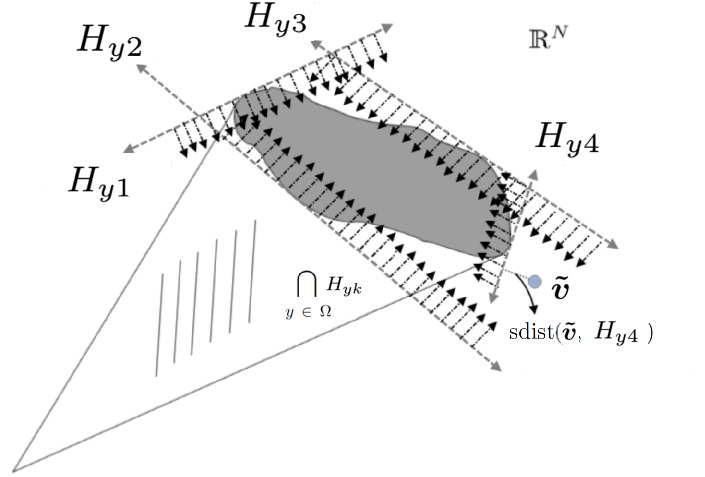}
      \includegraphics[width=0.49\textwidth]{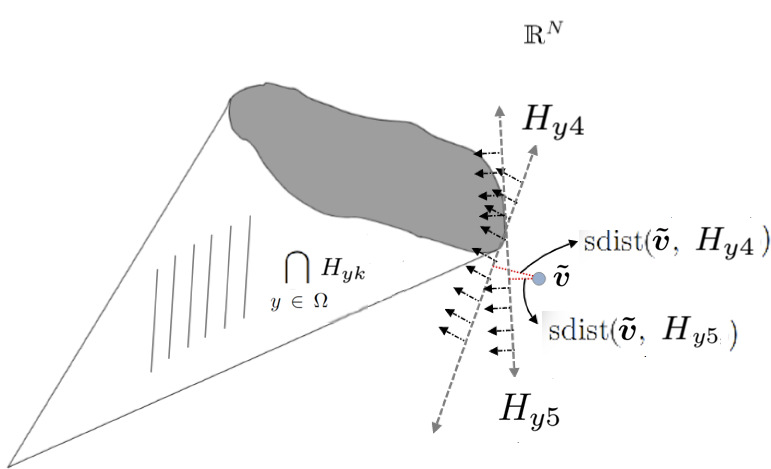}
  \end{center}
 \caption{Division of the space into half-spaces by the hyperplanes representing constraint boundaries defined by the set of points (say $y1, y2\cdots yn$) for the initial projection with coefficients $\bs{\tilde{v}}$. The algorithm greedily calculates the correction $s(x)$. Left: A geometrical visual of the distance calculation from $\bs{\tilde{v}}$ to the hyperplanes defining the boundaries of the constraints. {The hatched area inside the cone represents the feasible region}. Right: Projection of $\bs{\tilde{v}}$ on to $H_{y4}$. $H_{y4}$ is selected over $H_{y5}$ since it defines the violating hyperplane that is farthest away.}\label{fig:distandprojection}
 \end{figure}

{For a fixed $x \in \Omega$, 	let $H_x$ be a $(P-1)$-dimensional planar surface defined by an equality constraint corresponding to \Cref{eq:constraints} for a fixed $x$. Therefore, {$L(u) = \ell(x)$}.  Writing $u$ as an expansion in degrees of freedom $\tilde{v}$ that corresponds to a single linear equality constraint for the $\tilde{v}$, i.e., a $(P-1)$-dimensional plane. In a geometric sense, $H_x$ is a surface representing the constraint boundary in the space $\R^P$, dividing it into two hyperspaces as shown in \Cref{fig:distandprojection}.  One hyperspace represents the region that satisfies a linear inequality constraint and the other that does not. For a single linear constraint, $C$ can be defined in terms of hyperspaces as \Cref{eq:constraints}}. \\
{
\begin{align*}
  C = \bigcap_{x \in \Omega} \{ \textrm{feasible halfspace defined by plane } H_x\},
\end{align*}
}

{To obtain the full feasible set $C$, the filter iteratively  projects $\tilde{\bs{v}}$ onto the intersection of feasible halfspaces defined by individual hyperplanes $\{H_x\}$.  To describe the feasible regions further, let us consider a geometric interpretation.  Let the current state vector of the projection defined by $\tilde{\bs{v}}$ be a point shown by the blue dot in \Cref{fig:distandprojection}. The filter works by applying correction \Cref{eq:greedy-update} to $\tilde{\bs{v}}$ by computationally inspecting the signed distance function \Cref{eq:distform}}.

{
\begin{align}\label{eq:distform}
  s(x) \coloneqq \mathrm{sdist}(\bs{\tilde{v}}, C_x) = \left\{ \begin{array}{cc} -\mathrm{dist}(\bs{\tilde{v}}, H_x), & x \not \in C, \\
  +\mathrm{dist}(\bs{\tilde{v}}, H_x), & x \in C.
  \end{array}\right\}.
\end{align}
}

For the positivity-preservation example,  the process of ``inspection" means the ability to compute the global minimum of $s(x)$ to determine a region where $s$ is negative. Based on this inspection, the algorithm projects the state vector of the current iterate $\bs{\tilde{v}}$ onto $H_y$ for some $y \in \Omega$. This projection can cause a particular $H_{y1}$ with a positive $s(y1)$ to go negative. Therefore, we need to repeat the projection step until the solution $\bs{\tilde{v}}$ lies completely in (or on the boundary of) $C$. \cite[Algorithm 1]{p0} summarizes all the steps of the filter.\\ 

{
At each iteration,  the spatial point $x$ is found that minimizes the signed distance function given by \Cref{eq:global-minimization}.  Note that for more than one constraint we need to keep track of $x$ and the constraint specific to the minimum signed distance function.}

\begin{align}\label{eq:global-minimization}
 (x^\ast) \coloneqq \argmin_{x \in \Omega} s(x) 
\end{align}

{This greedy procedure is an infinite-constraint generalization of Motzkin's relaxation method \cite{motzkin_relaxation_1954}.
  The major computational work in the filtering procedure is to minimize the objective defined by \Cref{eq:sdist2} at each iteration.
\begin{align}\label{eq:sdist2}
  s(x) &= \lambda(x) \left( {L}_x(u) - \ell(x) \right), & \lambda^2(x) \coloneqq \frac{1}{\sum_{j=0}^{P-1} \left({L}_x(\psi_j)\right)^2}.
\end{align}
}


{Next, update $\bs{\tilde{v}}$ by projecting it onto the hyperplane corresponding to $x^\ast$:
\begin{align}\label{eq:greedy-update}
 \bs{\tilde{v}} \gets \bs{\tilde{v}} + \bs{h}(x^\ast) \min\left\{0, s(x^\ast)\right\}.
\end{align}
where $\bs{h}(x,)$ is the normal vector corresponding to hyperplane $H_x$ of a single constraint family that points toward $C$. This vector is readily computable from the orthonormal basis, see \cite{p0} for details. This procedure is repeated until $s(x^\ast)$ vanishes to within a numerical tolerance. }\\

Finding the violating hyperplane that is farthest away  corresponds to finding the global minimum of $s(x)$ on $\Omega$. Unlike the 1D procedure described in \cite{p0},  multidimensional global minimization requires a more expensive approach such as gradient descent (GD). {In this paper, we implement 2D and 3D global minimization using line-search-based GD with backtracking. We further provide a strategy to investigate the times taken for convergence and the accuracy of this GD approach by varying the values of the parameters. Then, we choose the best values for the parameters for the  accuracy of the method.} \Cref{sec:rootfinder} provides details of the GD variant chosen for filter implementation and corresponding parameter selection to efficiently solve the minimization problem.
\section{Algorithm}\label{sec:algor}
The structure of the solution is preserved by applying the filter as a postprocessing step that corrects $\bs{\hat{v}}$ obtained from the timestepping algorithm using iterative optimization. Each iteration of the filter involves the following crucial computational processes:

\begin{enumerate}

\item \label{it:it1}{A GD-based minimization to solve \Cref{eq:global-minimization}, which is a part of the global search for structural violations.} 
\item {A calculation of the correction to $\bs{\hat{v}}$ using \Cref{eq:greedy-update}}
\end{enumerate}
We briefly discuss the algorithmic details and the choices made to streamline the processes and their applications to the PDE timestepping.

\subsection{Finding the global minimum }\label{sec:rootfinder}

An important and the most expensive part of the optimization process is {the global minimization \Cref{eq:global-minimization}. For efficiency and usability of the filter, we need to fine-tune the most expensive step,  i.e., minimizing the scaled distance function \Cref{eq:sdist2}, which quantifies the difference between the current state of the coefficients and the coefficients that satisfy the constraints at a critical point of the current iteration.} The success of the proposed filter depends on accurately finding the points in the domain that violate the structure. {Since the problem formulation uses dG, it can be solved element by element on a domain with $E$ elements. Applying the filter at each timestep is essentially solving $E$ optimization problems at every timestep. Each optimization problem is over a single element. We can reduce the number of optimization problems by testing the solution for structural inconsistencies at a lattice of points defined by \Cref{eq:lattice}. For example, in the case of positivity preservation, if the minimum value of an element is nonpositive, the filter is applied to that element. } \\

In the 1D case, the filter finds the minimum value by evaluating the solution at the critical points on the element. If the subspace $V$ restricted to an element is a space of polynomials, then the critical points of \Cref{eq:sdist2} are defined by the roots of the derivative of \Cref{eq:sdist2}. For the root-finding problem, we employ the confederate matrix-generation approach. Relative to 1D,  minimization in higher dimensions is more complex and expensive.  Many variants of the gradient descent-based (GD-based) methods have been proposed to perform this operation efficiently and accurately.  These methods operate only up to a certain tolerance and, like all the GD-based methods, have the potential for stagnation.   \\

Since there is no multidimensional version of the confederate linearization approach that exactly solves the first-order optimality conditions, the problem of finding the critical points on the multidimensional domain is nonconvex.  Therefore, we cannot ensure that all such points will be flagged for correction by any particular minimization algorithm. From our analysis, one suitable GD approach for this job is to use adaptive descent size and the ability to retrace the steps.  For this reason, we choose the method of steepest descent with a backtracking line search strategy. {The backtracking line search starts with a relatively large step size and repeatedly shrinks it by a factor $\gamma$ until the Armijo–Goldstein condition \Cref{eq:armijogold} is fulfilled.} This widely used algorithm makes an informed choice about step size and direction at each iteration to efficiently arrive at the optimum value.  It avoids stagnation by tracing back its steps if the difference between descent values falls below a fixed tolerance. An essential part of this algorithm is choosing the parameters denoted by $ c \in (0,1)$ and $\gamma \in (0,1)$ that determine the step size and backtracking performed by the algorithm. For an objective function $f$, with a starting position $x$,  given parameters $c$ and $\gamma$ the Armijo–Goldstein condition is defined as \Cref{eq:armijogold}.

\begin{equation}\label{eq:armijogold}
     {\displaystyle f( {x} +\gamma_j {p} )\leq f( {x} )+\gamma_{j} \,c\,\nabla f(x)^T p\,} \hspace{2cm} \forall \textrm { } j \textrm{ until convergence,}
\end{equation}
where $\nabla f(x)^T$ is the slope in the direction of descent $p$ and $\gamma$ is calculated as
\[   
\gamma_j= 
     \begin{cases}
       {\gamma,} &\quad\text{if } j = 0\\
       \, c\gamma_{j-1} &\quad\text{ otherwise}.\\
        \end{cases}
\]

{$c$ and $\gamma$ play a vital role in testing the condition \cite{armijogoldstein}, which determines whether a step-wise movement from a current position to a modified position achieves an adequately corresponding decrease in the objective function. } Although many recommendations exist \cite{armijogoldstein,crockett1955gradient,curry1944method,bertsekas1997nonlinear,boyd2004convex} for an optimal way to choose $c$ and $\gamma$, these parameters must be tuned for maximum efficiency.\\

\subsection{Applying the structure-preserving filter to $\bs{\tilde{v}}$ obtained from the timestepper}

From the discussion in \Cref{sec:method}, the optimization that enforces a constrained solution is applied as a postprocessing part inside the timestepper,  independent of the choice of the timestepping routine.  As predicted by the `curse of dimensionality',  structure preservation and optimization for higher dimensional problems are more complex than their 1D counterparts. Since the central idea of the filter is agnostic to the dimensions of the problem, the \Cref{alg:greedy} remains the same as described in \cite{p0}. Similarly, the applications to PDE solutions follow the procedures established in \cite{p1} for 1D problems.  \\


\begin{algorithm}[H]
 \caption{Constrained PDE timestepping}
\label{alg:greedy}
\begin{algorithmic}[1]
\STATE{Input: Terminal time $T$, timestep size $\Delta t$, PDE solver spatial basis ${\phi}$}
\STATE {Define nsteps $= \frac{T}{\Delta t}$}
\FOR{$i= 0,\cdots,$ nsteps }
\STATE{Solve PDE and obtain the orthogonalized coefficients $\bs{\tilde{v_e}}^{i} \in \R^P$ for all elements in $e$ in $\Omega$}
	\FOR{Each element $e \in \Omega$} 
  		\WHILE{True}
     		\IF{ { $\exists (x \in \Omega)$ such that $s(x) < 0$}}
        		\STATE{{Compute $(x^\ast)$ as defined in \Cref{eq:global-minimization}} via GD line search  from \Cref{sec:rootfinder}}
    		\ELSE{}
        		\STATE{{break}}
    		\ENDIF
     \STATE \label{alg:greedy:update2}{Update $\bs{\tilde{v_e}}^{i+1}(t)$ via \Cref{eq:greedy-update}.}
    \ENDWHILE
 \STATE{Append to global $\bs{{v}}^{i+1} = [\bs{{v}}^{i+1} , \bs{\tilde{v_e}}^{i+1}]$}
 \ENDFOR
\ENDFOR
\end{algorithmic}
\end{algorithm}
Conversion from an input coefficient vector $\hat{\bs{v}}$ to corresponding coefficients $\tilde{\bs{v}}$ in an orthonormal basis is the first stage of the filtering procedure. \Cref{alg:greedy} presents a summary of steps taken by the PDE solver for the filter application. \\

A significant chunk of work in the iterative correction stage is done by evaluating the solution at the lattice points using basis interpolation, an expensive operation that involves processing and storage of large interpolation matrices. To make the processing computationally efficient, we use the Barycentric interpolation method from  \cite{Trefethen04},  expanded upon by \cite{laughton2021fast}. This method reduces the number of calculations and storage required to evaluate the basis functions at arbitrary points in the domain, thereby reducing the cost of the iterative correction stage. Although the implementation details of Barycentric interpolation method are beyond the scope of this paper, for the minimization part of the numerical experiments in \Cref{sec:results},  we observe a significant speed-up using the Barycentric approach.

\section{Main results}
\label{sec:results}

This section presents four-part numerical results of the filter's application to different types of problems in 2D and 3D. {Since applying multiple constraints together can be reduced to one optimization problem on a smaller solution space, the numerical results presented here focus on positivity preservation.  For 1D examples that preserve multiple structural constraints simultaneously, refer to the results presented in \cite{p1}.}\\

\Cref{sec:res1} details the process to choose the parameter values for the GD line search algorithm described in \Cref{sec:rootfinder}.  \Cref{sec:res2} shows the results of applying the filter to a single-element function projection in multiple dimensions.   \Cref{sec:advsolution} presents the results and analysis of the application of the filter using the dG formulation of the advection problems defined on composite and homogeneous domains, respectively. Finally, in \Cref{sec:results3} we present the results of the FEM solution to the mathematical-biology problem of platelet aggregation and blood coagulation (PAC) with and without the application of the filter. The PAC model is posed as an advection-diffusion-reaction (ADR) system of PDEs.

\subsection{Parameter selection for {backtracking line search used in GD-based minimization} }\label{sec:res1}

{To fix the values for $c$ and $\gamma$ \Cref{eq:armijogold} for the numerical results in \Cref{sec:results}, we first choose a set of functions. The choice of the functions depends on factors such as the constraint type and the domain geometry. For example, for an application where positivity is the primary structural concern, a list of functions with values close to zero is a good choice. For applications involving discontinuities in the function or the domain, discontinuous functions would be a good choice. The feasible space for $c$ and $\gamma$ is $(0,1)$. A non-comprehensive sampling on $(0,1)$ is performed for both parameters, and the GD is called with the parameters set to permutations of the sampled values. GD\_linesearch denotes the call to the GD with line search. As its output, the GD\_linesearch routine provides the number of iterations (denoted by $niter$) taken by the algorithm. The error refers to the difference between the minimum value found and the known minimum value.\\}

{As a large number of GD iterations per filter iteration makes the filtering process more expensive, choosing the parameter values for which GD converges in the least number of iterations is desirable. Another important quantity to consider here is the error. It is desirable to have GD return with a minimum value as close to the known (or golden) minimum value as possible. With these quantities ($niter$ and $err$) as the metrics, we prescribe a procedure to select the ideal value of $c$ and $\gamma$ in \Cref{alg:candgam} [Appendix \ref{sec:apdx}]. The selection criteria narrow the samples to a range of $c$ and $\gamma$ with the minimal number of iterations and then pick from that subset the ones with the smallest error. In the case of a tie, all the $c$ and $\gamma$ values are included in the return variables. The procedure can be summarized as taking a list of functions that fairly represent the complexity of space and returning two ranges of ideal values: one for each of the parameters $c$ and $\gamma$, respectively. }{Since the quadrilateral and hexahedron represent richer space than their other canonical 2D and 3D counterparts, \Cref{alg:candgam} [Appendix \ref{sec:apdx}] is applied and analyzed only on a quadrilateral for 2D and a hexahedron for 3D}.\\

Future developments in the approaches for finding optima in multiple dimensions may improve the efficacy and efficiency of the filter.  To this end and to provide robustness,  we employ the object-oriented principle by modularizing the minimization part of the filter.\\

We consider different example functions \Cref{eq:fns} on the 2D domain $[-1,1] \times [-1,1]$  and \Cref{eq:fns3} on the 3D domain $[-1,1]\times[-1,1]\times[-1,1]$, which have varying complexities and projection orders to set the GD parameters comprehensively and fairly.

 \begin{align}\label{eq:fns}
 f_0(x,y) &= (x + 0.6)^2 + (y - 0.2)^2. \nonumber \\
 f_1(x,y) &= -\sin((x-0.1)+0.5\pi)\cos(y-0.2). \nonumber \\
 f_2(x,y) &=
    \begin{cases}
       \text{1,} &\quad\text{if } x\leq0 \text{ and  }{y\leq0}, \\
        \text{0,} &\quad\text{ otherwise. }  \\
        \end{cases} 
\end{align}
 
 \begin{align}\label{eq:fns3}
 f_3(x,y,z) &= (x + 0.6)^2 + (y - 0.2)^2 + (z+0.1)^2.  \nonumber \\
 f_4(x,y,z) &= -\sin((x-0.1)+0.5\pi)\cos(y-0.2)\cos(z-0.2). \nonumber \\
 f_5(x,y,z) &= 
     \begin{cases}
       \text{1,} &\quad\text{if } x\leq0 \text{ and  }{y\leq0} \text{ and  }{z\leq0}, \\
        \text{0,} &\quad\text{ otherwise.}\\
        \end{cases} 
 \end{align}
 
Using the orthogonal basis for \Cref{eq:greedy-update},  $\bs{\tilde{{v}}}$ is obtained, which is then used by \Cref{alg:candgam} to find the values of parameters $c$ and $\gamma$. For functions that have an exact analytical minimum, calculating the error between the exact minimum and the optimum returned by the GD is straightforward. For the discontinuous functions such as $f_2$ and $f_5$, an approximate \textit{golden} minimum is calculated by projecting the function using polynomial order $N = 8$ on a dense grid of points in $\Omega$ and finding the least value of the projected polynomial.  In all our test cases,  GD converges to the exact (or golden) minimum within an acceptable tolerance (numerical 0 set to $10^{-7}$). Therefore,  the selection of the parameters is made based on the fewest iterations taken by GD ($niter$) to converge.  \\
\begin{table}
\caption{Results for the 2D experiment on a quadrilateral domain to find the optimal ranges of the gradient descent line search parameters: $c$ and $\gamma$.  M denotes the polynomial order for projected functions \Cref{eq:fns}. The number of iterations taken by GD is denoted by $niter$. The number of quadrature points for projection is constant ($Q = 11$)}\label{tab:f1res}
  \includegraphics[width=\textwidth]{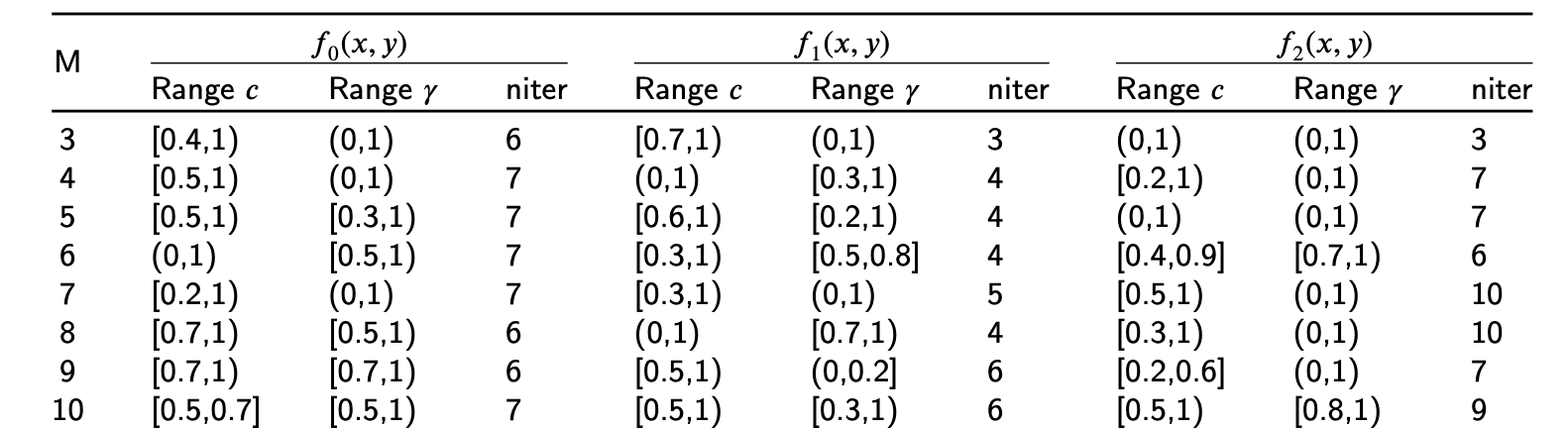}
\end{table}
%

Consider  $k$ discrete equispaced samples of $\in (0,1)$. For $k = 9$, by applying \Cref{alg:candgam}, we get the results shown in \Cref{tab:f1res}.\\

Based on the results of these tests,  we infer that for the numerical experiments, the appropriate choices for $c$ and $\gamma$ in 2D are 0.7 and 0.7, respectively.  Following a similar procedure, \Cref{alg:candgam} on a hexahedron for functions defined in \Cref{eq:fns3},  the choices of $c$ and $\gamma$ for 3D are 0.2 and 0.7, respectively.  These values are used for all the numerical experiments in the numerical results. \\

The C++ implementation of the filter in Nektar++ \cite{cantwell2015nektar++} supports change in $c$ and $\gamma$ as parameters in the configuration file.

\subsection{Projection examples on 2D and 3D elements and application of the structure-preserving filter}\label{sec:res2}

\begin{figure}[h]
\begin{center}
\includegraphics[width=0.35\textwidth]{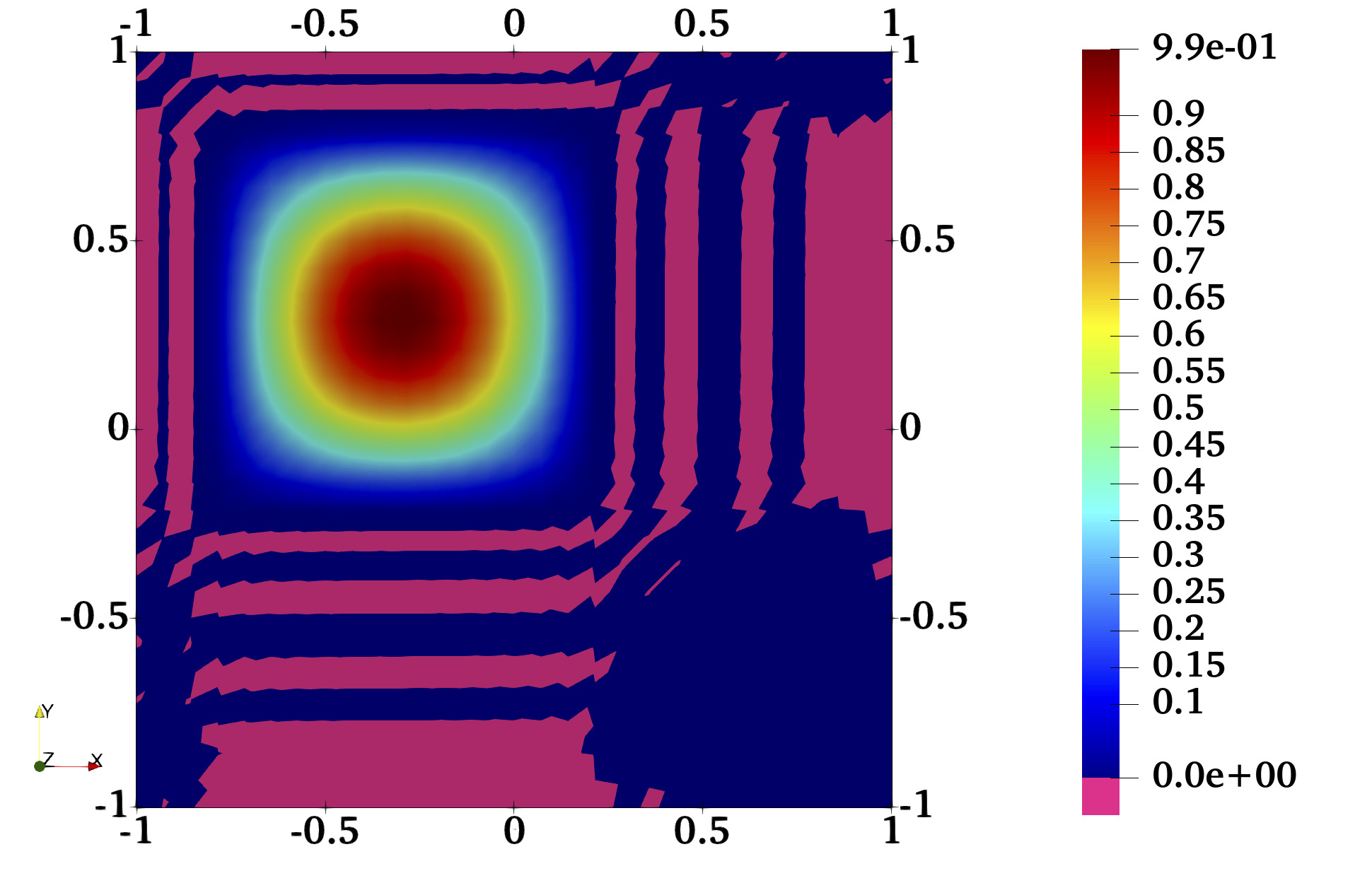}
\includegraphics[width=0.25\textwidth]{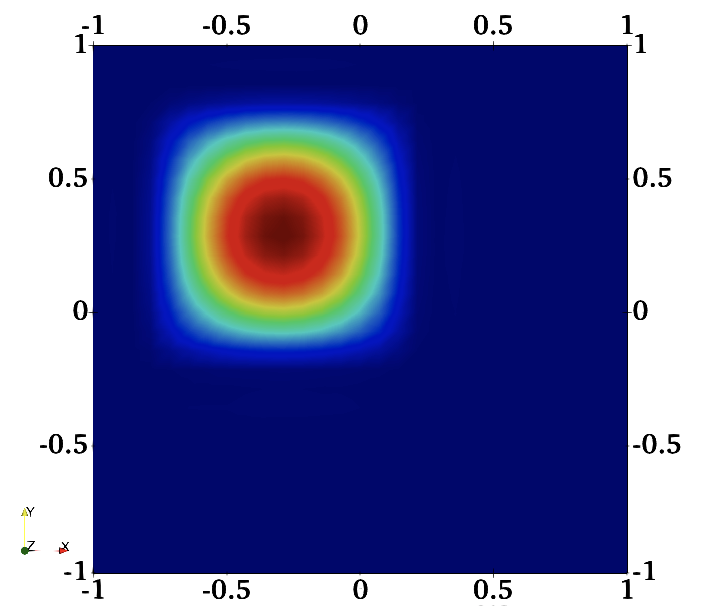}\\
\includegraphics[width=0.45\textwidth]{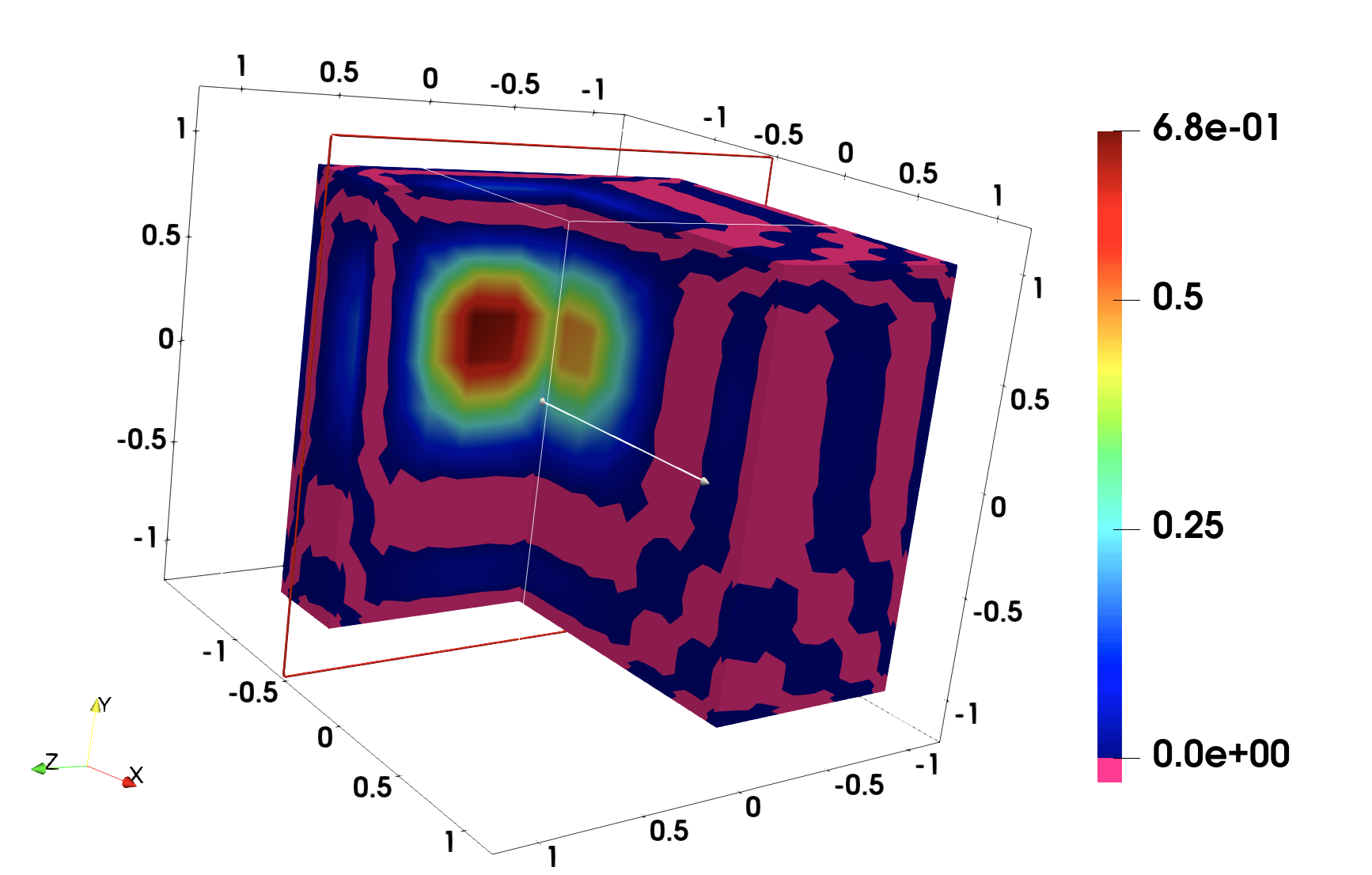}
\includegraphics[width=0.29\textwidth]{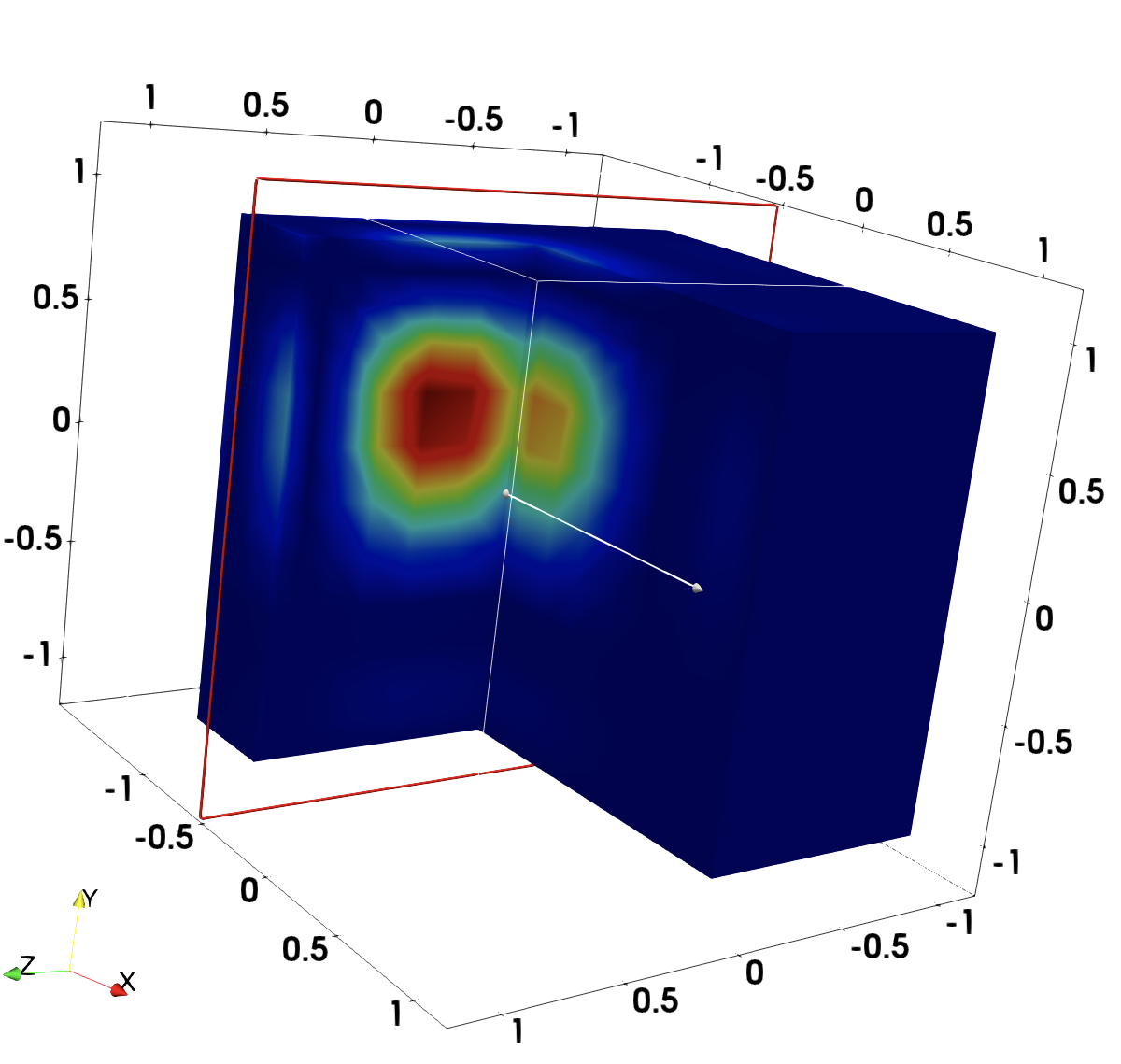}
\end{center}
\caption{Galerkin projection of $f(x,y)$ \Cref{eq:projfunction} and $f(x,y,z)$ \Cref{eq:func3} on a quadrilateral and hexahedron element, respectively. The 2D projection uses polynomial order $ = 7$, and the 3D projection uses polynomial order $ = 5$.  \textit{Left:} Unfiltered version showing areas where the desired structure (positivity) is lost. \textit{Right:} After the filter is applied, the solution is non-negative up to tolerance. }
 \label{fig:projexpt}
\end{figure}

Consider a 2D function \Cref{eq:projfunction} and a 3D function  \Cref{eq:func3} that are both discontinuous clamped versions of a smooth sinusoidal function. The initial projections of $f(x,y)$ on a quadrilateral and $f(x,y,z)$ on a hexahedron are shown in \Cref{fig:projexpt}. The discontinuity produces oscillations similar to Gibb's phenomenon upon projection, which is interesting for the application and analysis of the filter.

\begin{equation} 
\label{eq:projfunction}
f(x,y) = \nu(x,y) \sin\Big(\pi(0.2-x)\Big)\sin\Big(\pi(y+0.2)\Big),
\end{equation}
where
\[   
\nu(x,y)= 
     \begin{cases}
       \text{1,} &\quad\text{if } x\in [-0.8 ,0.2] \text{ and  }  y\in [-0.2,0.8],\\
        \text{0,} &\quad\text{ otherwise.}\\
        \end{cases}
 \]

\begin{equation} 
\label{eq:func3}
f(x,y,z) = \nu(x,y,z) \sin\Big(\pi(0.2-x)\Big)\sin\Big(\pi(y+0.2)\Big)\sin\Big(\pi(z+0.2)\Big),
\end{equation}
where
\[   
\nu(x,y,z)= 
     \begin{cases}
       \text{1,} &\quad\text{if } x \in [-0.8, 0.2] \text{ and  }y\in [-0.2, 0.8] \text{ and  } z\in [-0.2, 0.8]\\
        \text{0,} &\quad\text{ otherwise.}\\
        \end{cases}
\]

\begin{figure}[h]
\begin{center}
\includegraphics[width=0.48\textwidth]{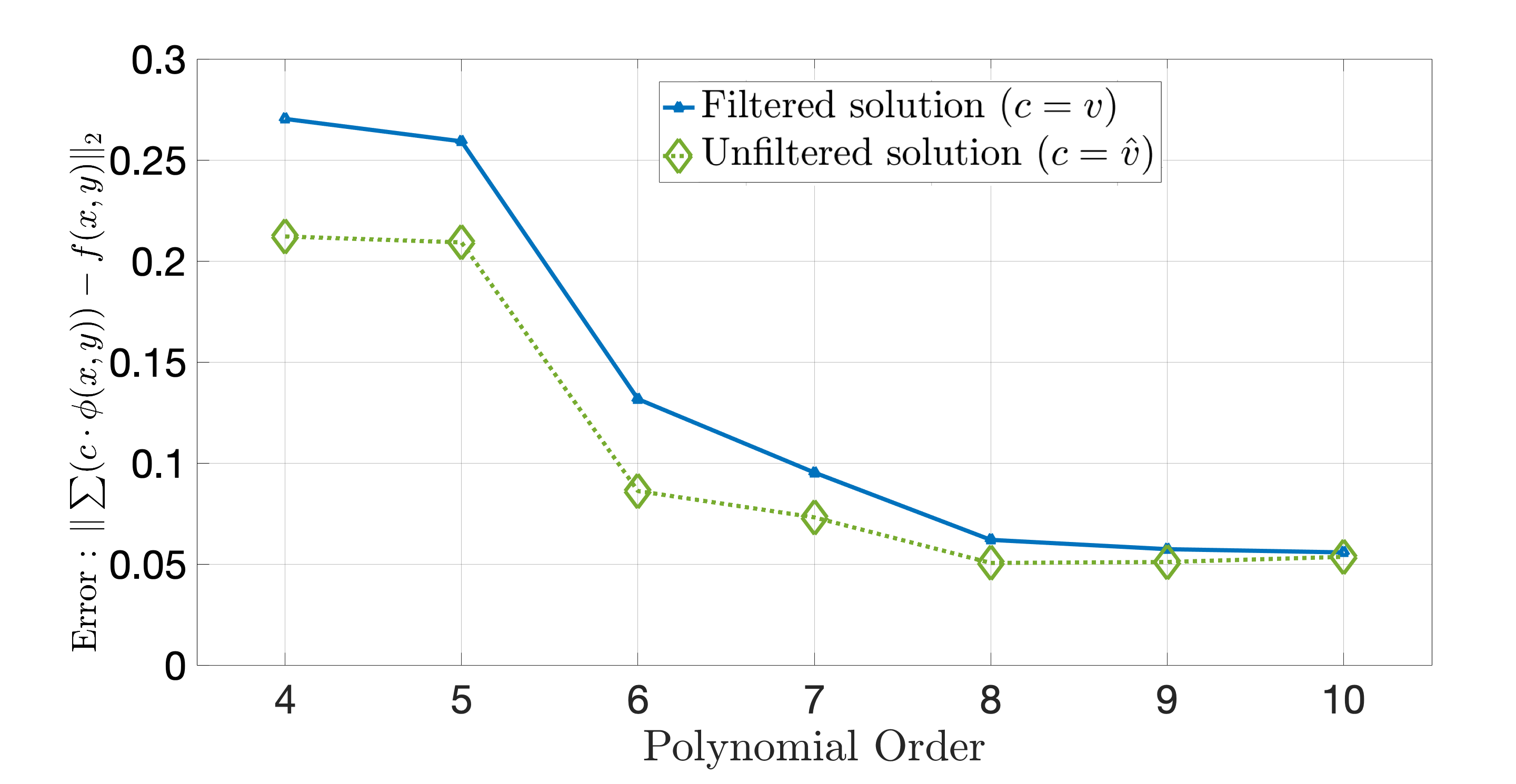} 
\includegraphics[width=0.48\textwidth]{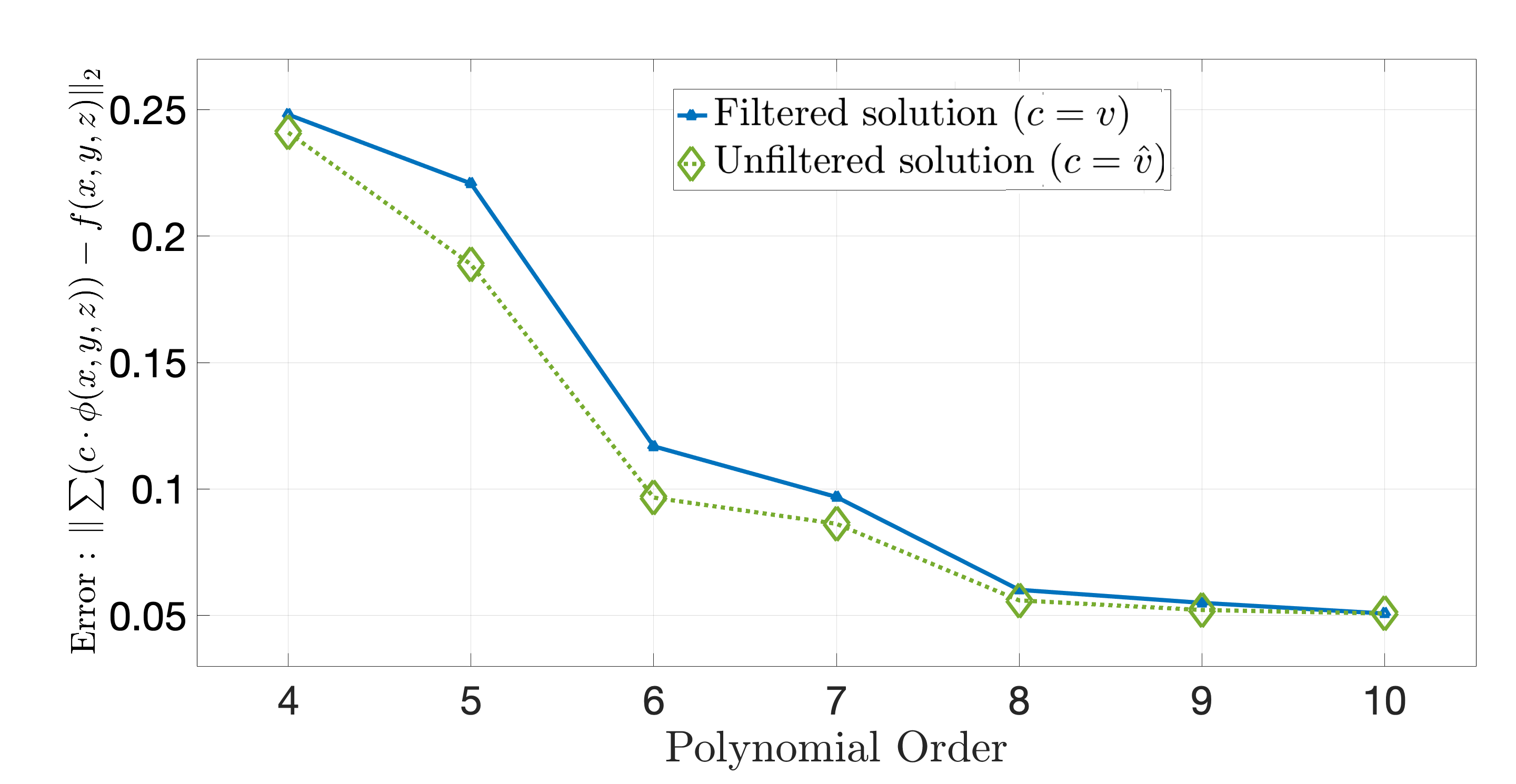}
\end{center}
\caption[A p-convergence study and cost (in s) of the filtered and unfiltered Galerkin projections from \Cref{fig:projexpt}  ]{A p-convergence study of the filtered and unfiltered Galerkin projections from \Cref{fig:projexpt}.
\textit{Left:} $f(x,y)$ defined in \Cref{eq:projfunction}  and \textit{Right:}$  f(x,y,z)$ defined in \Cref{eq:func3}}
\label{fig:propconv2d}
\end{figure}

{It is evident from the left subfigures in \Cref{fig:projexpt} that the projection of discontinuous non-negative functions leads to the negative values as shown in the highlighted regions. As seen in the right subfigures of \Cref{fig:projexpt}, the application of the structure-preserving filter restores the non-negative structures of \Cref{eq:projfunction} and \Cref{eq:func3}, respectively.} To analyze the p-convergence, the experiment is repeated for different polynomial orders, as shown in \Cref{fig:propconv2d}.\\

{\Cref{fig:propconv2d} presents $L_2$ error comparison between the filtered and unfiltered version of projections of nonsmooth functions in 2D and 3D. The choice of a nonsmooth function is to demonstrate the effect of filter application in worst-case scenarios.}

\subsection{Structure preservation in 2D and 3D advection problems} \label{sec:advsolution}

\begin{figure}[t]
\begin{center}
\includegraphics[width=0.5\textwidth]{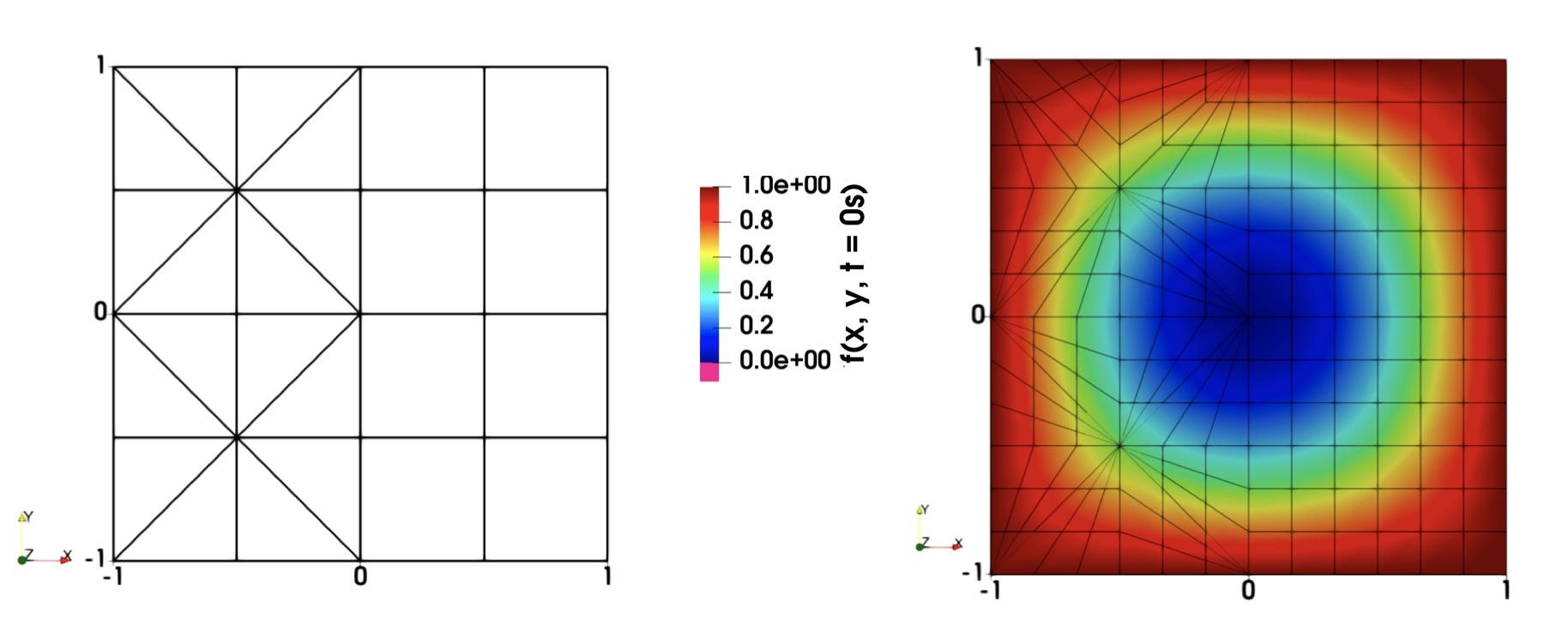}\\
\includegraphics[width=0.5\textwidth]{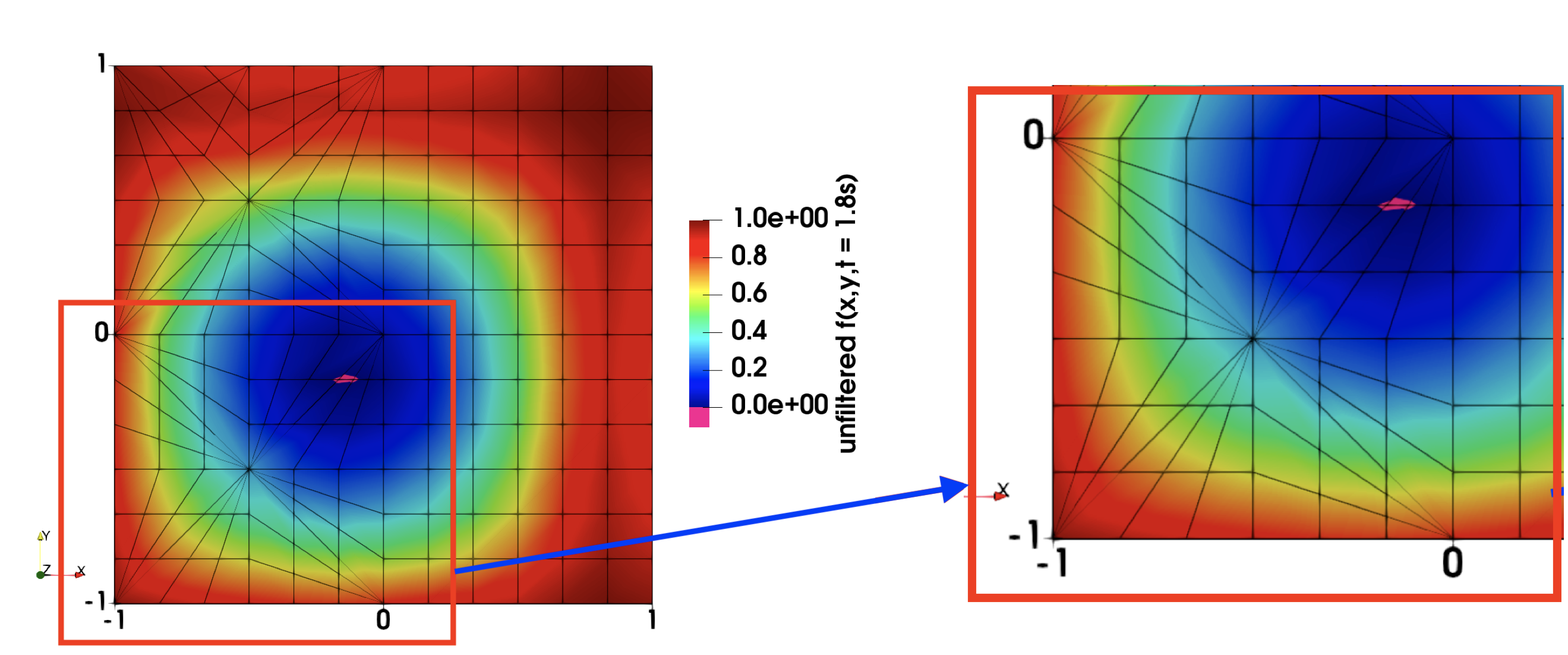}\\
\includegraphics[width=0.5\textwidth]{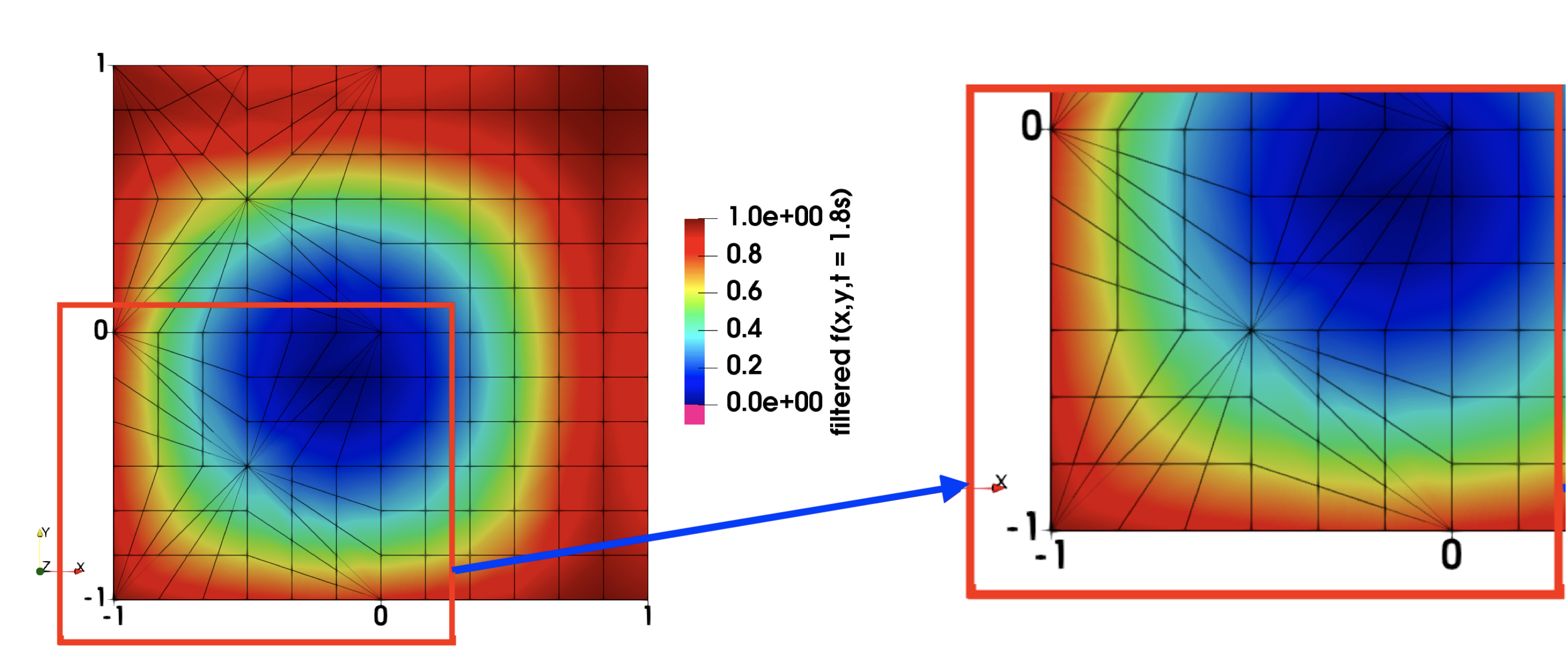}
\caption{\textit{Row 1, Left:} A 2D composite mesh used for solving \Cref{eq:adv2d}. The mesh contains 16 triangles and 8 quadrilaterals. \textit{Row 1, Right:} Initial state of the system, which is a projection of \Cref{eq:para} on the mesh.
\textit{Row 2}: Unfiltered solution at timestep 1800.
\textit{Row 3}: Filtered solution at timestep 1800.}	
 \label{fig:para2dadv}
\end{center}
\end{figure}

Consider a dG solution to the advection problem \Cref{eq:adv2d} on a 2D domain as shown in the first row of \Cref{fig:para2dadv}.
\begin{equation}
\label{eq:adv2d}
{u}_t + \textbf{a} \cdot \nabla{u} = 0,
\end{equation}

Consider a function $f$ defined on $[-1,1]\times[-1,1]$ as \Cref{eq:para}.
\begin{equation}\label{eq:para}
f(x,y) = 1 - \cos\Big(\frac{\pi x}{2}\Big)\cos\Big(\frac{\pi y}{2}\Big) 
\end{equation}

For a given initial condition to be $f(x,y)$, $\textbf{a} = [1,1]$, and the periodic boundary conditions, we formulate the problem in the discontinuous Galerkin framework. The first step is projecting $f$ on a set of $E$ elements $\{e_0,e_1\cdots e_{E-1}\} \in \Omega$ using the typical nonorthogonal (hats and bubbles) basis $\phi$. Locally, for an element $e$, we have

$$f_e(x,y) = \sum_{i = 0}^{P-1} \hat{u}_i \phi_i(x,y),$$
where $\Omega$ is the mesh shown in \Cref{fig:para2dadv} consisting of a mix of quadrilaterals and triangles. {The filter steps do not change for different element types as emphasized by the choice of composite mesh.}\\

The solution for the advection problem using the timestep {$\Delta t = 1e-3$,} polynomial order 4, RK-4 integration scheme, and upwind flux calculations is shown in \Cref{fig:para2dadv}. Row 2 of the figure shows the simulation state at a particular timestep and highlights the negative values in the domain. In row 3 the non-negative structure is restored after the application of the filter, 
The filter is applied at each timestep to preserve the structure (positivity) of the solution on a lattice of points of interest. The lattice is a set of points defined in \Cref{eq:lattice}. If a structure violation is found at any point in the lattice, the parent element is flagged for filtering. Since the boundaries are periodic, the final state of the simulation looks similar to the initial state.\\

For a similar analysis of the 3D advection problem, consider the initial state as a smooth continuous function \Cref{eq:3dadvsimple} on a cube mesh of 64 hexahedron elements defined in the domain $[-1,-1,-1]\times[1,1,1]$. The advection velocity is defined by \textbf{a} $ = [1,1,1]$, and all the boundary conditions are periodic. The integration method used is RK-4 with a timestep of {$1e{-3}$}. The timestepping is performed for a total of 2000 timesteps and the flux calculation is upwind.
Following a similar selection procedure as in 2D, we discover the elements that violate the structure at each timestep and apply the structure-preserving filter to those elements.\\
\begin{equation}\label{eq:3dadvsimple}
f(x,y,z) = 1 - \cos\Big(\frac{\pi x}{2}\Big)\cos\Big(\frac{\pi y}{2}\Big)\cos\Big(\frac{\pi z}{2}\Big)  
\end{equation}

\begin{figure}[h!]
\centering
\includegraphics[width=0.7\textwidth]{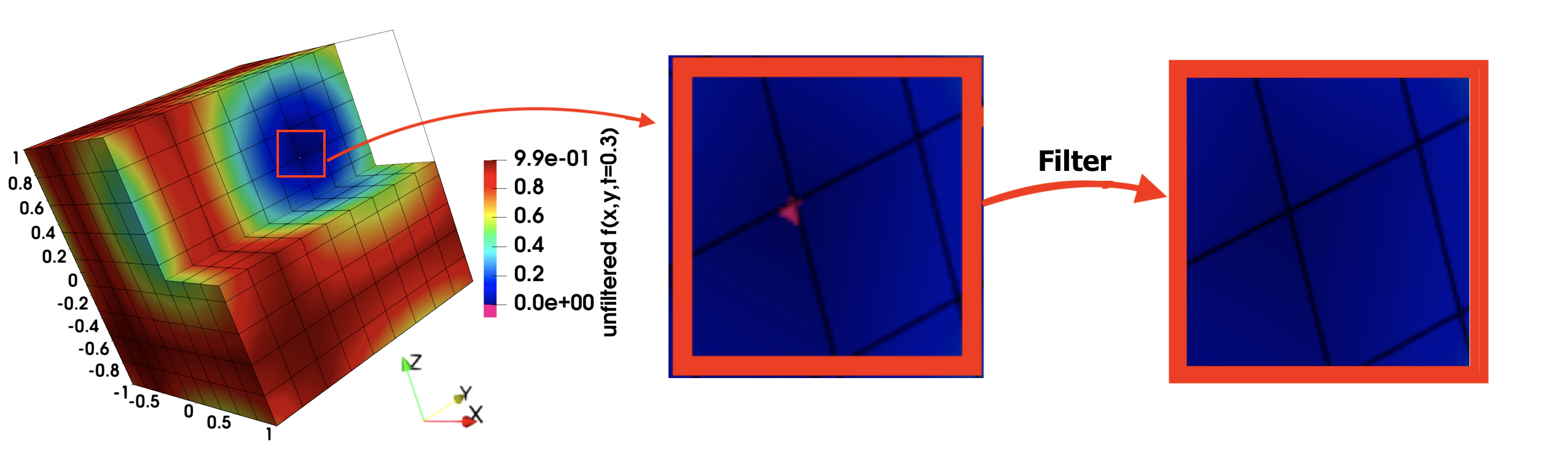}
\caption{The state of the simulation at t = 0.3 seconds ({timestep = 300}) using initial condition \Cref{eq:3dadvsimple}. \textit{Left and middle:} The highlighted region and its zoomed counterpart show the points in the elements that lose the positivity structure. \textit{Right:} The filtered values at the same elements that preserve the positivity structure.}
  \label{fig:simple3dadv}
\end{figure}

\Cref{fig:simple3dadv} shows an instance during the advection process where a loss of structure is encountered, i.e., negative intermediate values are found, as shown by the highlighted region in the left and middle subfigures of \Cref{fig:simple3dadv}. 
The process of discovering the structurally nonconformal elements and application of the filter to those elements adds to the total time cost of the advection solution, which is an important aspect to consider when deciding on the criteria for applying the filter. \\

\begin{figure}[h!]
\centering
\includegraphics[width=0.48\textwidth]{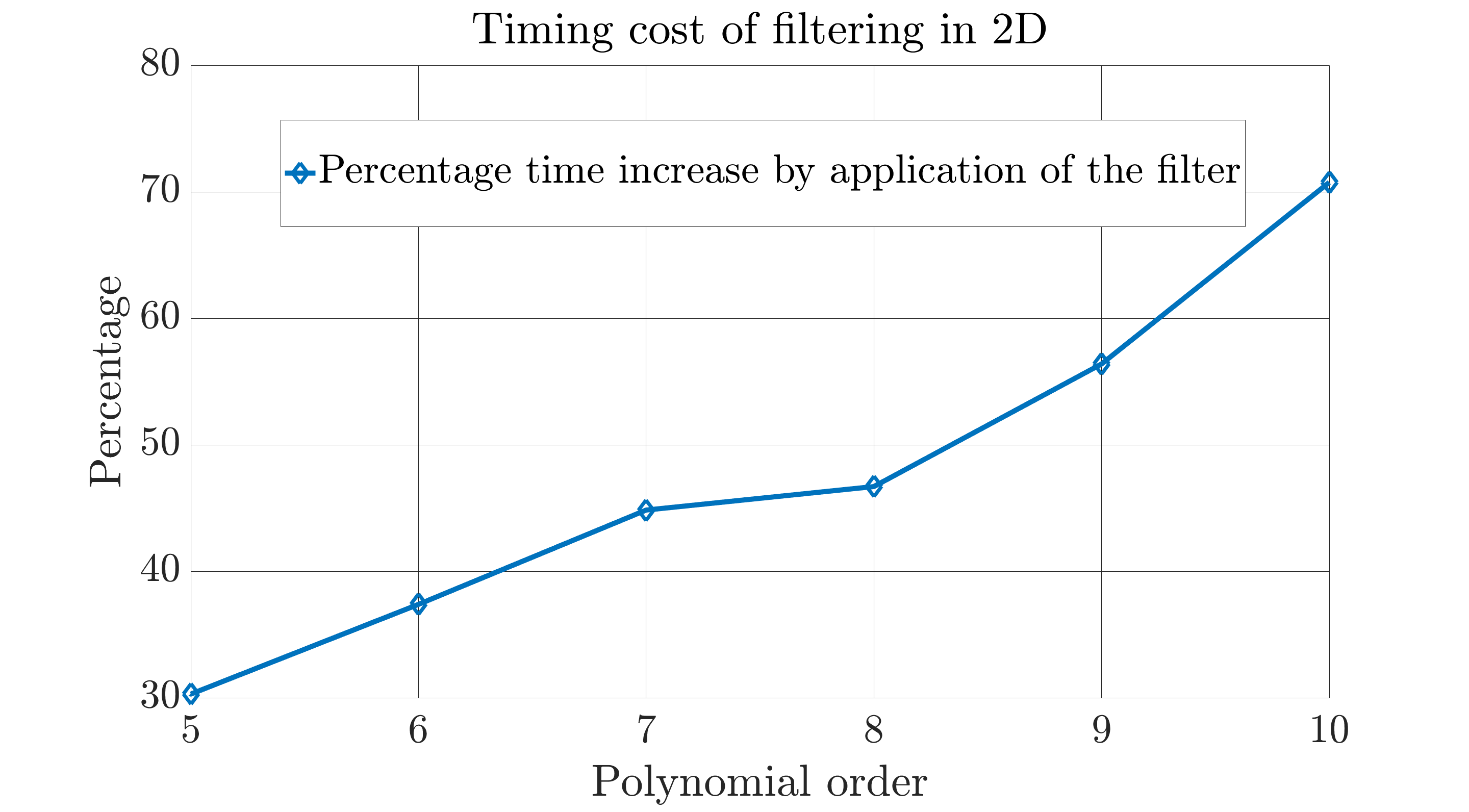}
\includegraphics[width=0.48\textwidth]{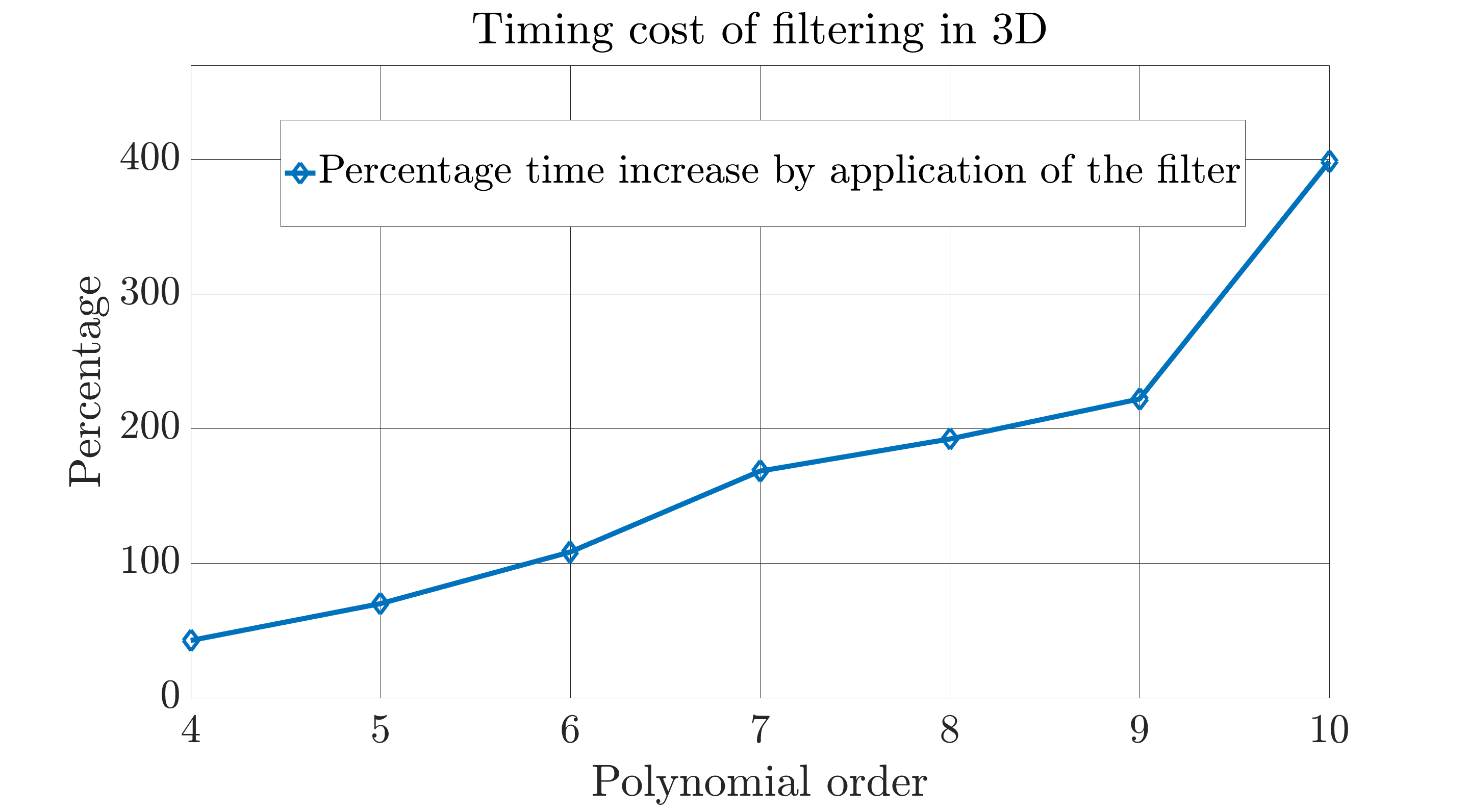}
\caption{Percentage increase in time to complete the experiment by application of the filter. \textit{Left:} 2D \textit{Right:} 3D}
  \label{fig:alltimings}
\end{figure}

The percentage increase in total time taken by the experiment by the application of structure-preserving filter is shown in \Cref{fig:alltimings}. Note that for 2D, at low orders (N=2,3,4), the cost of filtering is between 50-70$\%$ of the total time.  This is due, in part, to the efficiency at these orders of the unfiltered solver such as caching effects, as well as up-front costs associated with the optimization.  We find that for 2D the ratio of filtered to unfiltered is lowest at N=5, and then proceeds to climb linearly as shown in the left subfigure of \Cref{fig:alltimings}. The time taken by the filter in 3D is notably higher than the time taken by the filter in 2D because of the higher complexity of finding the global minimum in 3D. The cost of overall filtering per timestep depends on the number of elements that the filter operates upon. Therefore, the procedure of selective application of the filter becomes increasingly important. \\
%

\subsection{{Structure preservation of the canonical rotating solid body test}}\label{sub:levequetests}

\begin{figure}[h]
\centering
\includegraphics[width=0.5\textwidth]{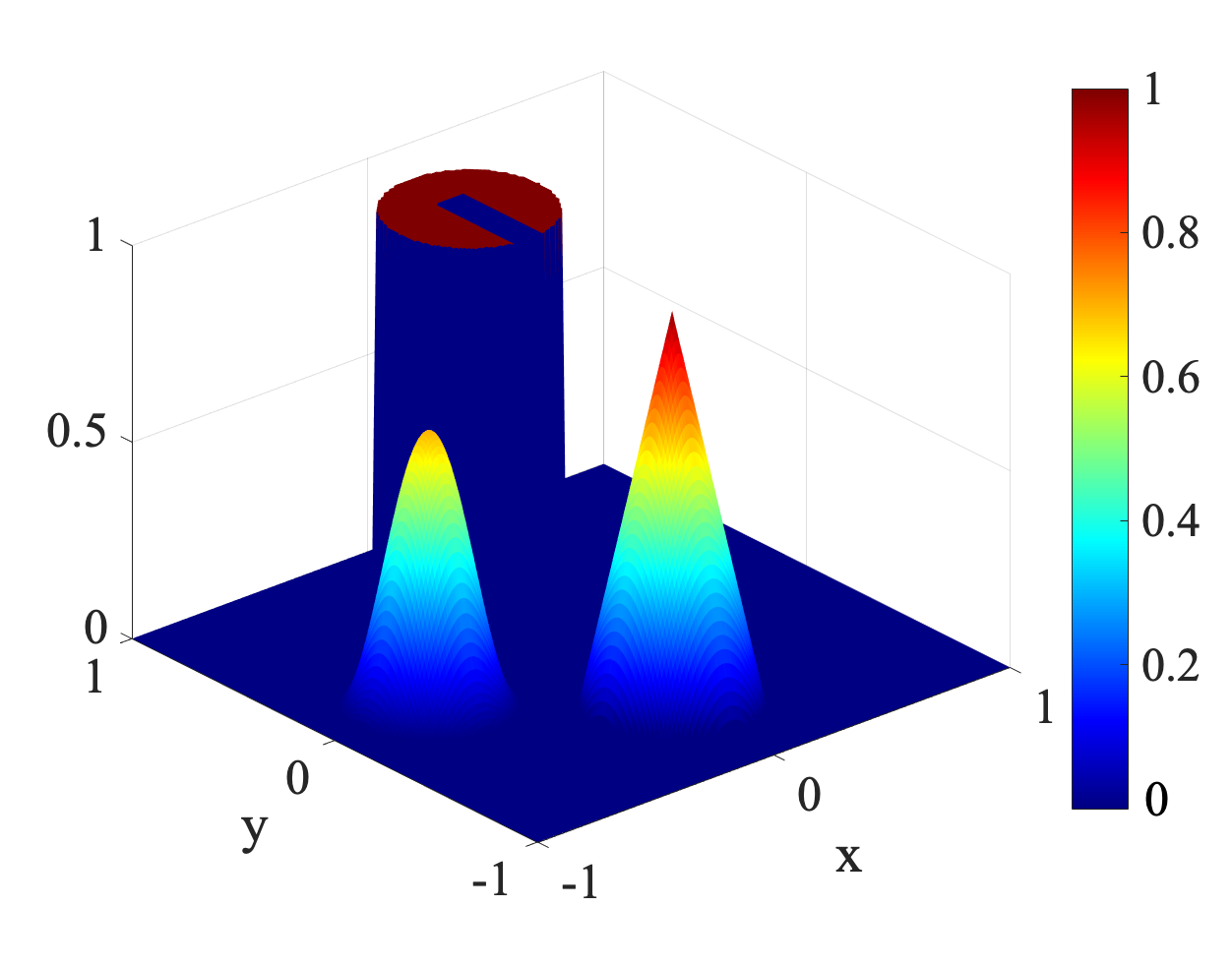}
\caption{Initial state for solid body rotation tests. }
  \label{fig:newfn1}
\end{figure}

We now investigate the application of the filter on the solid body rotation experiment using a discontinuous initial condition defined as a combination of a notched cylinder, a cone, and a smooth hump. Consider the initial data shown in \Cref{fig:newfn1}. This example is extensively used in advection and structure preservation literature \cite{leveque1996high}.  The parameters to reproduce the initial state: the cylinder, cone, and hump have a radius of 0.3 each. The centers of the cylinder, cone, and hump are (0, 0.5), (0,-0.5), and (-0.6,0) respectively.  Without changing the domain and boundary conditions used for the 2D advection experiment shown in \Cref{fig:para2dadv},  the advection velocity is changed to circular such that in one time period, the solution returns to its original state.  For all tests in this section, the domain has 144 quadrilateral elements and the timestep {$(\Delta t)$ is $5e-4$.}

\begin{figure}[t]
\centering
\includegraphics[width=0.32\textwidth]{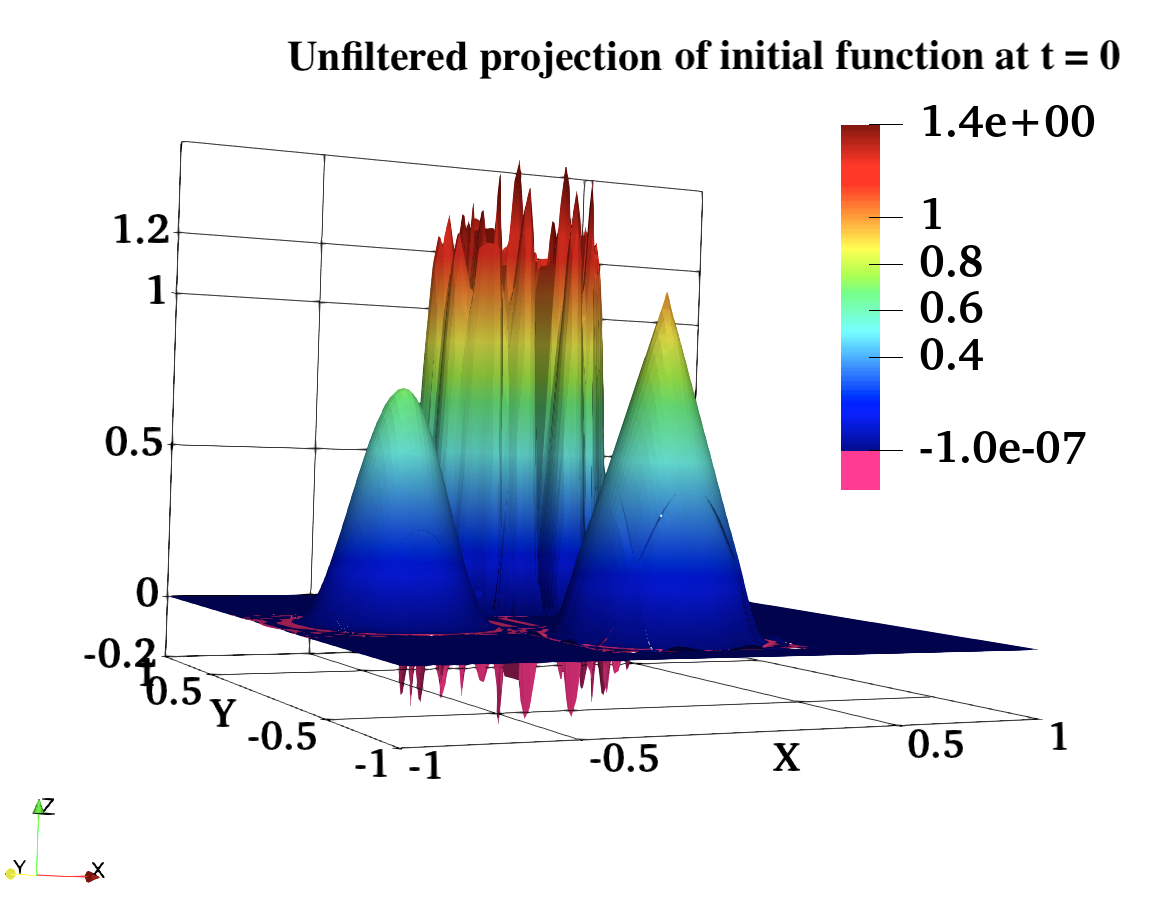}
\includegraphics[width=0.33\textwidth]{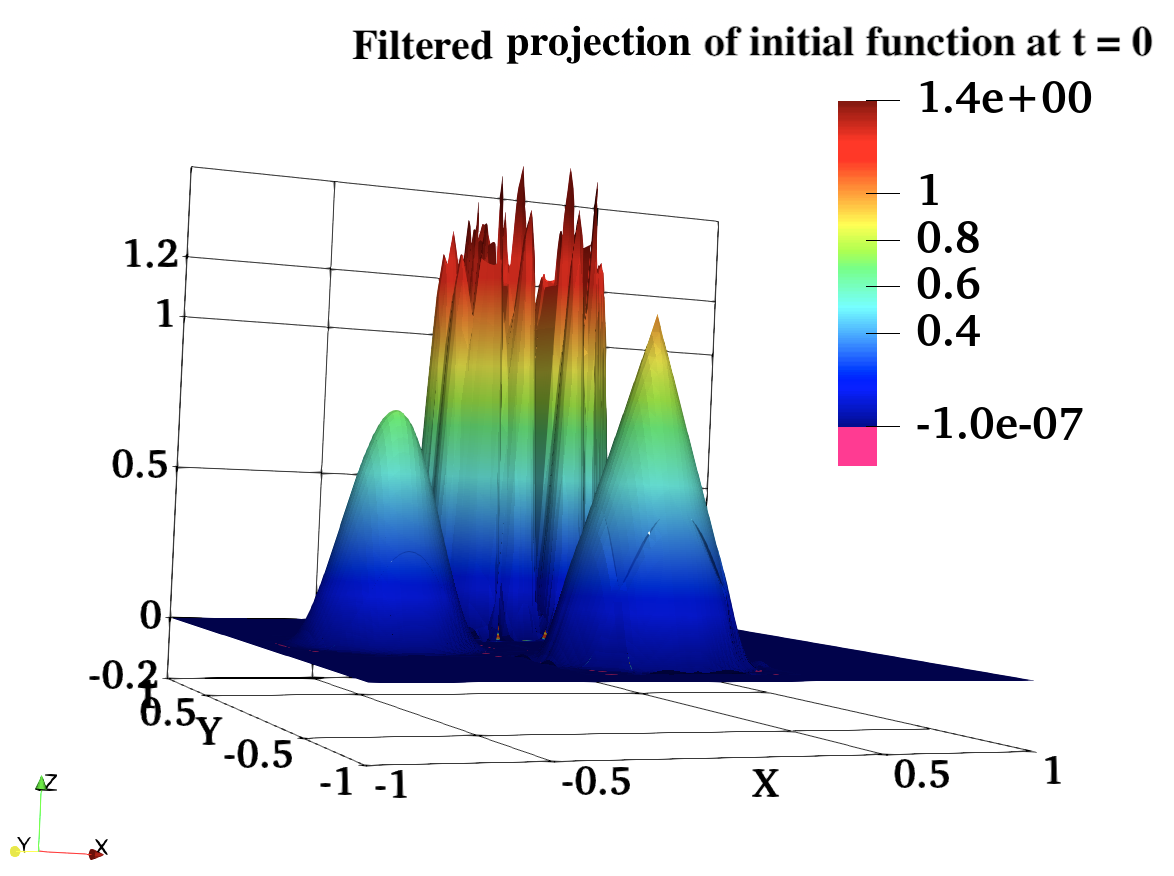}\\
\includegraphics[width=0.32\textwidth]{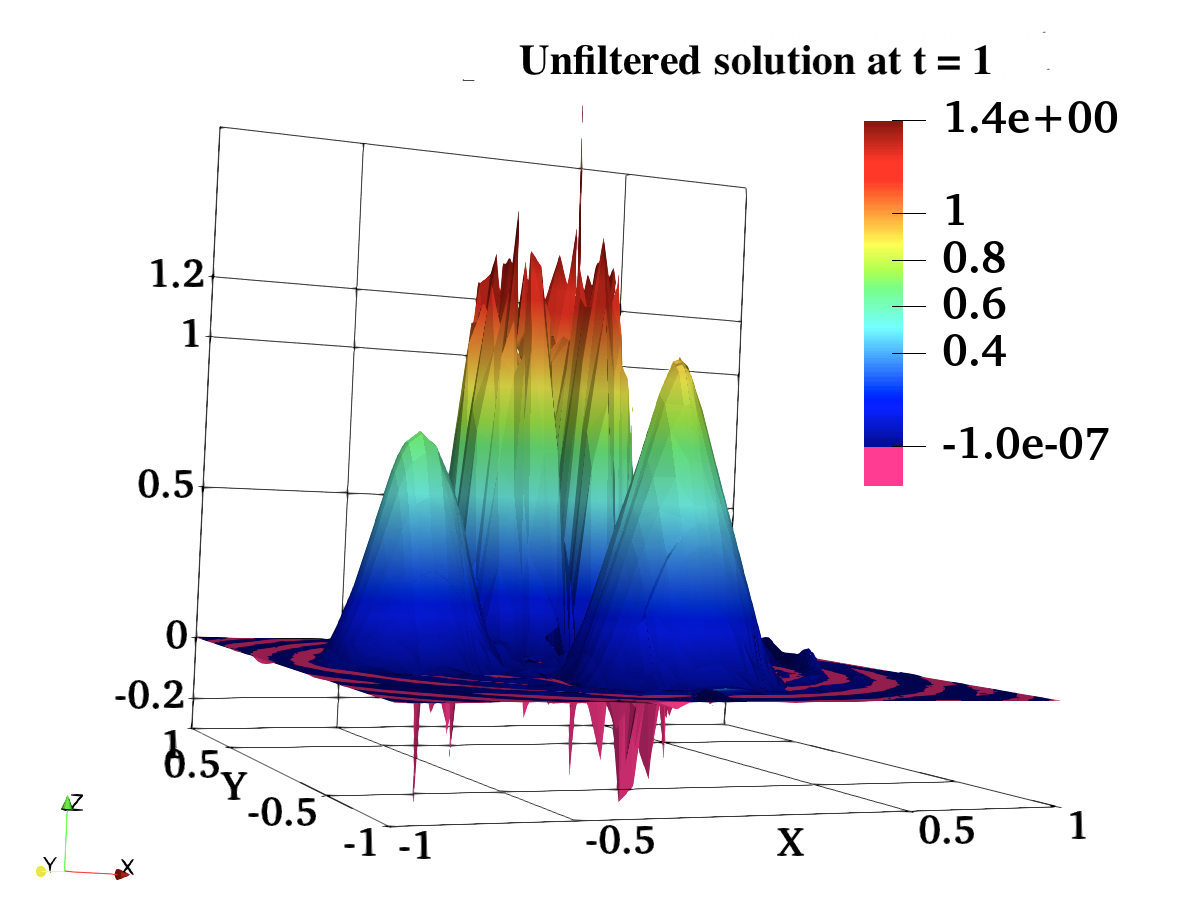}
\includegraphics[width=0.33\textwidth]{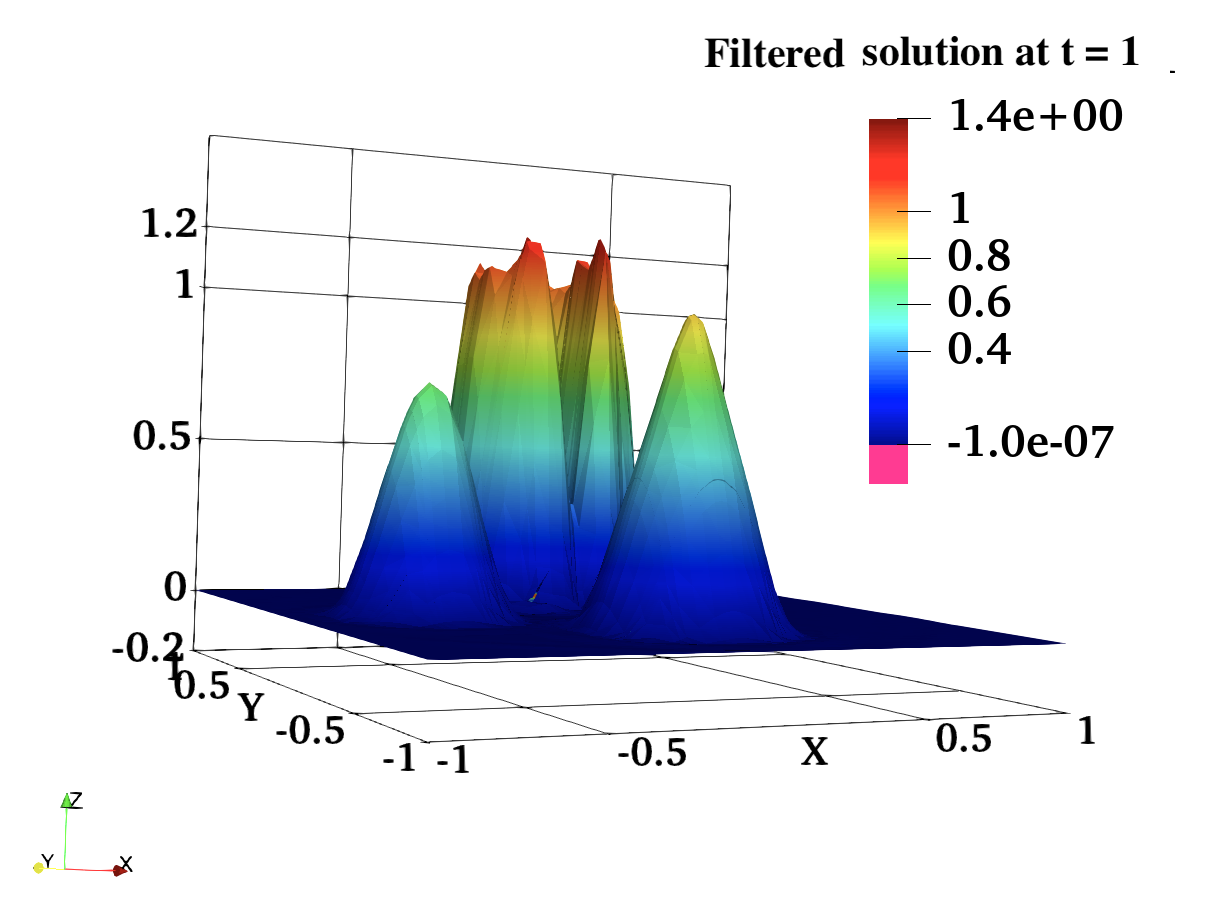}
\caption{{Snapshots at the beginning and end of advection. The highlighted region shows the values below tolerance for negativity $(10^{-7})$.  At t = 1s, the solid body finishes one rotation and returns to the original position.  Parameters for the test: {Timestep = $5e-4$,} polynomial order = 5. 
\textit{Row 1:} {Initial function projection at time = 0s. } \textit{Row 2:} Solution at time = 1s. \textit{Column 1:} Unconstrained. \textit{Column 2:} Constrained.  }}
  \label{fig:newfn}
\end{figure}

\begin{figure}[h!]
\centering
\includegraphics[width=0.36\textwidth]{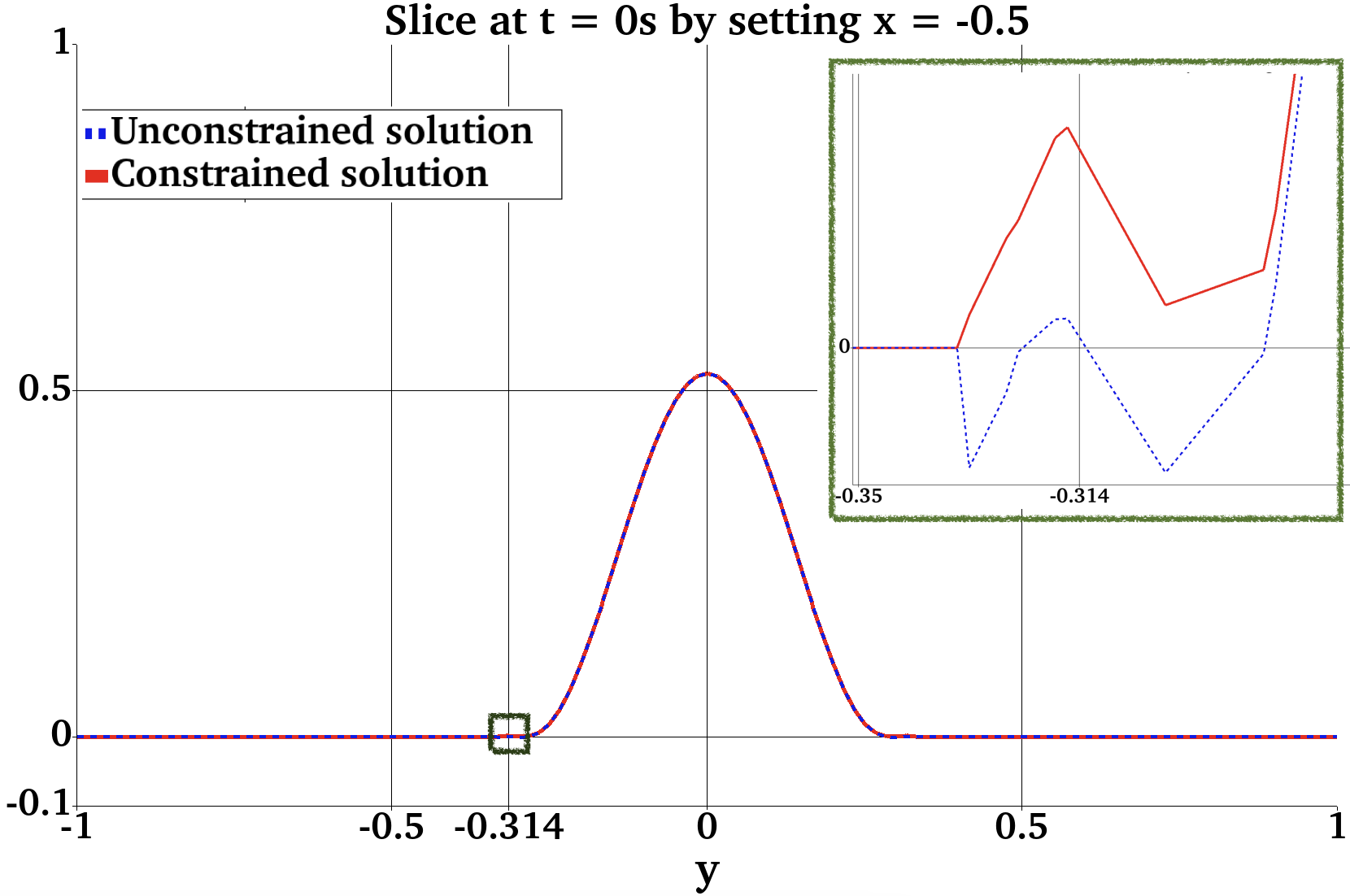}
\includegraphics[width=0.36\textwidth]{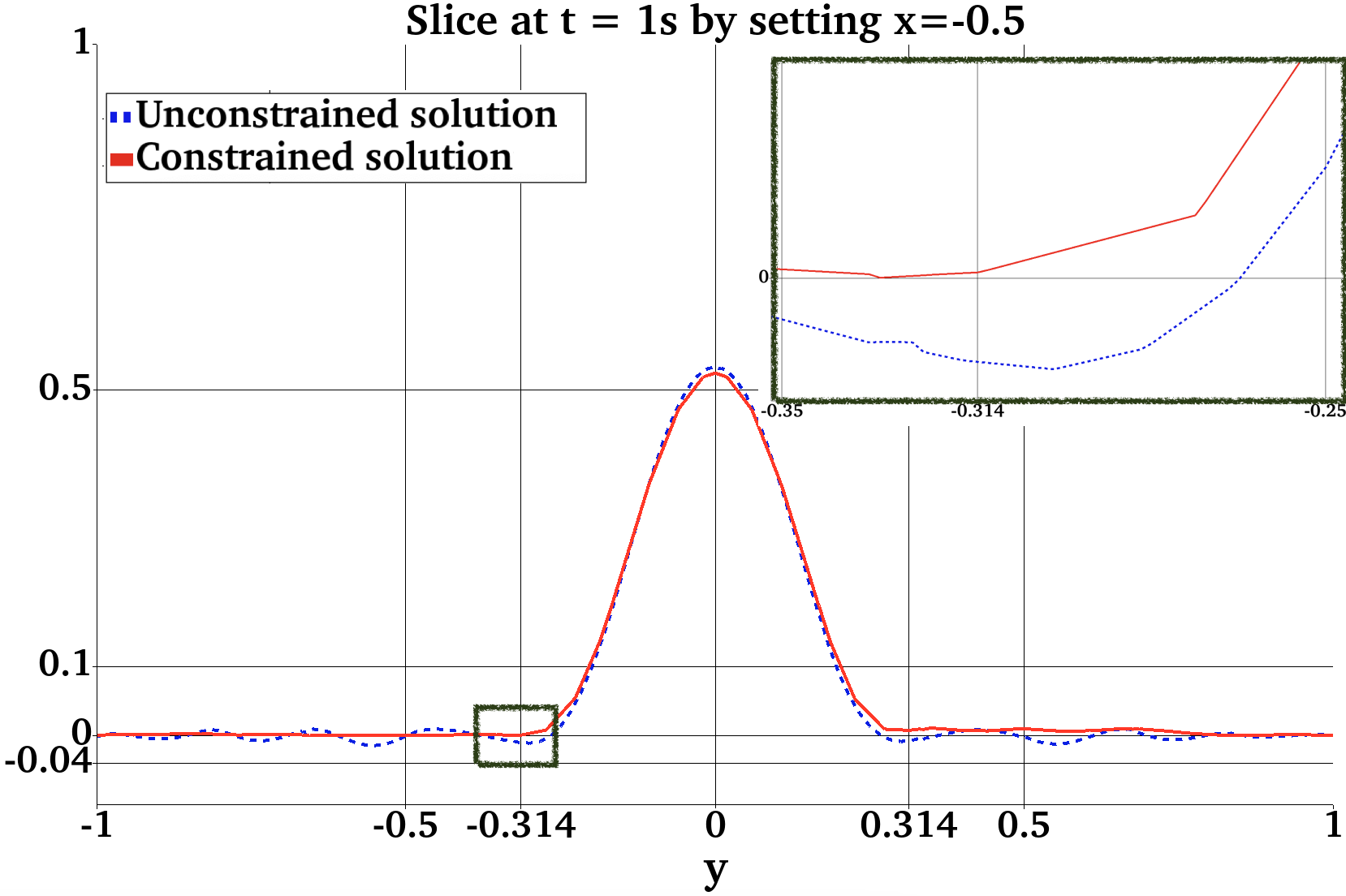}\\
\includegraphics[width=0.36\textwidth]{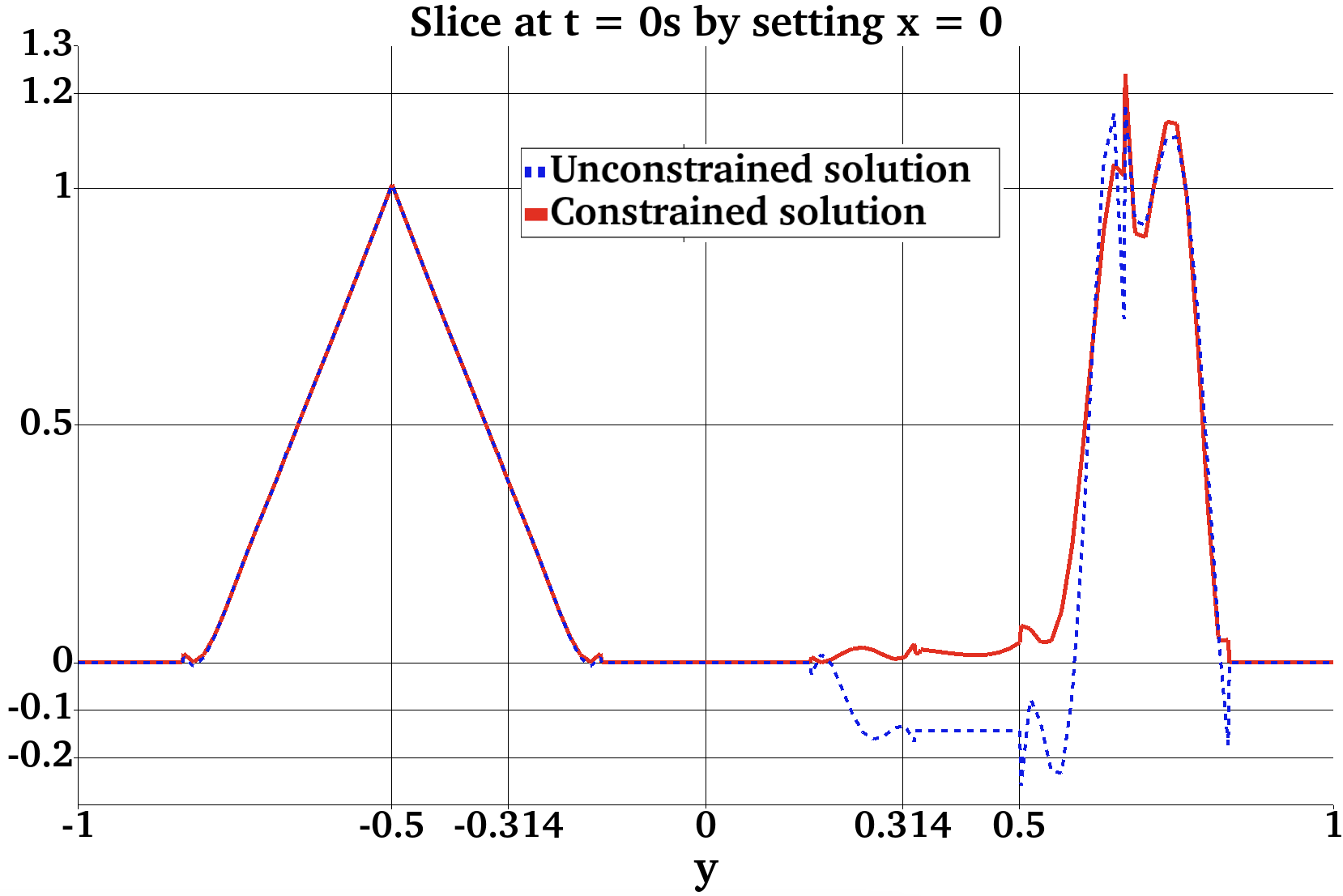}
\includegraphics[width=0.36\textwidth]{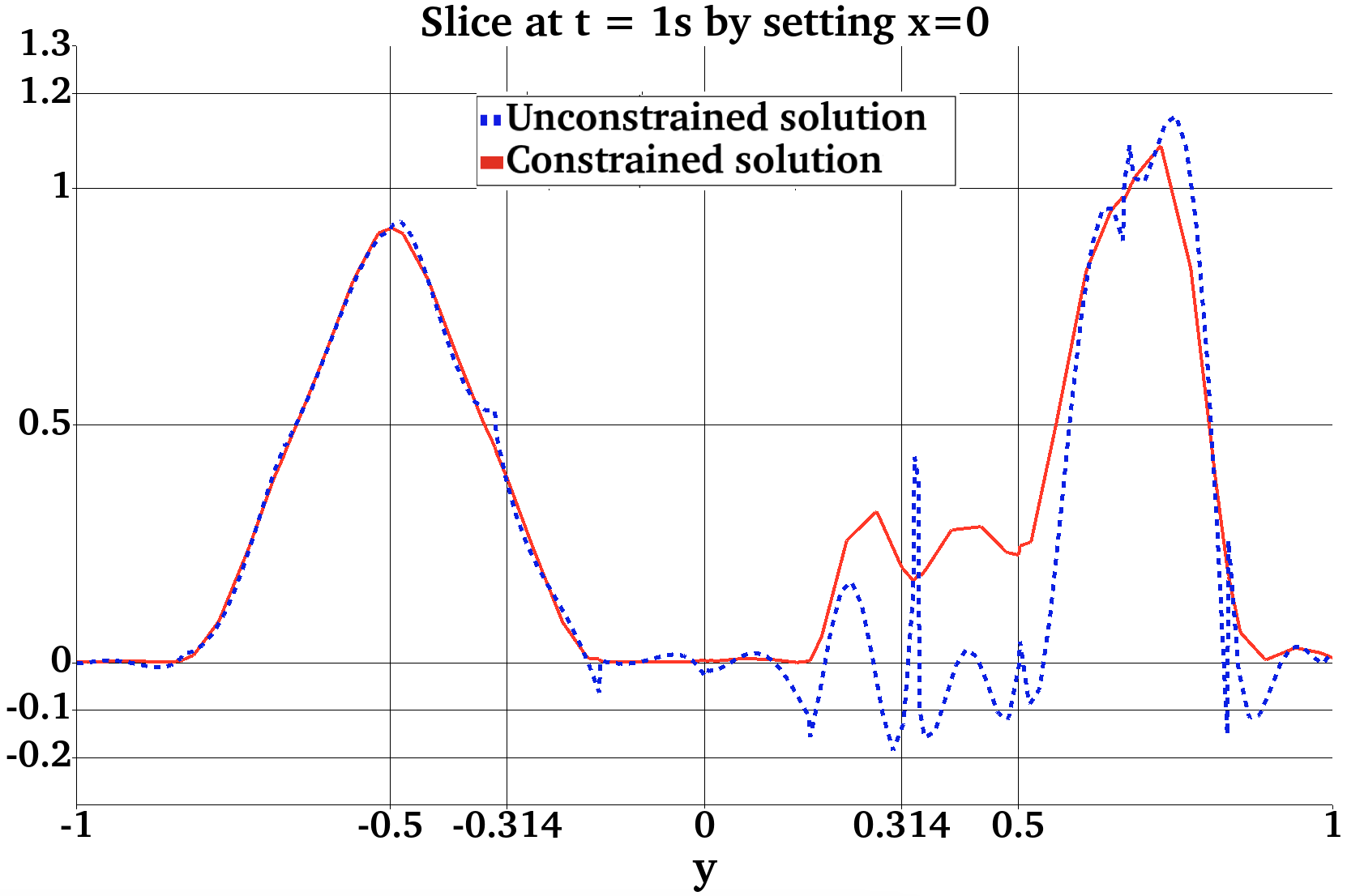}
\caption{{Snapshots at various points in time during advection. {$\Delta t = 5e-4$}, time period = 1 and polynomial order = 5.  Some subfigures show insets with the zoomed-in area of interest.
\textit{Column 1:} State of the system at t = 0. \textit{Column 2:} State of the system at t = 1s.  \textit{Row 1:} Slice at x = -0.5 and \textit{Row 2:} Slice at x = 0.}}
  \label{fig:adv3dlevequeslice1}
\end{figure}

{Since the initial condition for the experiment is defined by a complicated discontinuous function, the solution suffers from a larger loss of structure as compared to the case with a smooth initial condition}.
{\Cref{fig:newfn} shows the initial and final states during rotation of the body shown in \Cref{fig:newfn1}, with and without the filter application. }\\

\begin{figure}[h!]
\centering
\includegraphics[width=0.36\textwidth]{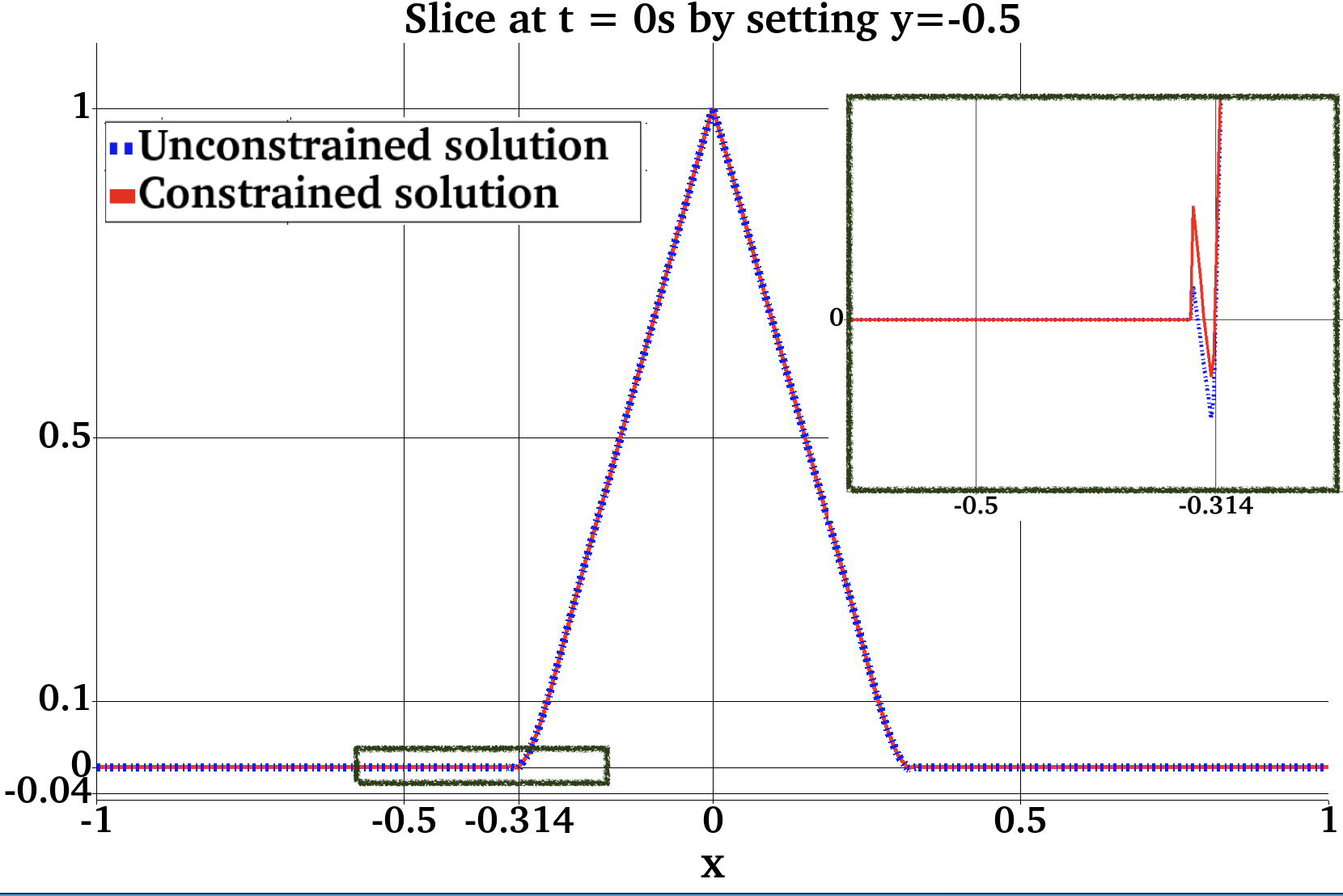}
\includegraphics[width=0.36\textwidth]{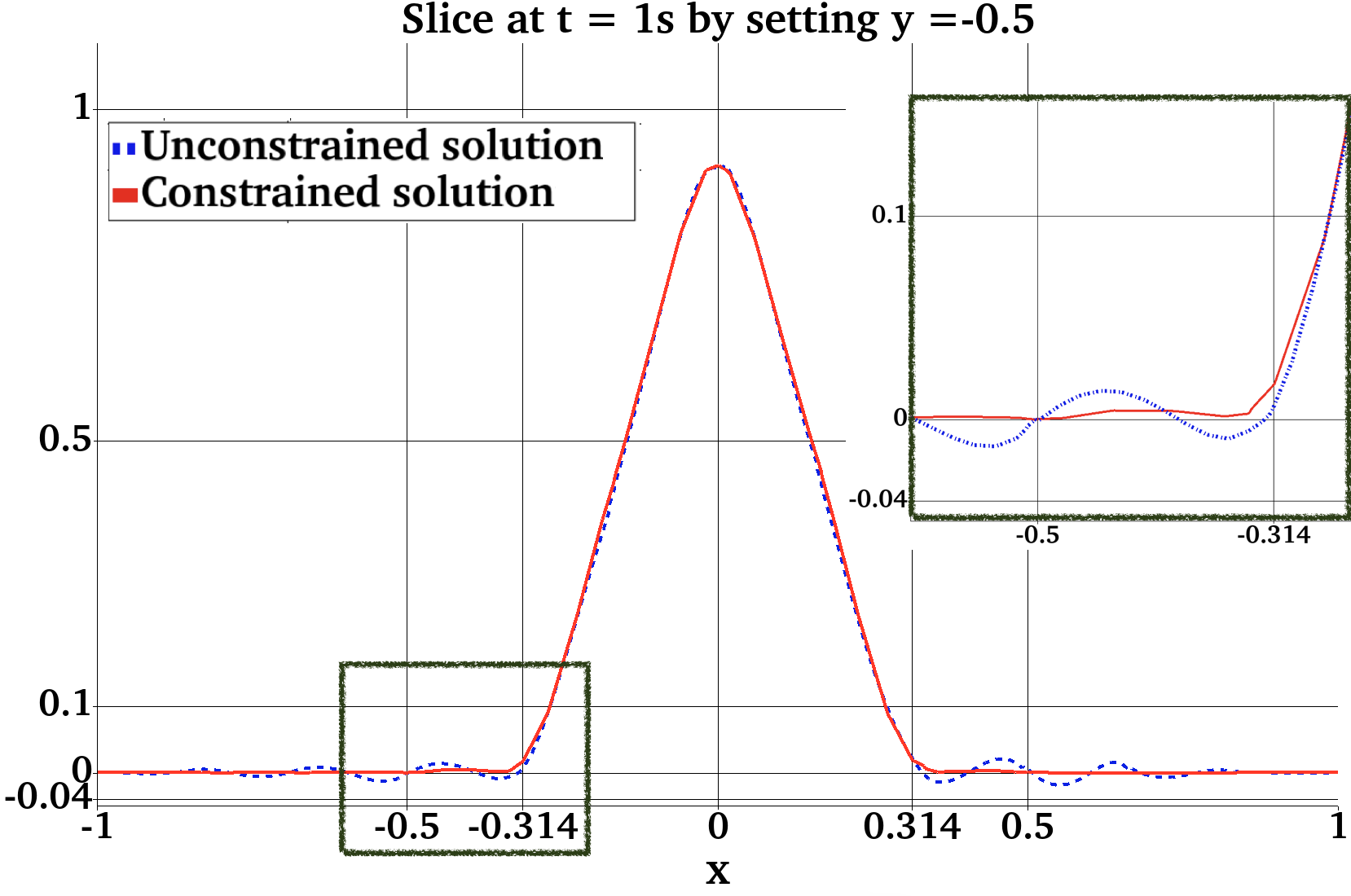}\\
\includegraphics[width=0.36\textwidth]{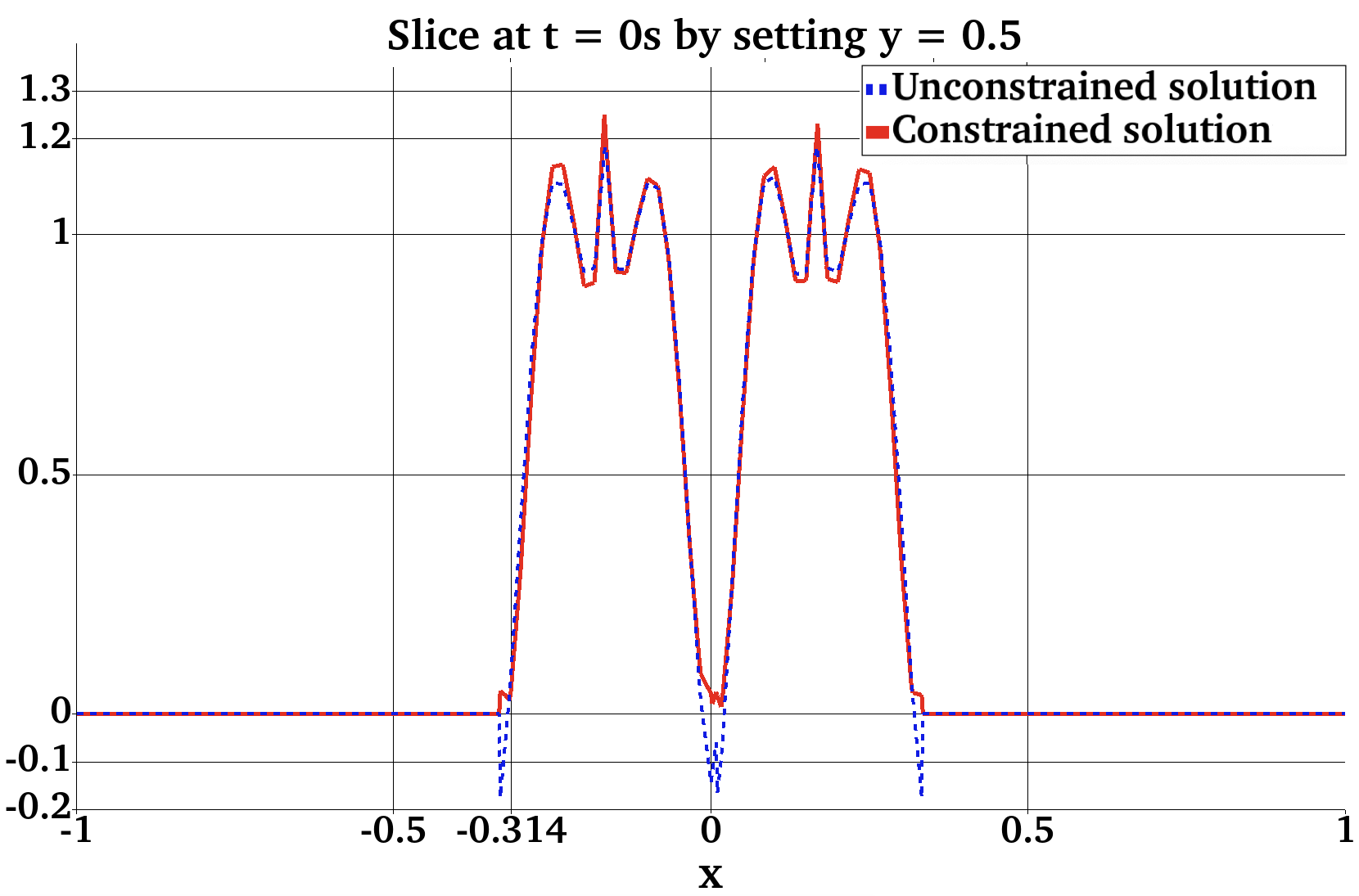}
\includegraphics[width=0.36\textwidth]{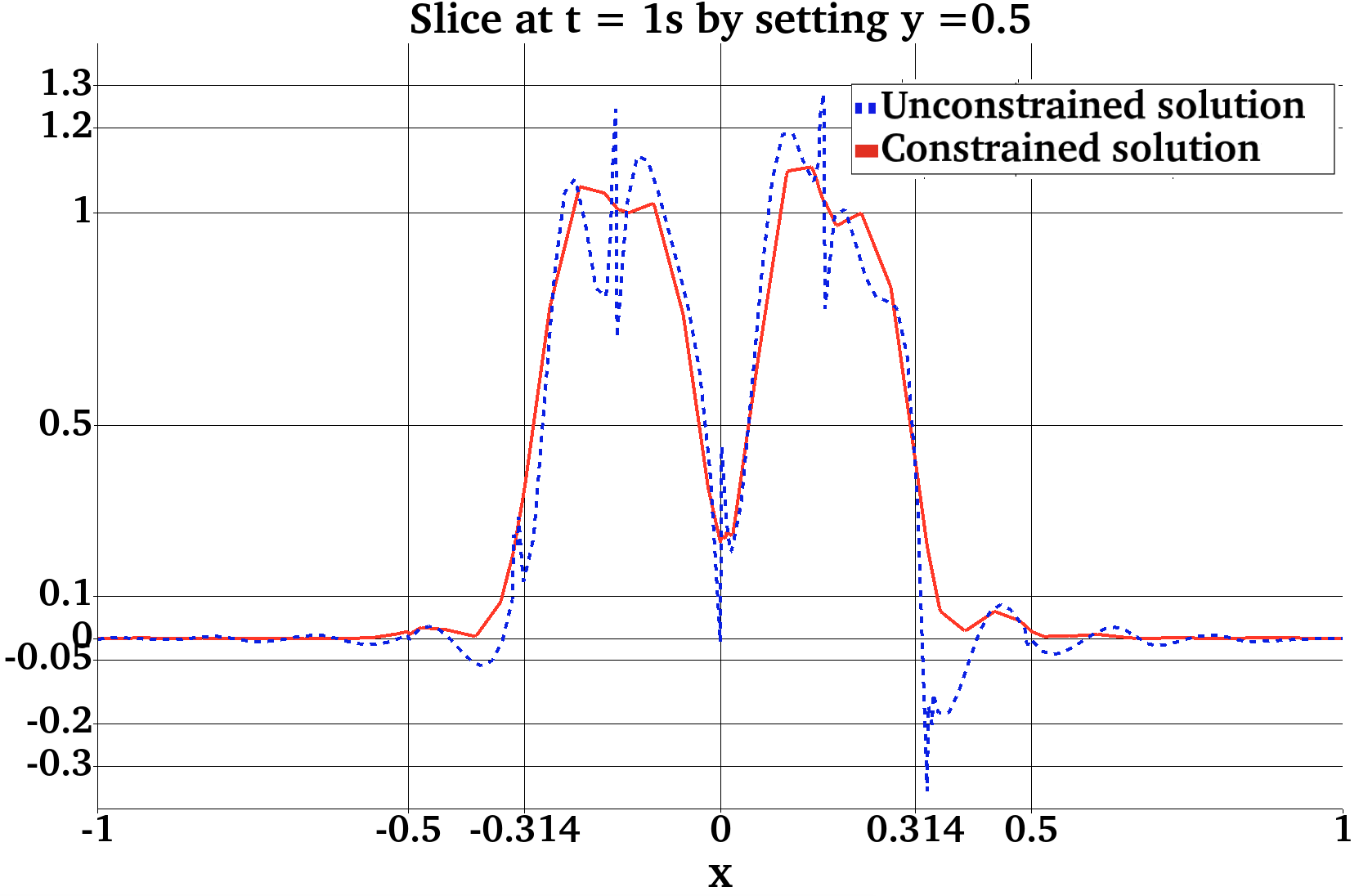}
\caption{{Snapshots at various points in time during advection. {$\Delta t = 5e-4$}, time period = 1 and polynomial order = 5.  Some subfigures show insets with the zoomed-in area of interest.
\textit{Column 1: }State of the system at t = 0. \textit{Column 2:} State of the system at t = 1. \textit{Row 1: }Slice at y = -0.5 and \textit{Row 2: }Slice at y = 0.5}}
  \label{fig:adv3dlevequeslice2}
\end{figure}
\begin{figure}[h!]
\centering
\includegraphics[width=0.49\textwidth]{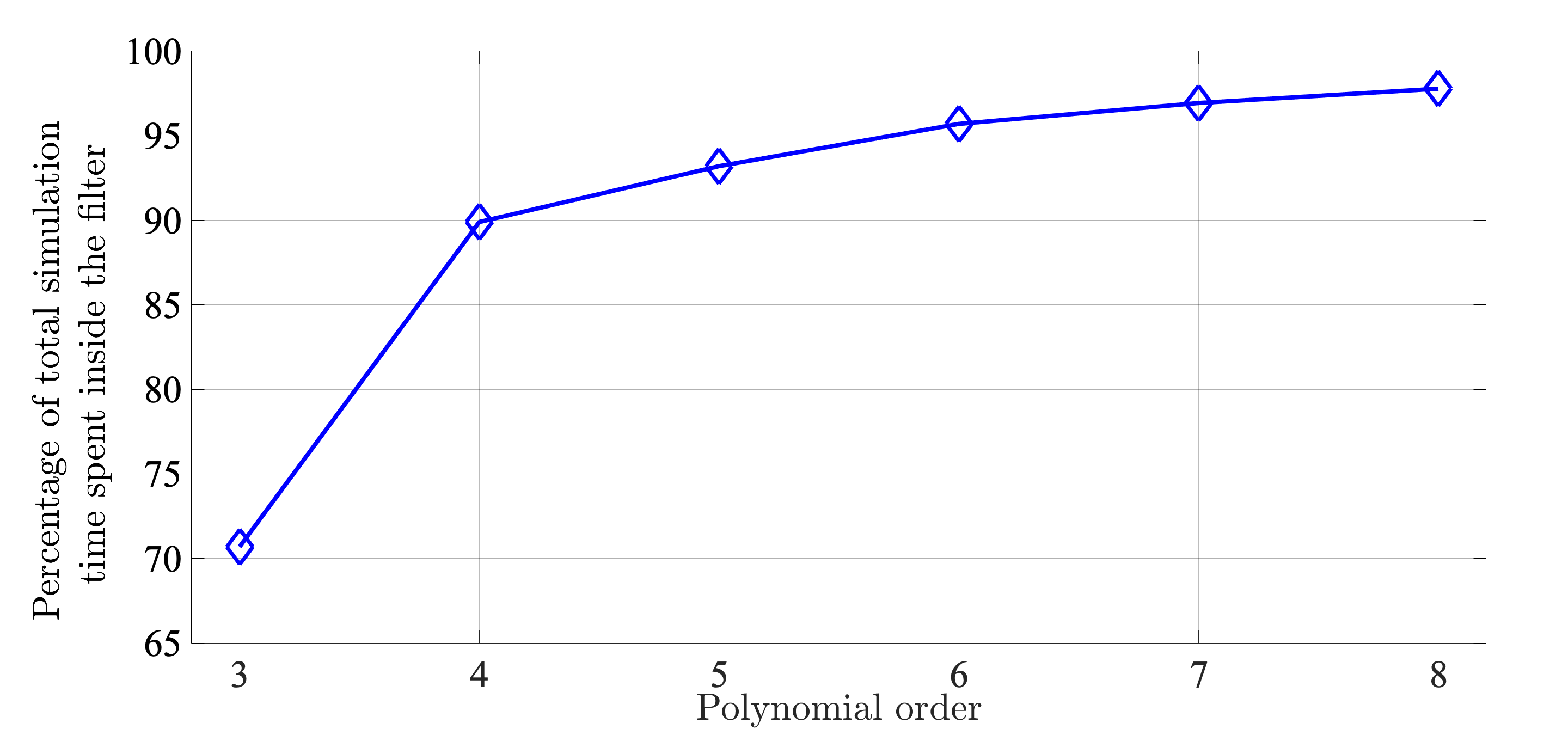}
\includegraphics[width=0.49\textwidth]{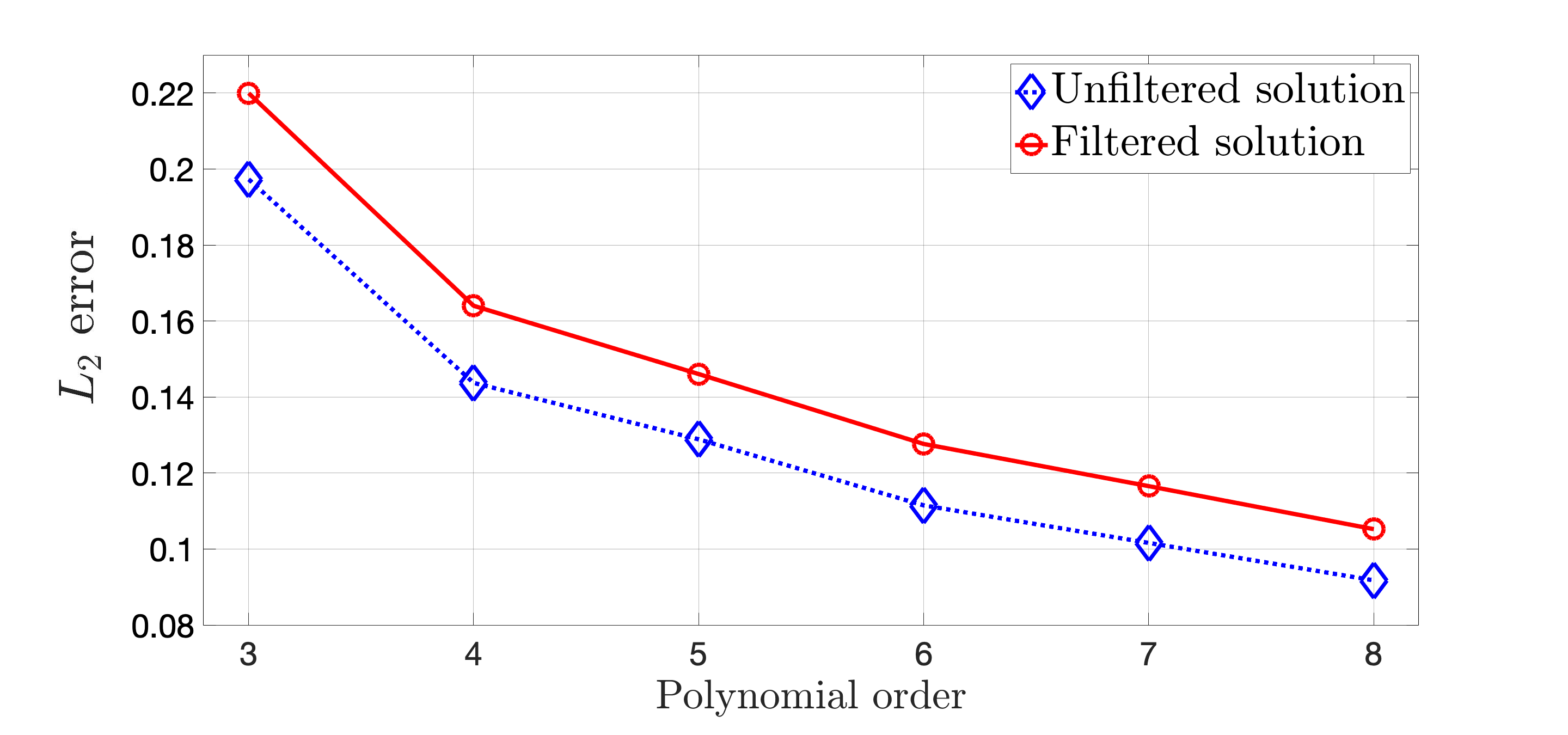}\\
\includegraphics[width=0.49\textwidth]{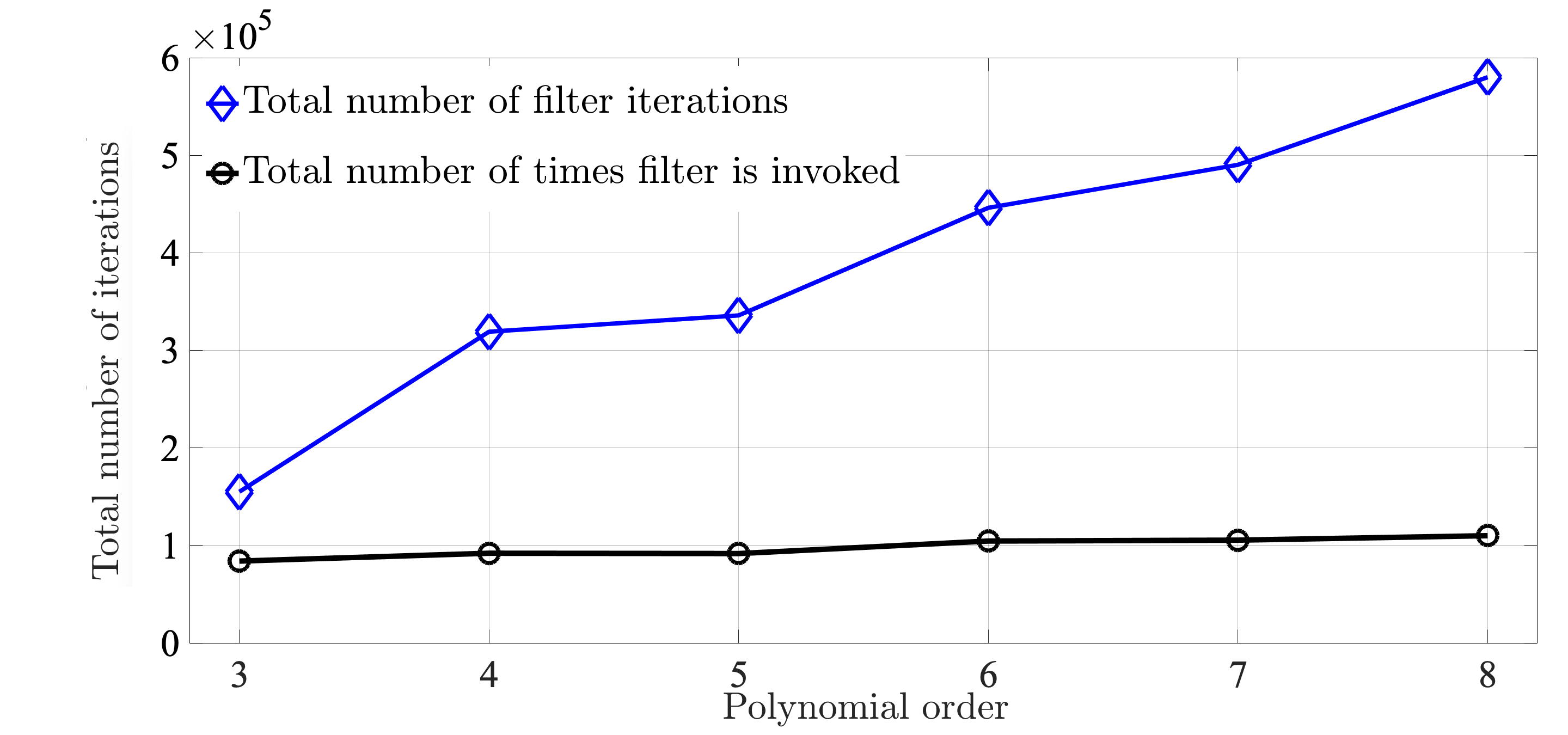}
\includegraphics[width=0.49\textwidth]{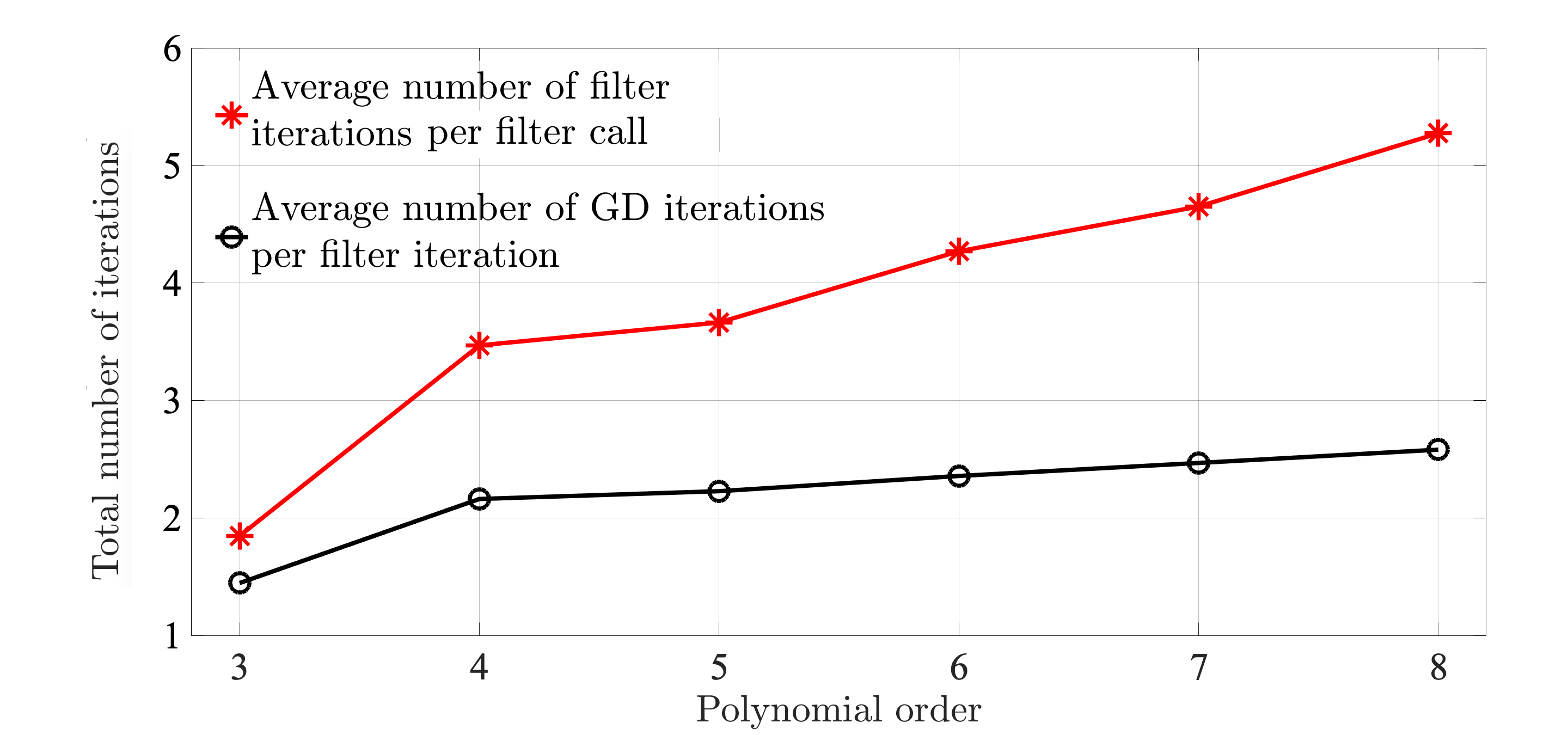}
\caption{{Analysis of filter application to the experiment shown in \Cref{fig:newfn1}.  \textit{Top left:} Percentage of total simulation time spent inside the filter,  \textit{Top right:} p-convergence of $L_2$ errors for filtered and unfiltered solution at t=1s. \textit{Bottom left:} Total number of filter iterations vs the total number of times the filter is called for various polynomial orders. \textit{Bottom right:} Average number of filter iterations per call vs the average number of GD line-search algorithm per filter iteration for various polynomial orders.}}
 \label{fig:resultsnewfn}
\end{figure}

 {\Cref{fig:adv3dlevequeslice1,fig:adv3dlevequeslice2} show different cross sections of the solution with and without filtering at the initial and final time. The total time taken by the filter as a function of the entire simulation time is shown in \Cref{fig:resultsnewfn}.\\} 

{As evident from the top left panel in \Cref{fig:resultsnewfn}, the filter takes proportionally longer to restore the structure as the loss of structure increases. In many practical applications, such as the one presented in \Cref{sec:results3}, we observe that paying the extra filtering cost results in the successful termination of the simulation, which otherwise fails due to the presence of structural inconsistencies. In such cases, the tradeoff between the time cost vs guarantee of structure preservation leans towards the latter. The filtered and unfiltered $L_2$ errors after one rotation of the solid body are shown in the top right panel of \Cref{fig:resultsnewfn}. Whereas filtering changes the final state of the solution, thus contributing to the errors,  the difference between filtered and unfiltered errors is still comparable.  The filtered version of the tests using different polynomial order follows the same p-convergence as the unfiltered counterpart.  A comparison between the number of times the filter is invoked and the number of iterations taken per invocation is shown in the bottom left panel of \Cref{fig:resultsnewfn}. We notice that the number of iterations taken per call to the filter increases with the increase in polynomial order because of the larger magnitude of structural inconsistencies observed as a result of the oscillations in the solution. The largest contributor to the cost of applying the filter occurs from the GD line-search to find the global minimum.  At each iteration of the filter, we find the minimum.  Each call to find the minimum, in turn, takes a few iterations of GD. On the bottom right panel of \Cref{fig:resultsnewfn} is a comparison between the average number of GD iterations per iteration of the filter.  Also shown is the average number of filter iterations per call to the filter.}

\subsection{Structure preservation in advection tests on various 3D domains}

The procedure to filter the elements that exhibit loss of structure remains agnostic to the type of element. To demonstrate this, consider different tesselations of the domain $ [-1,-1,-1]\times[1,1,1]$ and solve an advection problem with the nontrivial initial state \Cref{eq:f3adv}. The setup remains the same as the previous 3D example, i.e., \textbf{a} $ = [1,1,1]$ in all directions, periodic boundaries, RK-4 integration scheme using timestep {$1e{-3}$}, and for a total of 2000 timesteps.

\begin{figure}[h]
\centering
\includegraphics[width=0.26\textwidth]{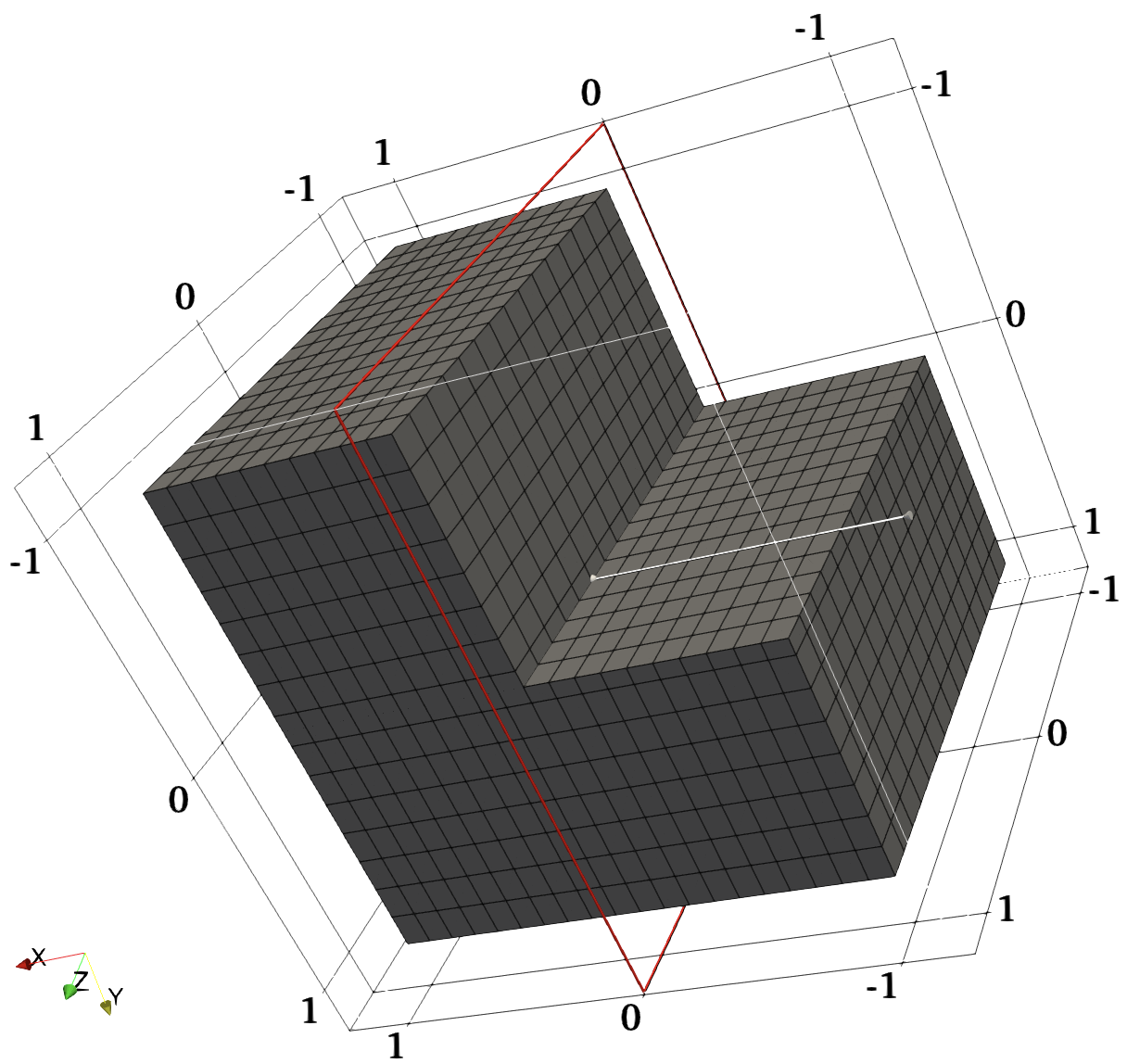}
\includegraphics[width=0.31\textwidth]{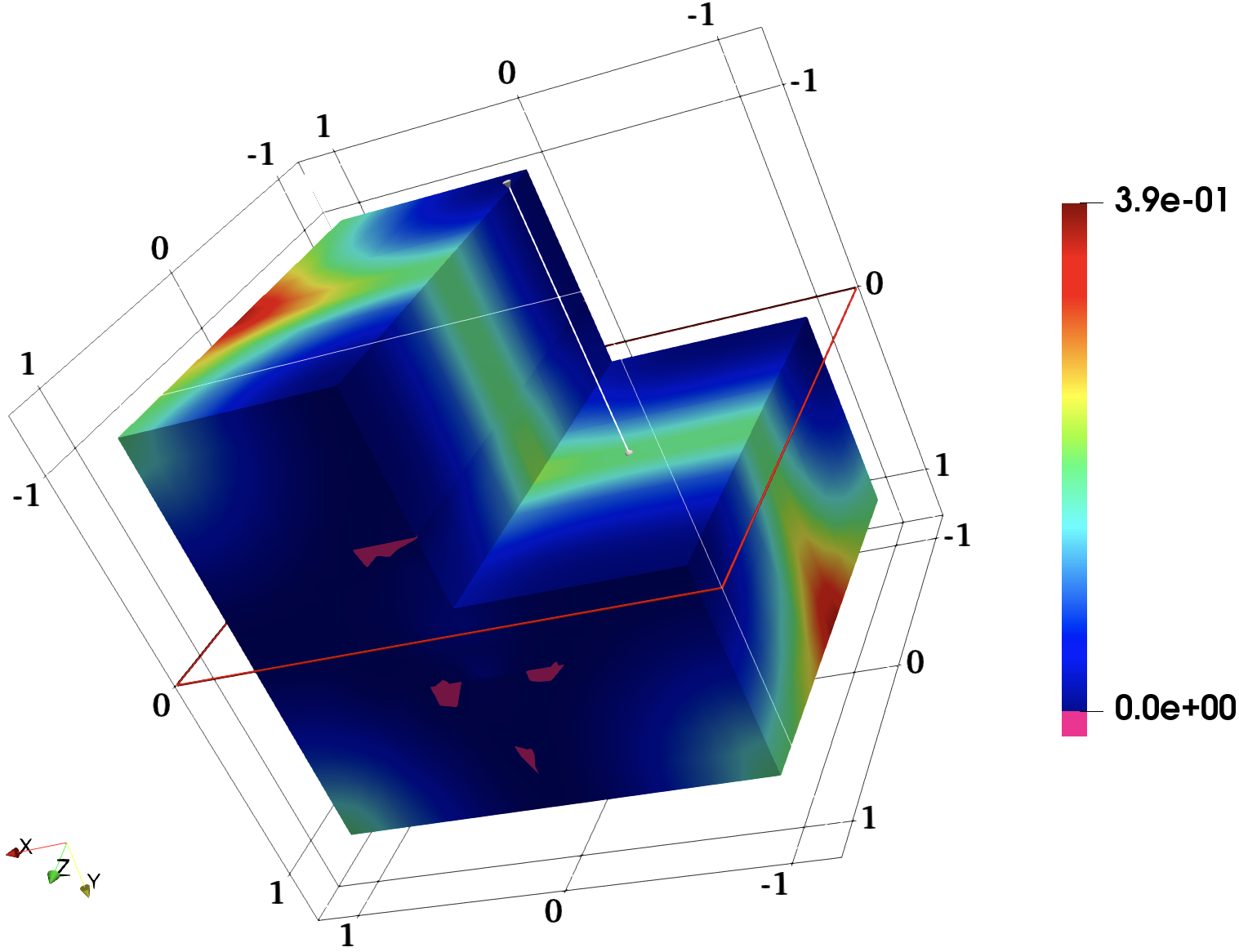}
\includegraphics[width=0.26\textwidth]{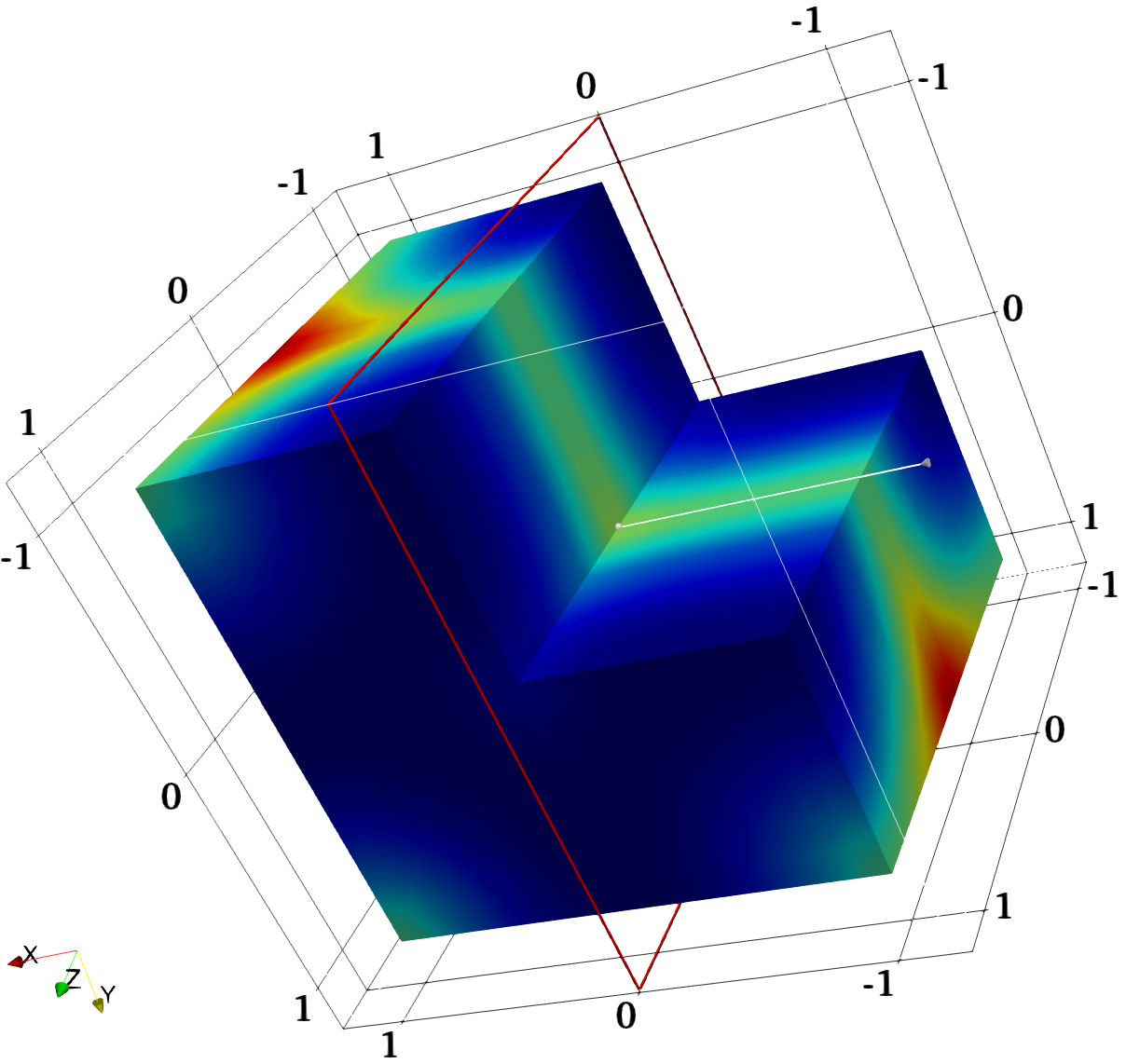}\\	
\includegraphics[width=0.26\textwidth]{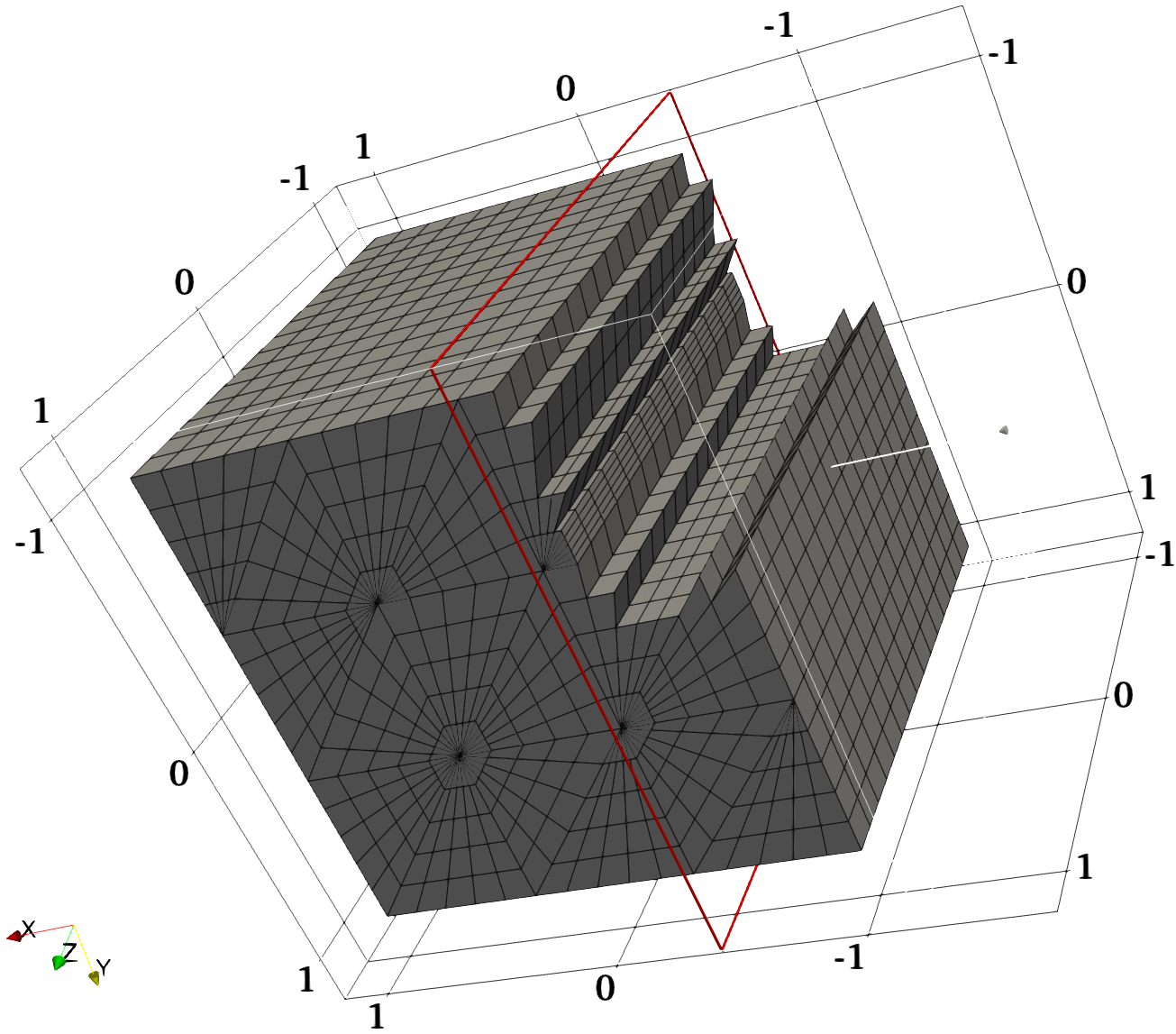}
\includegraphics[width=0.31\textwidth]{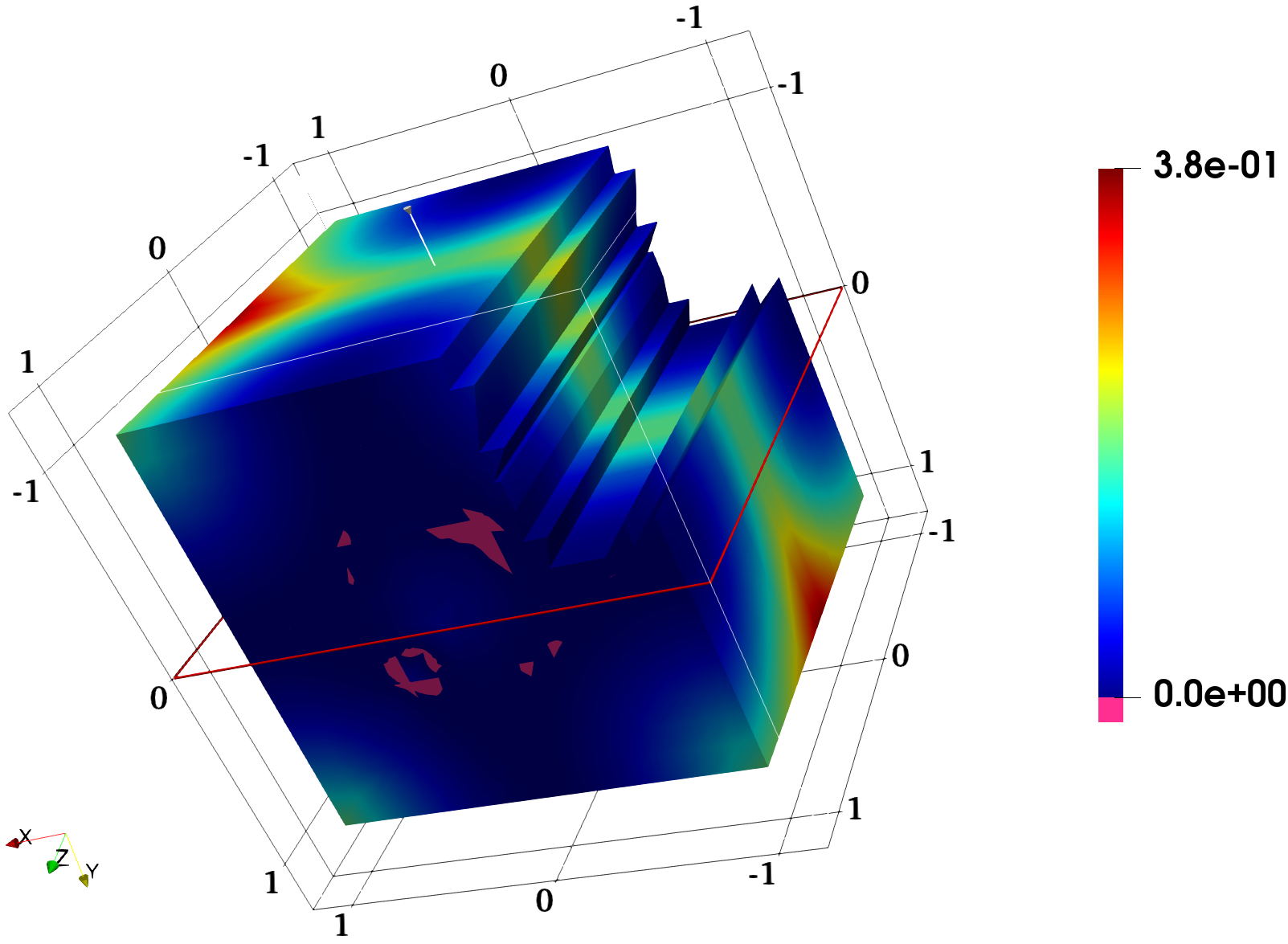}
\includegraphics[width=0.26\textwidth]{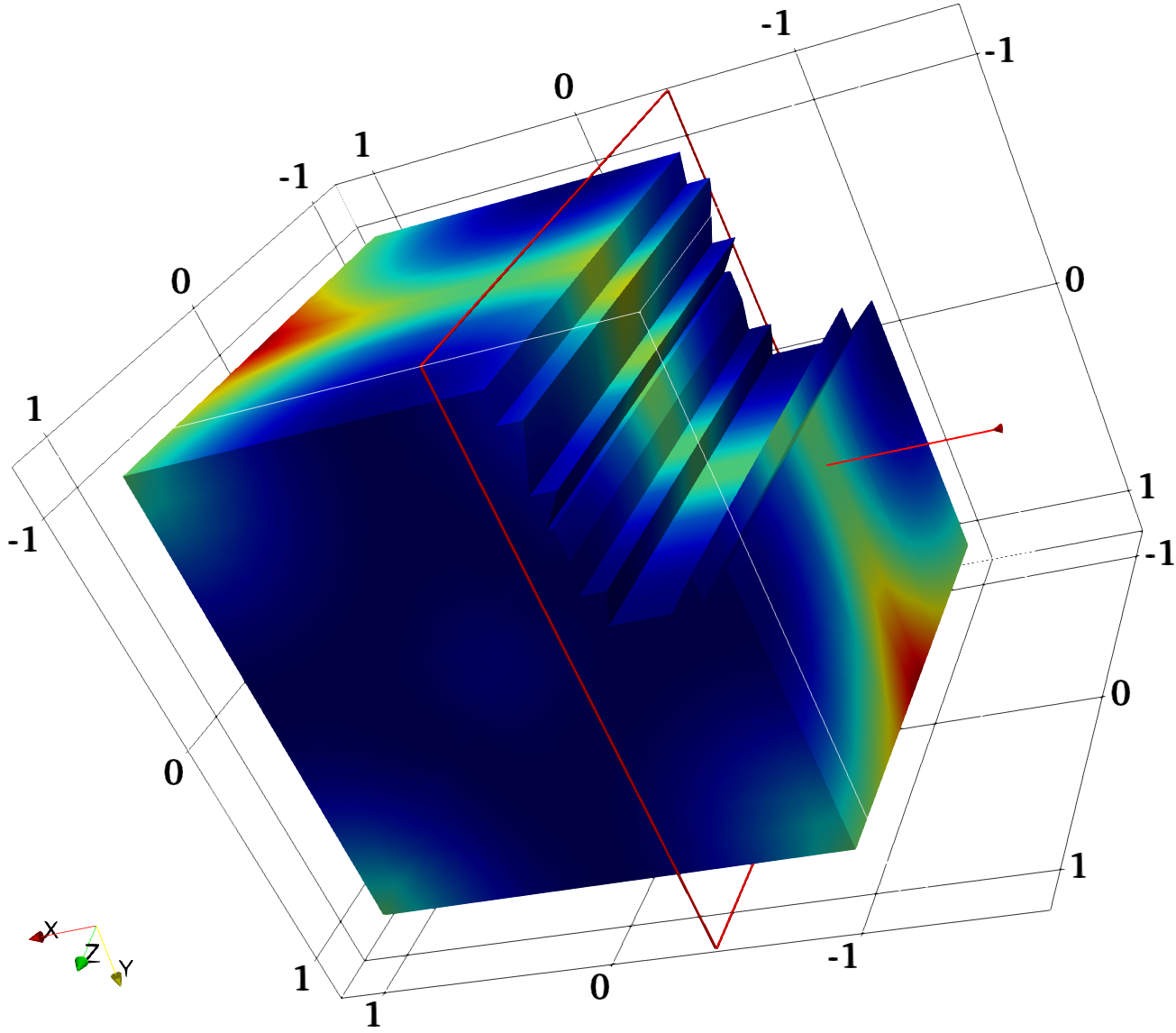}\\
\includegraphics[width=0.26\textwidth]{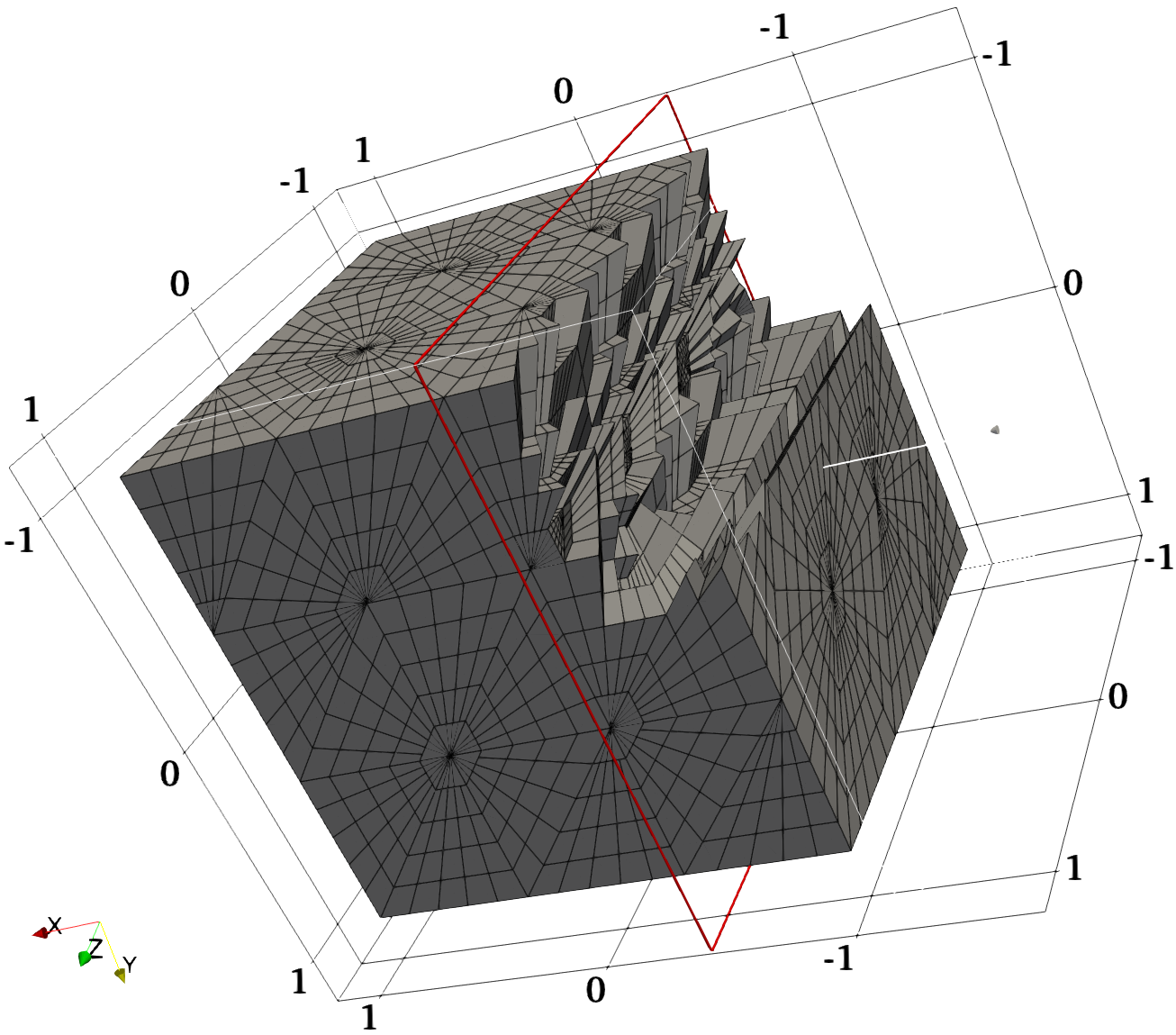}
\includegraphics[width=0.31\textwidth]{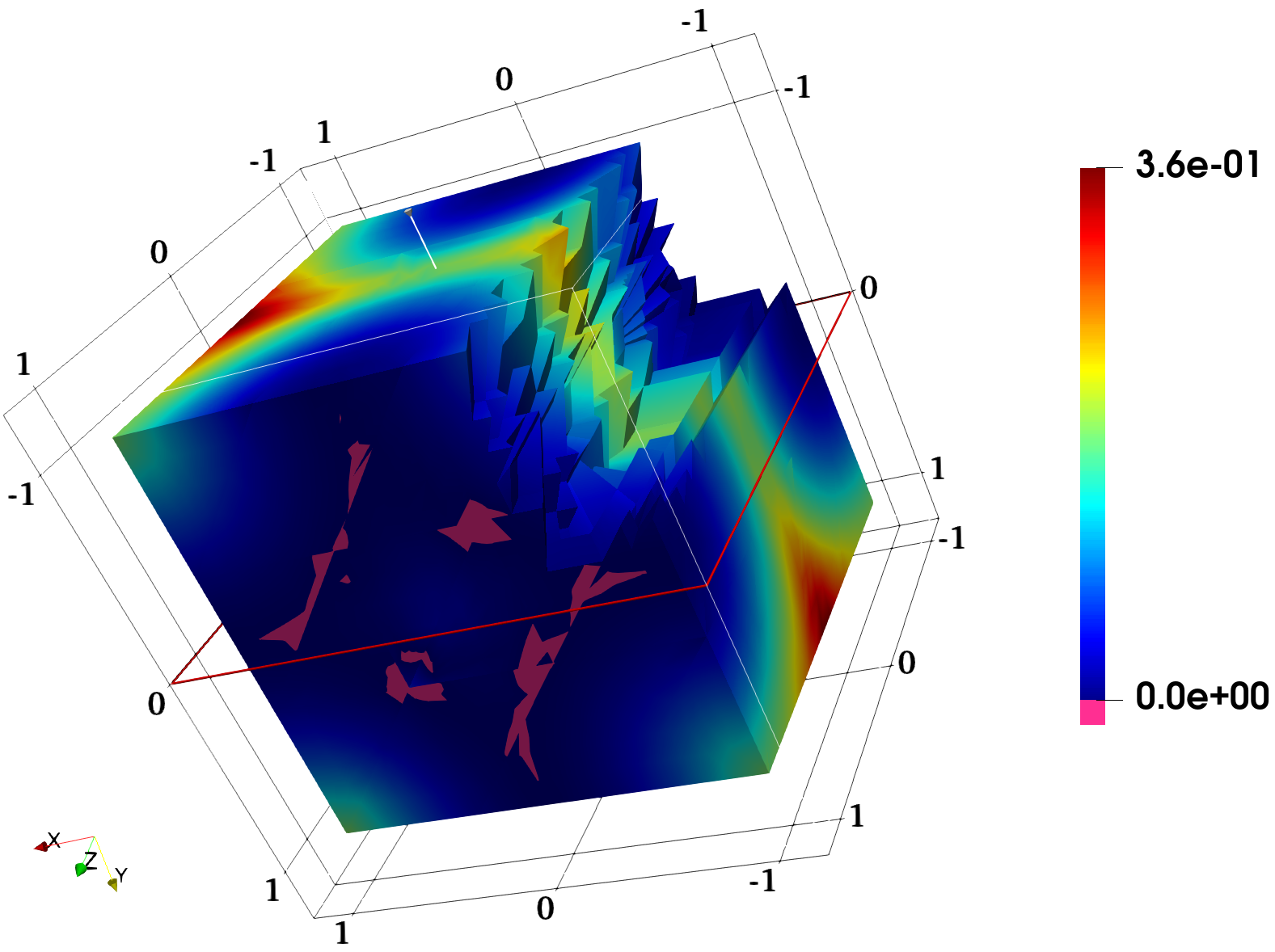}
\includegraphics[width=0.26\textwidth]{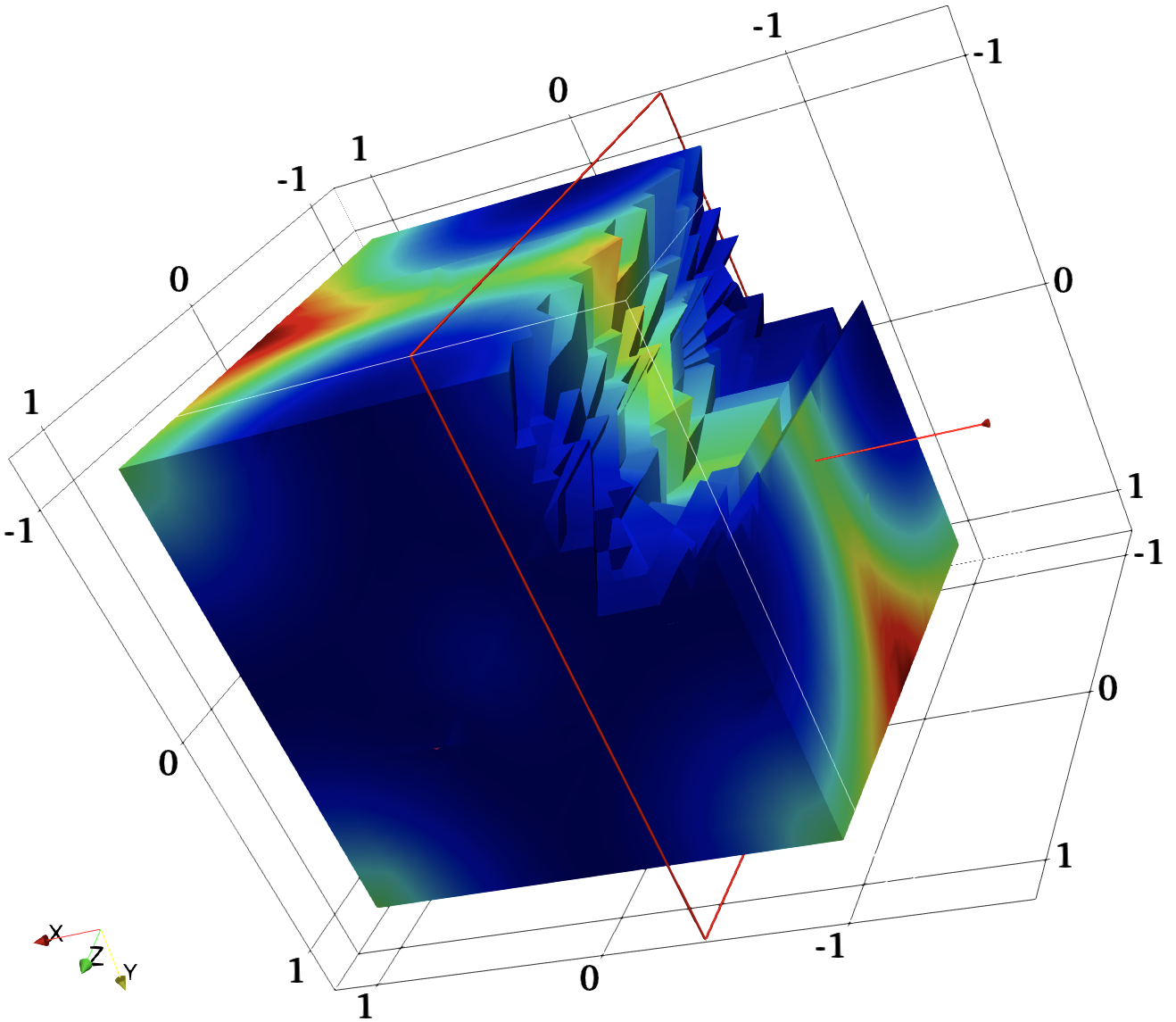}
\caption{A summary of the advection experiment using polynomial order $ = 4$ and {$\Delta t =1e-3$} on the 3D cube domains meshed using different canonical elements individually.
\textit{Left to right}: The mesh profile, the unfiltered solution, and the filtered solution at timestep = 1000, respectively.
\textit{Top row}: Cube mesh with 27 hexahedra, \textit{ Middle row}:  Cube mesh with 54 prisms,
\textit{Bottom row}:  Cube mesh with 162 tetrahedra.
}
 \label{fig:3dmeshall}
\end{figure}

\begin{equation}\label{eq:f3adv}
f(x,y,z, t = 0) = 0.2 \Big( [1-\sqrt{(x^2 + y^2)}]^2 +z^2 \Big)
\end{equation}

\begin{figure}[!b]
\centering
\includegraphics[width=0.42\textwidth]{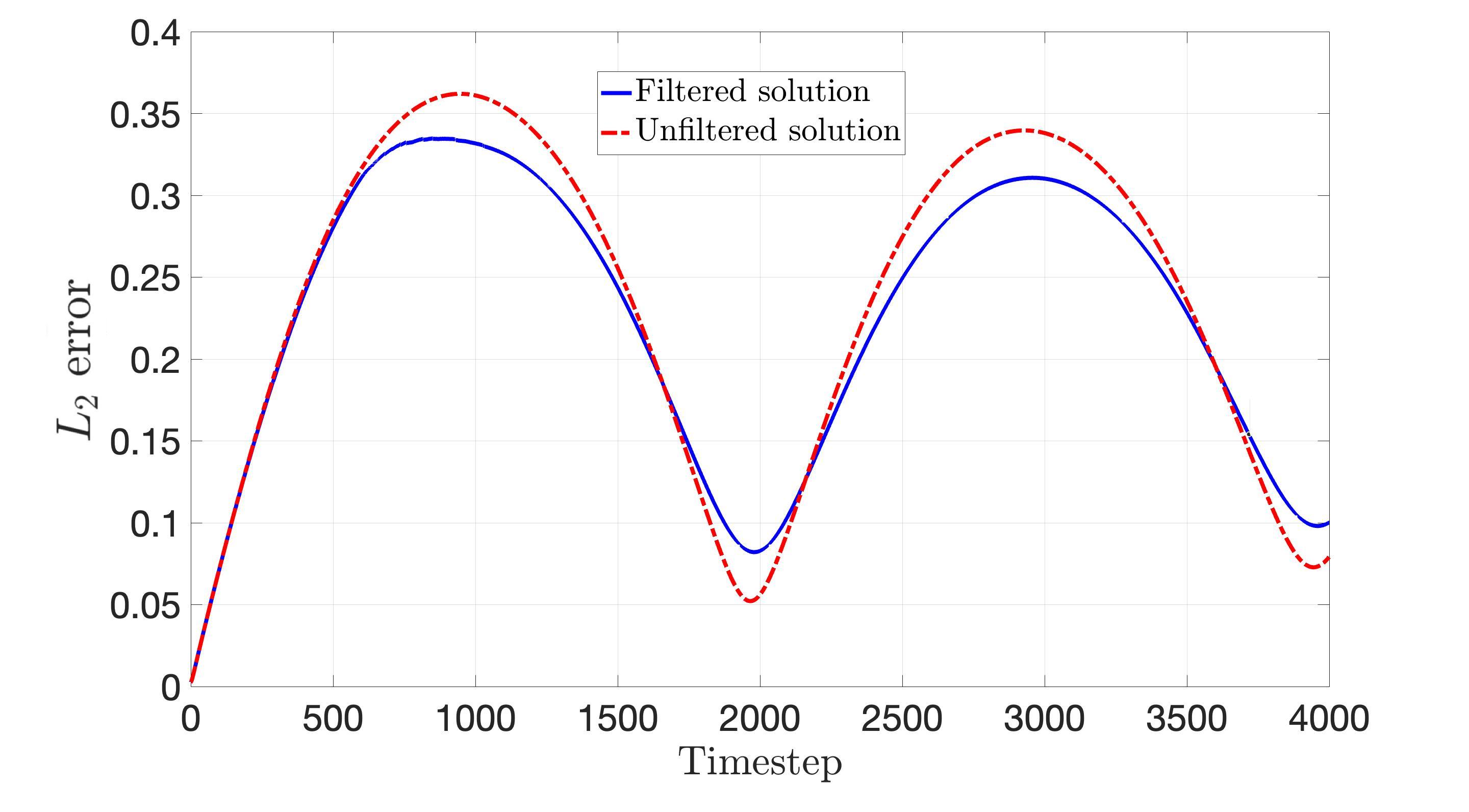}
\includegraphics[width=0.39\textwidth]{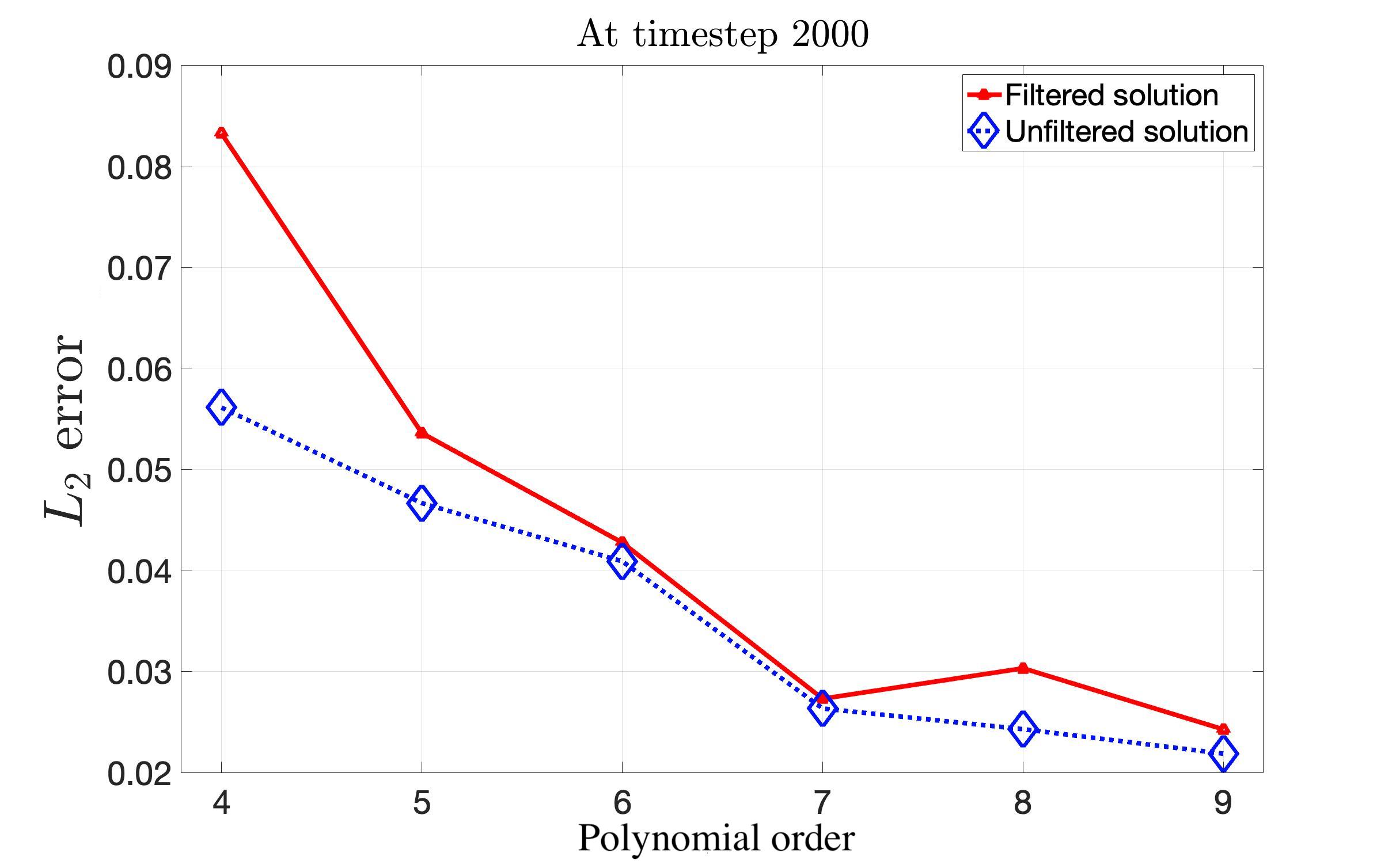}
\caption{\textit{Left:} Per timestep $L_2$ error for advection problem for polynomial order$ = 4$ on a cube mesh of 24 tetrahedral elements. The {timestep = $1e-3$}, and the total steps $ = 4000$. At each timestep, for the $L_2$ error computation, we consider 10 quadrature points in each direction. \textit{Right:} p-convergence for filtered and unfiltered $L_2$ error at timestep 2000 for the advection problem on a cube mesh of 24 tetrahedral elements. }
  \label{fig:3dpertstet}
\end{figure}

The initial state looks like a torus inside a cube. The advection experiment is repeated on different homogeneous meshes with element types hexahedra,  tetrahedra, and prisms individually. \Cref{fig:3dmeshall} shows the results before and after the filter application at a particular timestep for all these meshes.\\

The analysis of the filtering effects on the accuracy and convergence of the advection process, especially on the meshes with hexahedron and prism elements, reveals that the difference between the $L_2$ errors is too small to warrant a convergence comparison. 
However, for the mesh with tetrahedral elements, the $L_2$ errors perceptibly vary at each timestep, as shown by a longer run of the same experiment (total 4000 timesteps) in \Cref{fig:3dpertstet}.

\subsection{Filter application to dG FEM implementation of the PAC model}\label{sec:results3}

A prominent scenario in which the problem of structure-preservation is encountered is in the mathematical-biology domain. A detailed mathematical model for platelet aggregation and blood coagulation (PAC) by \cite{leiderman2011grow} is considered for this experiment.  The model describes the process of evolution, interaction, and decay of the chemical species involved in the process of thrombosis. \\

Although the details of the individual species in the PAC process are beyond the scope of this paper,  \cite{leiderman2011grow} has a detailed explanation of the nature and the evolution of the chemical species.  The model can be summarized as a comprehensive collection of ODEs and PDEs that track all the chemical species in various phases during the process. The results in \cite{leiderman2011grow} use a finite-difference method with a specified limiter to truncate nonpositive concentration values to zero. In order to implement this model using the FEM without loss of structure, an alternative approach to truncation, such as the one presented in this paper, is needed. For our experiments, we use the structure-preserving filter instead of a truncation limiter. The focus is limited to the system of equations to track the chemical species that advect, diffuse, and react. One such species is fluid-phase thrombin (FP $e_2$), an essential component of the thrombosis process. For this section, we track the evolution of thrombin by advection, diffusion, and decay in the fluid-phase. For the detailed results of the evolution of all the chemical species, including thrombin, under various circumstances in the original model, refer to \cite{leiderman2011grow}.\\


We set up a version of the PAC model that solves the species evolution problem by a combination of advection, diffusion, and reaction (ADR) PDEs. Unlike the original model, the velocity $\vu$ of the fluid medium is not reported by a modified Navier-Stokes solver. Instead, the velocity is constant $\vu = [u_x, u_y] = [5, 0]$. To solve this system of PDEs on a sample blood vessel represented by a rectangular domain $\Omega = [0,300]\times[0,20]$, we employ the dG FEM.  The domain is tessellated using 2048 quadrilaterals. The bottom wall has an injury site spanning from $x = 20$ to $x = 60$. We use polynomial order $=4$  and timestep of {$\Delta t = 1e{-2}$} to solve the problem. \\

 The total number of species tracked in our implementation of the model is 56. Since the focus is on the behavior of thrombin, we present the PDE describing the evolution of fluid-phase thrombin (FP $e_2$) \Cref{eq:thrombosis}. For details on the other chemical species involved in \Cref{eq:thrombosis}, refer to \cite{leiderman2011grow}. The boundary conditions vary depending on the chemical being tracked. For FP $e_2$, we have a Dirichlet zero boundary on the left, and a no-flux boundary on the right and top. The bottom boundary has two kinds of conditions: a robin boundary condition on the injury site and Dirichlet zero everywhere else.  
\begin{equation}\label{eq:thrombosis}
\begin{split}
\frac{\partial \ve_{2}}{\partial t} &=  \underbrace{-\vu \cdot \nabla e_{2} + \nabla \cdot (D\nabla e_{2})}_\text{Transport by Advection Diffusion}\\
& + \underbrace{ k_{e_{2}}^{on}e_2(N_2P^{b,a}+N_2P^{se,a} - z_2^{mtot} - e_2^{mtot}) + k_{e_2}^{off}e_2^m}_\text{Binding to platelet receptor}\\
& + \underbrace{ (k_{z_5:e_{2}}^{cat} + k_{
z_5:e_2}^-) [Z_5:E_2] - k_{z_5:e_2}^+z_5e_2}_\text{Activation of V}\\
& + \underbrace{ (k_{z_7:e_{2}}^{cat}  + k_{
z_7:e_2}^-) [Z_7:E_2] - k_{z_7:e_2}^+z_7e_2}_\text{Activation of VII}\\
& + \underbrace{(k_{z_8:e_{2}}^{cat}  + k_{
z_8:e_2}^-) [Z_8:E_2] - k_{z_8:e_2}^+z_8e_2}_\text{Activation of VIII}
\end{split}
\end{equation}
\begin{figure}[h]
  \begin{center}
\includegraphics[width=0.8\textwidth]{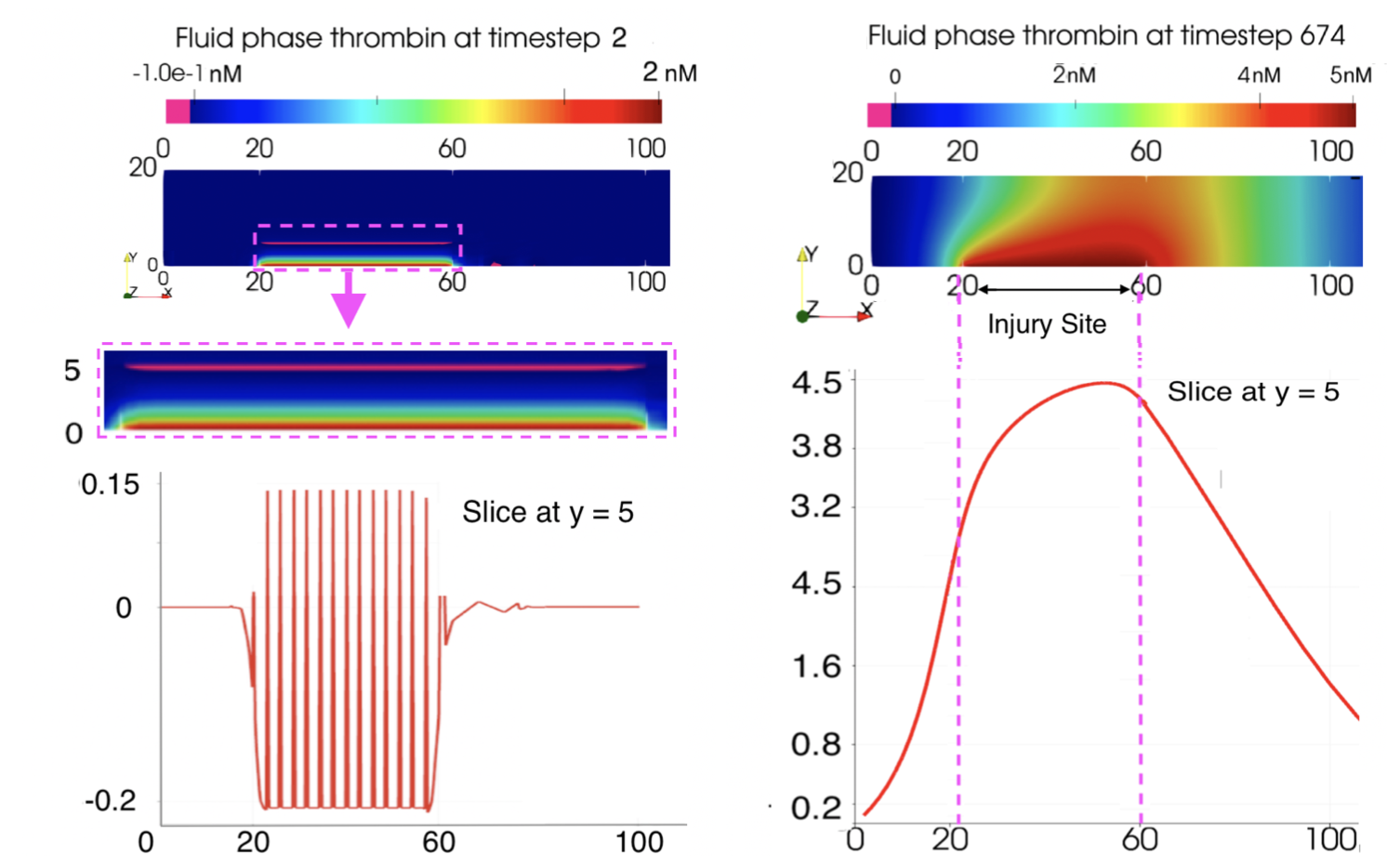}
\caption{\textit{Top left:} The unconstrained simulation solution becomes invalid, as shown in the inset at \textbf{timestep 2}, which causes a blowup at timestep 87, leading to the failure of the simulation.
\textit{Top right:} The constrained simulation runs for\textbf{ 674 timesteps} until the depletion of the chemicals causes the model to terminate naturally.
\textit{Bottom left:} Profile of fluid phase thrombin in unconstrained simulation at timestep 2 on a domain slice at $y=5$. \textit{Bottom right:} Profile of fluid phase thrombin in constrained simulation at timestep 674 on a domain slice $y=5$}
 \label{fig:PAC2con}
\end{center}
    \end{figure}

\Cref{fig:PAC2con} shows filtered and unfiltered runs of our implementation of this model and the resulting concentration profile of FP $e_2$. At timestep = 2, we get the invalid (negative) values around the injury site because of the sharp concentration changes (refer to the slice at $y=5$ in the bottom left part of \Cref{fig:PAC2con}). Therefore, the fluid phase thrombin concentration does not adhere to the desired structure (positivity), and the simulation does not succeed. If this experiment is allowed to continue, the effect of invalid values will compound over time, rendering the simulation results nonphysical and useless. At timestep 38, the simulation starts becoming unstable because of structural discrepancies.  If it is allowed to run further, the accumulation and propagation of negative values results in the simulation blow-up at timestep 87, at which point the code reports invalid values such as \textit{NaN} (datatype: Not a Number). \\

When the filter is applied, it runs for significantly more timesteps (total 674) and terminates naturally. The right side of  \Cref{fig:PAC2con} shows the state of FP $e_2$ at the depletion point of the chemicals involved.

\section{Conclusions}

We present a formalism that solves the problem of structure-preservation for PDE solutions in 2D and 3D. The construction and design of a postprocessing structure-preserving filter are detailed and applied to multidimensional dG FEM solutions to different PDEs. A geometric interpretation of the mathematical foundation behind the filter is presented, followed by an algorithm to apply the proposed filter in a timestepping PDE framework. At the core of the filter lies the expensive requirement for global minimization on a weighted distance function that corresponds to the objective we optimize. We employ gradient descent with backtracking line search for the minimization to reduce the cost of the filtering routine. To this end, we detail an investigative procedure to reduce the minimization cost by precomputing specific parameter values for the gradient descent approach. Using numerical examples, we compare the convergence rates with and without the application of the filter for different problem sizes and domain structures. The percentage increase in time taken by the simulation is computed to understand the cost of the filter and it is observed that the cost scales with the order and size of the problem. In the end, using a filtered solution to a mathematical-biology problem of platelet aggregation and blood coagulation, we provide evidence of the proposed method's efficacy and utility.  {One future direction for investigation could attempt to understand when inter-element flux preservation is beneficial in these types of numerical simulations. For example, comparing two filtered solutions, one with and one without flux preservation, could form the basis for more experiments.}
\label{sec:conclusions}

\section*{Acknowledgments}

V. Zala and  R.M. Kirby acknowledge that their part of this research was sponsored by ARL under cooperative agreement number W911NF-12-2-0023. The views and conclusions contained in this document are those of the authors and should not be interpreted as representing the official policies, either expressed or implied, of ARL or the U.S. Government. The U.S. Government is authorized to reproduce and distribute reprints for Government purposes notwithstanding any copyright notation herein.  A. Narayan was partially supported by NSF DMS-1848508. This material is based upon work supported by both the National Science Foundation under Grant No. DMS-1439786  and the Simons Foundation Institute Grant Award ID 507536 while A. Narayan was in residence at the Institute for Computational and Experimental Research in Mathematics in Providence, RI, during the Spring 2020 semester.

\bibliographystyle{siamplain}
\bibliography{references}
\clearpage
\setcounter{section}{0}
\renewcommand{\thesection}{Appendix \Alph{section}}

\section{}\label{sec:apdx}

\setcounter{figure}{0}
\setcounter{algorithm}{0}
\renewcommand{\thealgorithm}{A.\arabic{algorithm}}
The parameter selection algorithm for backtracking line search used in GD-based minimization. \Cref{eq:fns} and \Cref{eq:fns3} define the set of functions used in the selection process.
\begin{algorithm}[H]
 \caption{Experiment to determine ideal values of $c$ and $\gamma$ given a set of functions $f_i$, for $i = 0,1,\cdots,n$}
\label{alg:candgam}
\begin{algorithmic}[1]
\STATE{$G \leftarrow {P}$ samples in $(0,1)$ }
\FOR{each $f_i$ $\forall$ $i =0,1,
\cdots,n$}
\STATE{$H\leftarrow \{\}$ }
\FOR{each $\gamma_j$ in G}
\FOR {each $c_k$ in G}
\STATE{ $(niter, err) \leftarrow $ {GD\_linesearch} ($f_i, c_k, \gamma_j$)}
 \STATE{$H \leftarrow \Big\{H; \{c_k, \gamma_j, niter, err \}\Big\}$ }
 \ENDFOR
\ENDFOR
\ENDFOR
\STATE{$H_1 \leftarrow$ select ranges of $c$ and $\gamma$ with least $niter$ from $H$}
\STATE{$H_2 \leftarrow$ select ranges of $c$ and $\gamma$ with least $err$ from $H_1$}
\RETURN {Ranges for $c$ and $\gamma$ from $H_2$}
\end{algorithmic}
\end{algorithm}

\section{}\label{sec:apdxb}

\renewcommand{\thefigure}{B.\arabic{figure}}
Additional results for experiment described in \Cref{sub:levequetests} for polynomial orders 3 and 7 are attached below:\\

\begin{figure}[!h]
\centering
\includegraphics[width=0.44\textwidth]{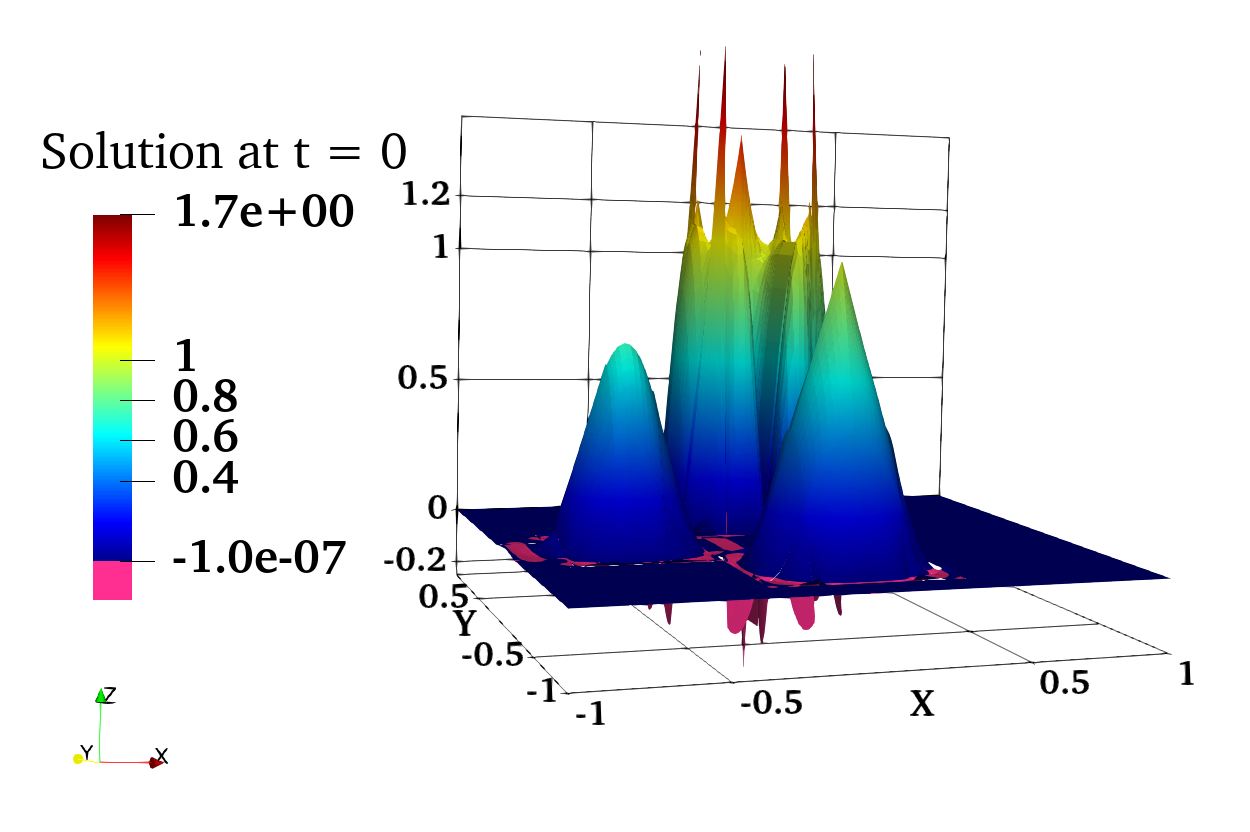}
\includegraphics[width=0.44\textwidth]{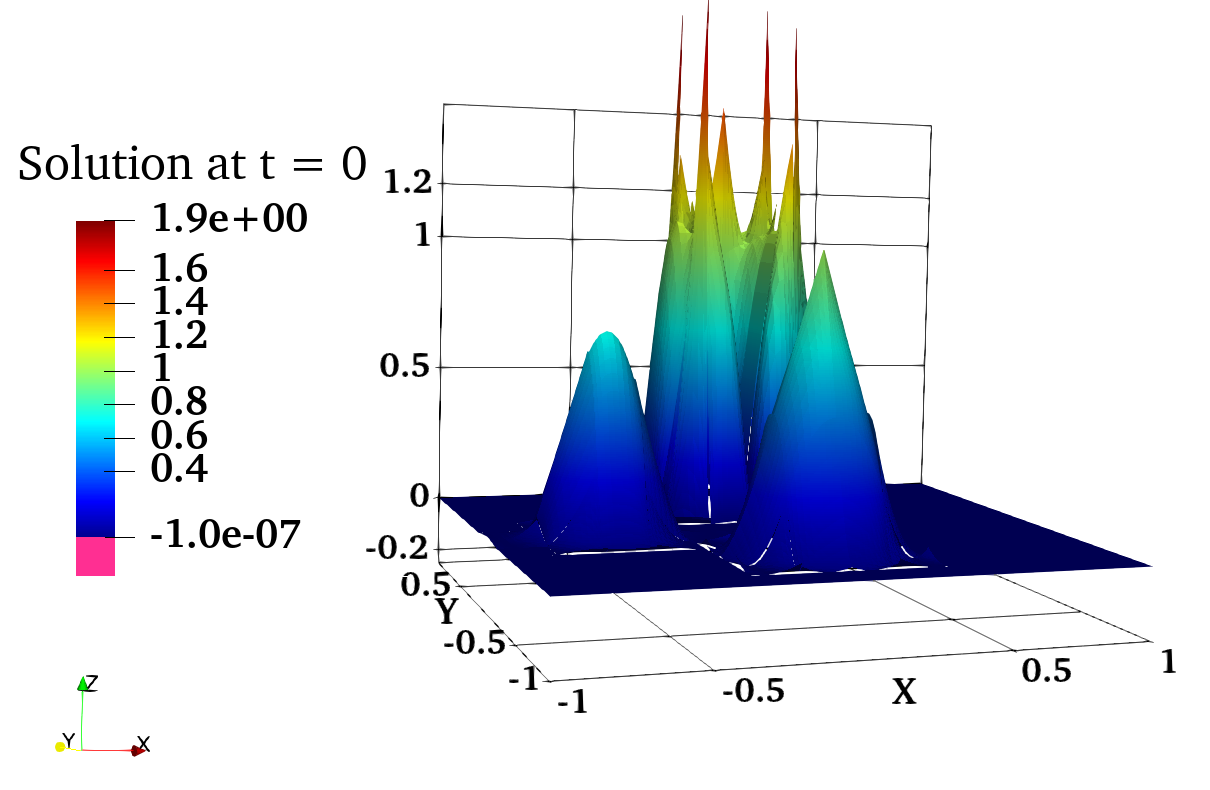}\\
\includegraphics[width=0.44\textwidth]{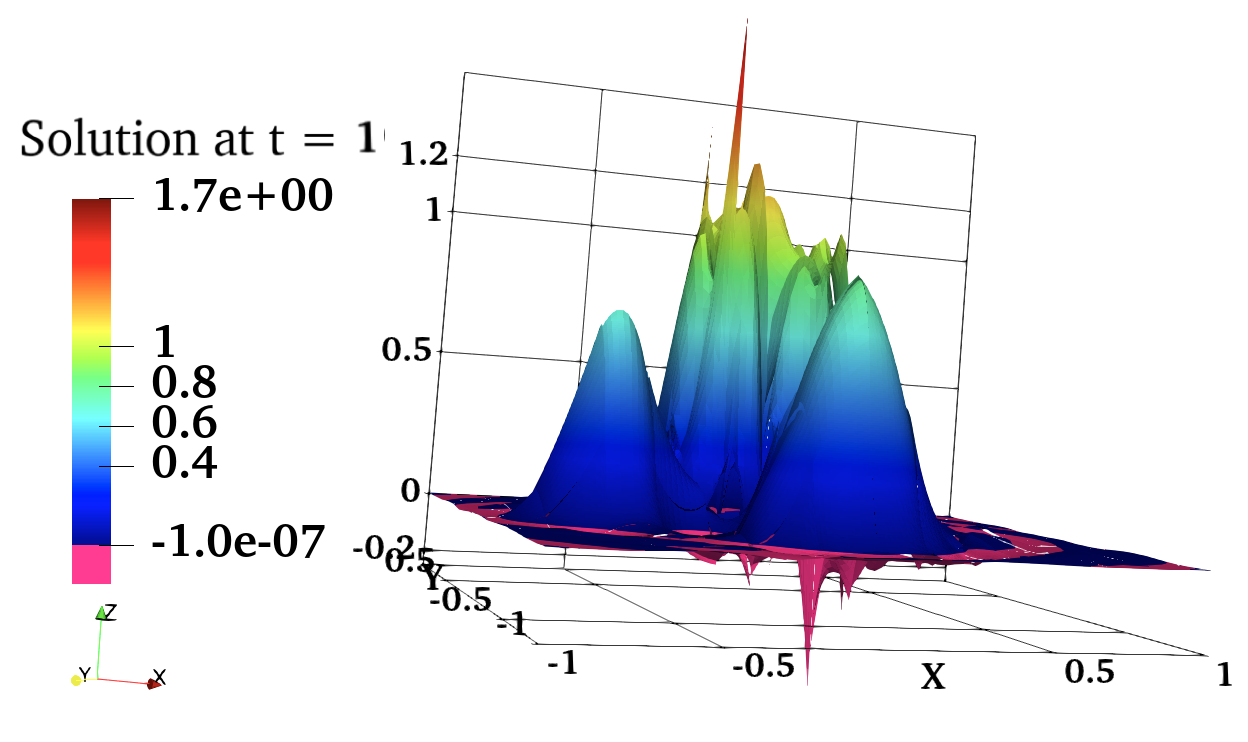}
\includegraphics[width=0.44\textwidth]{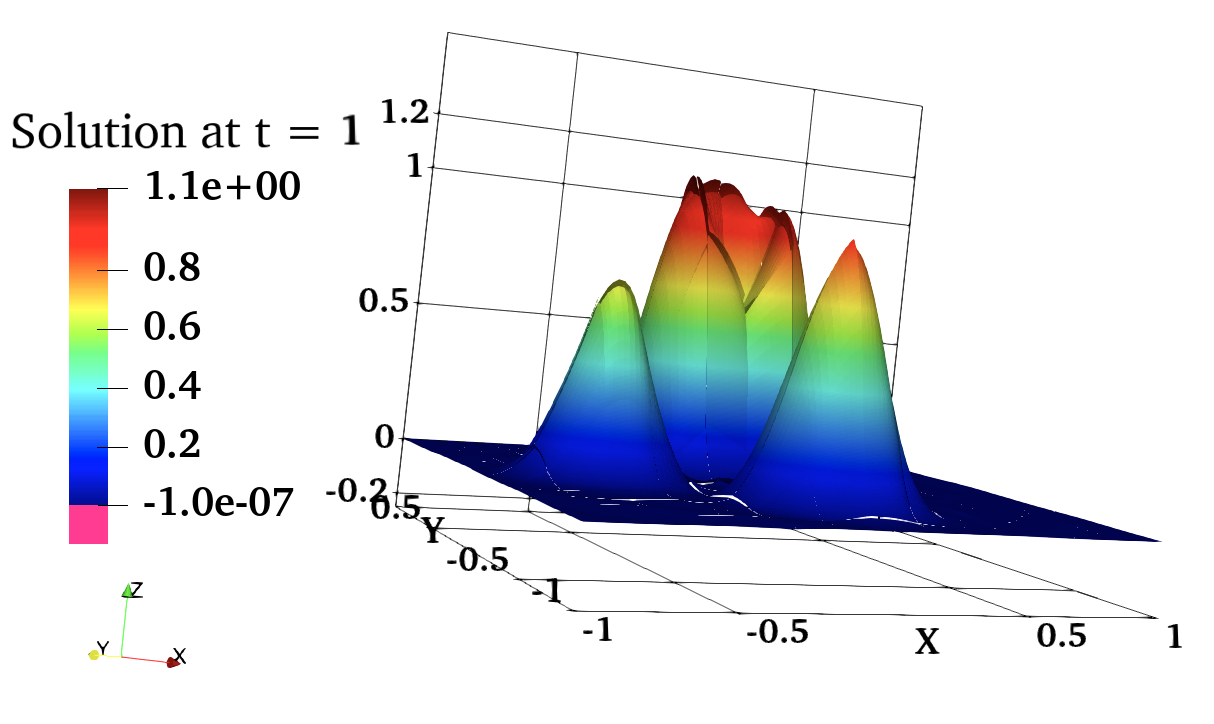}
\caption{{Snapshots at the beginning and end of advection. The highlighted region shows the values below tolerance for negativity $(10^{-7})$.  At t = 1s, the solid body finishes one rotation and returns to the original position.  Parameters for the test: Timestep = $5e-4$, polynomial order = 3. 
\textit{Row 1:} Time = 0s. \textit{Row 2:} Time = 1s. \textit{Column 1:} Unconstrained solution. \textit{Column 2:} Constrained solution.  }}
  \label{fig:newfn_apb}
\end{figure}

\begin{figure}
\centering
\includegraphics[width=0.36\textwidth]{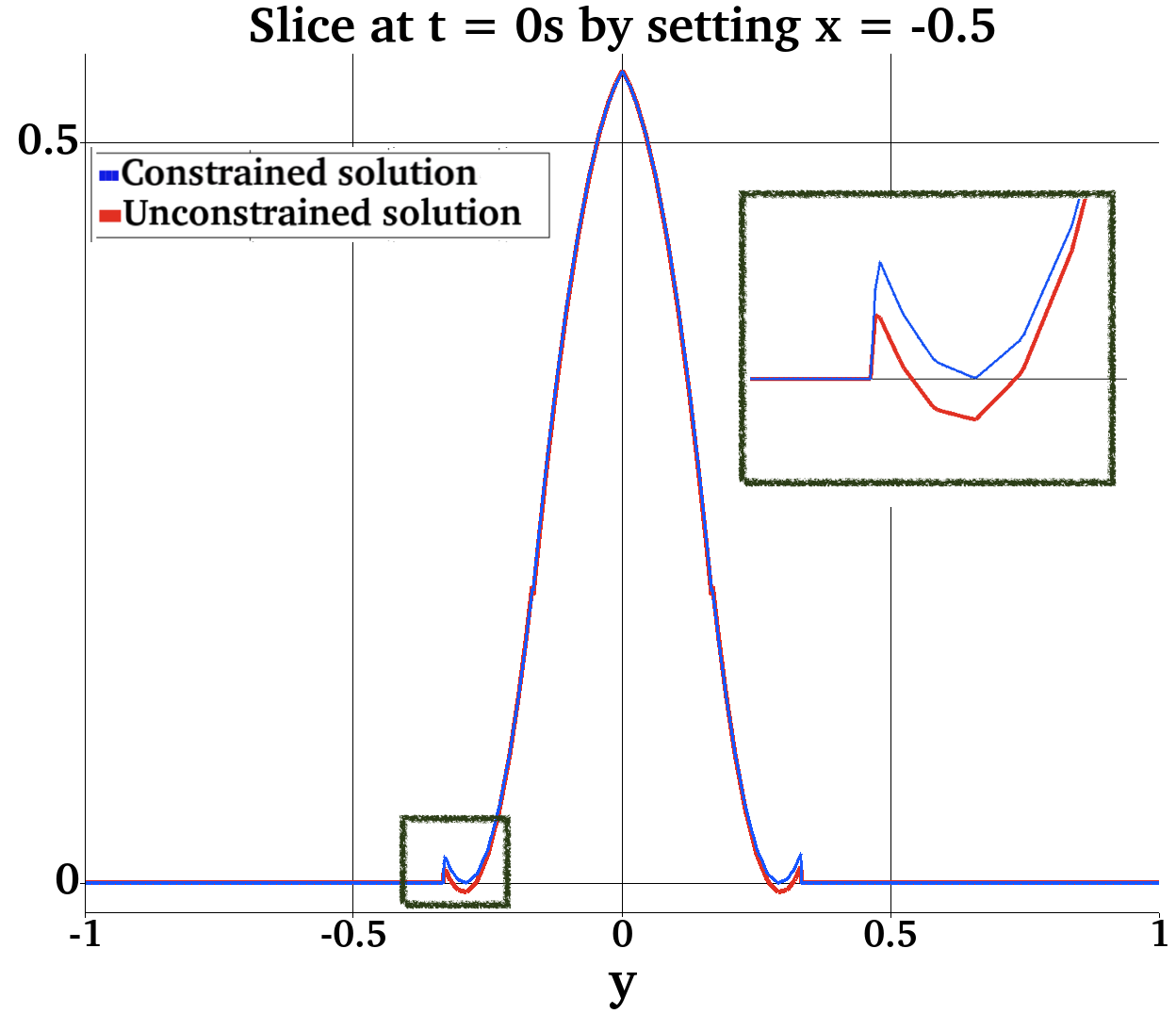}
\includegraphics[width=0.36\textwidth]{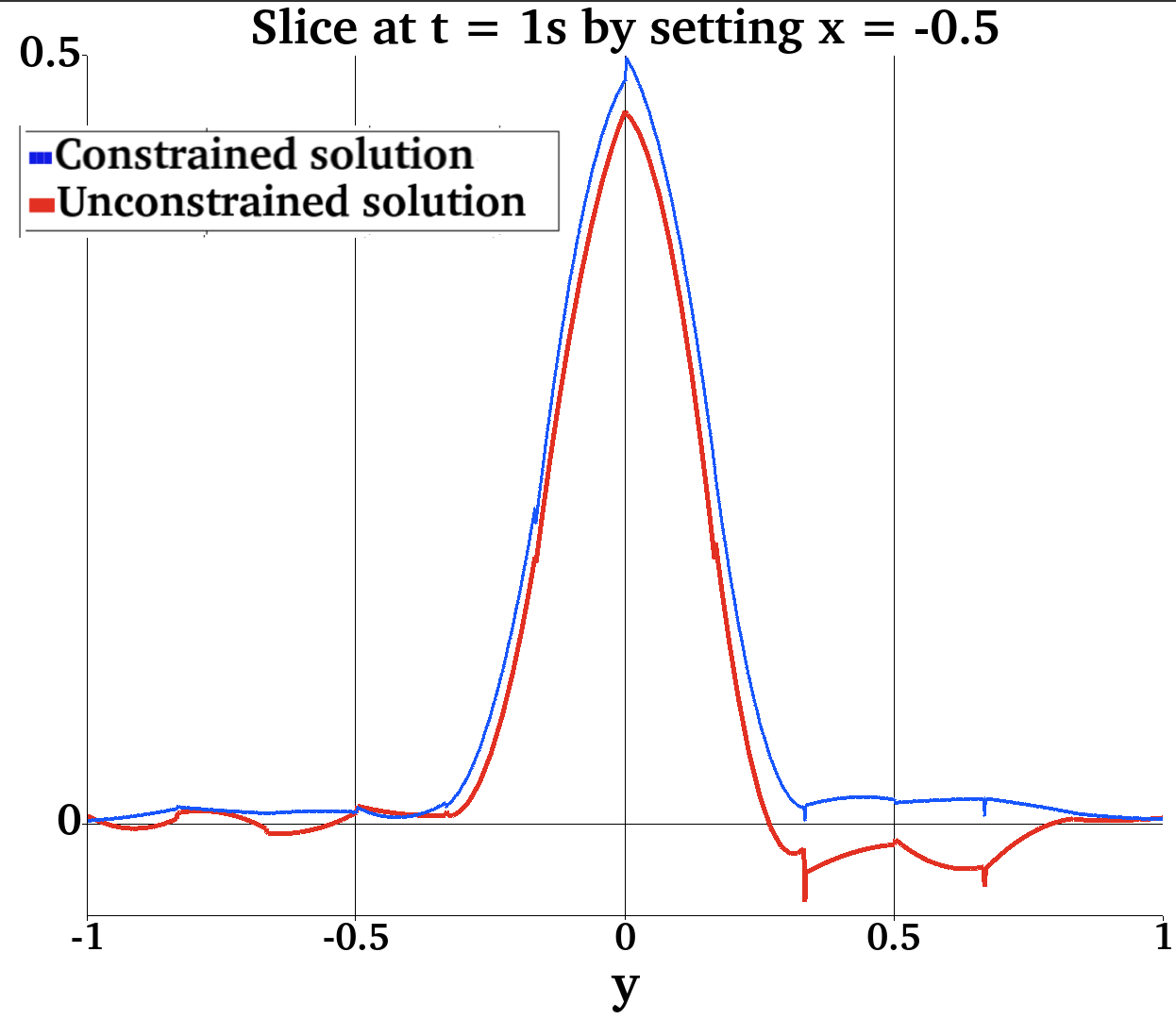}\\
\includegraphics[width=0.36\textwidth]{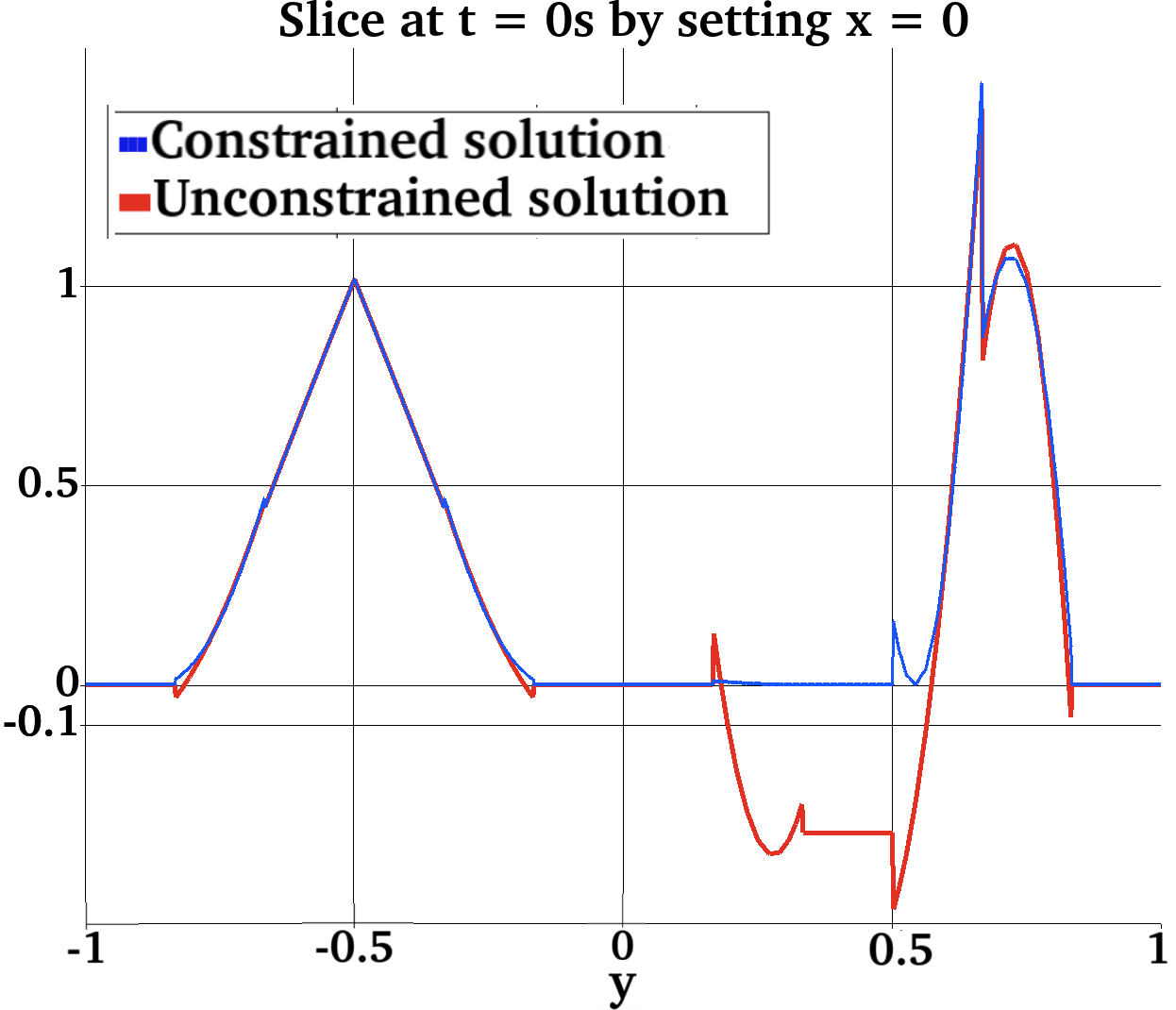}
\includegraphics[width=0.36\textwidth]{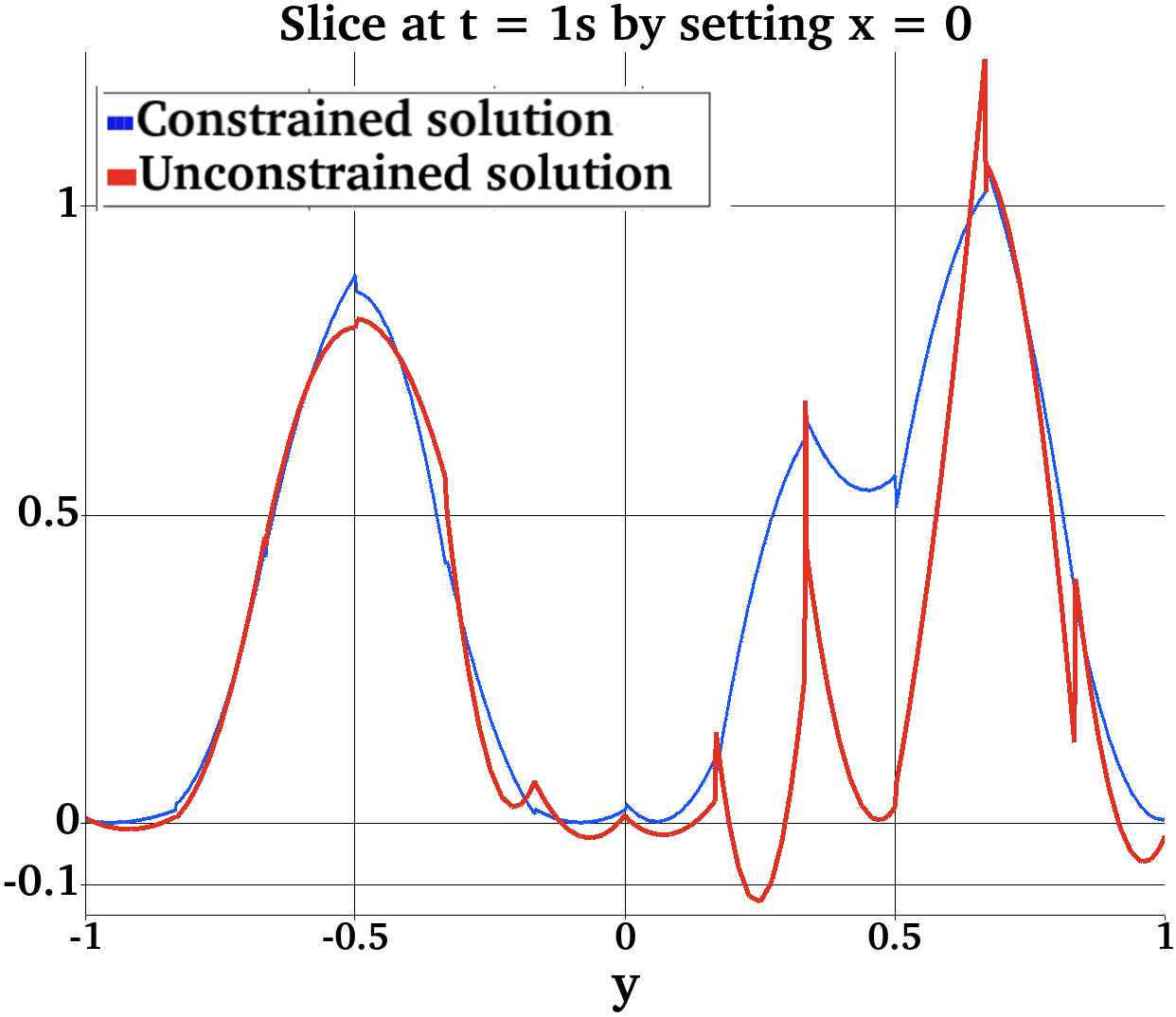}\\
\includegraphics[width=0.36\textwidth]{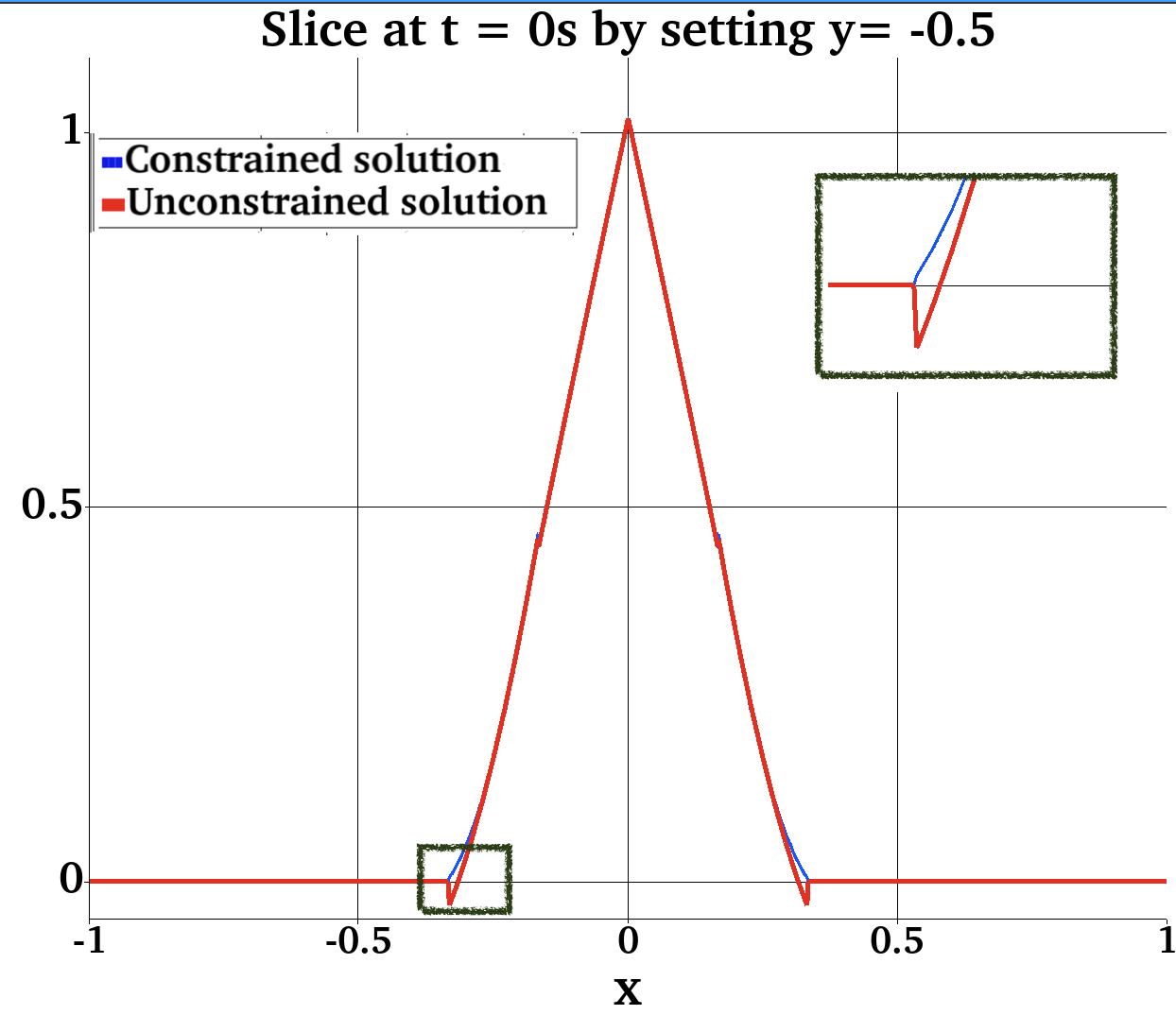}
\includegraphics[width=0.36\textwidth]{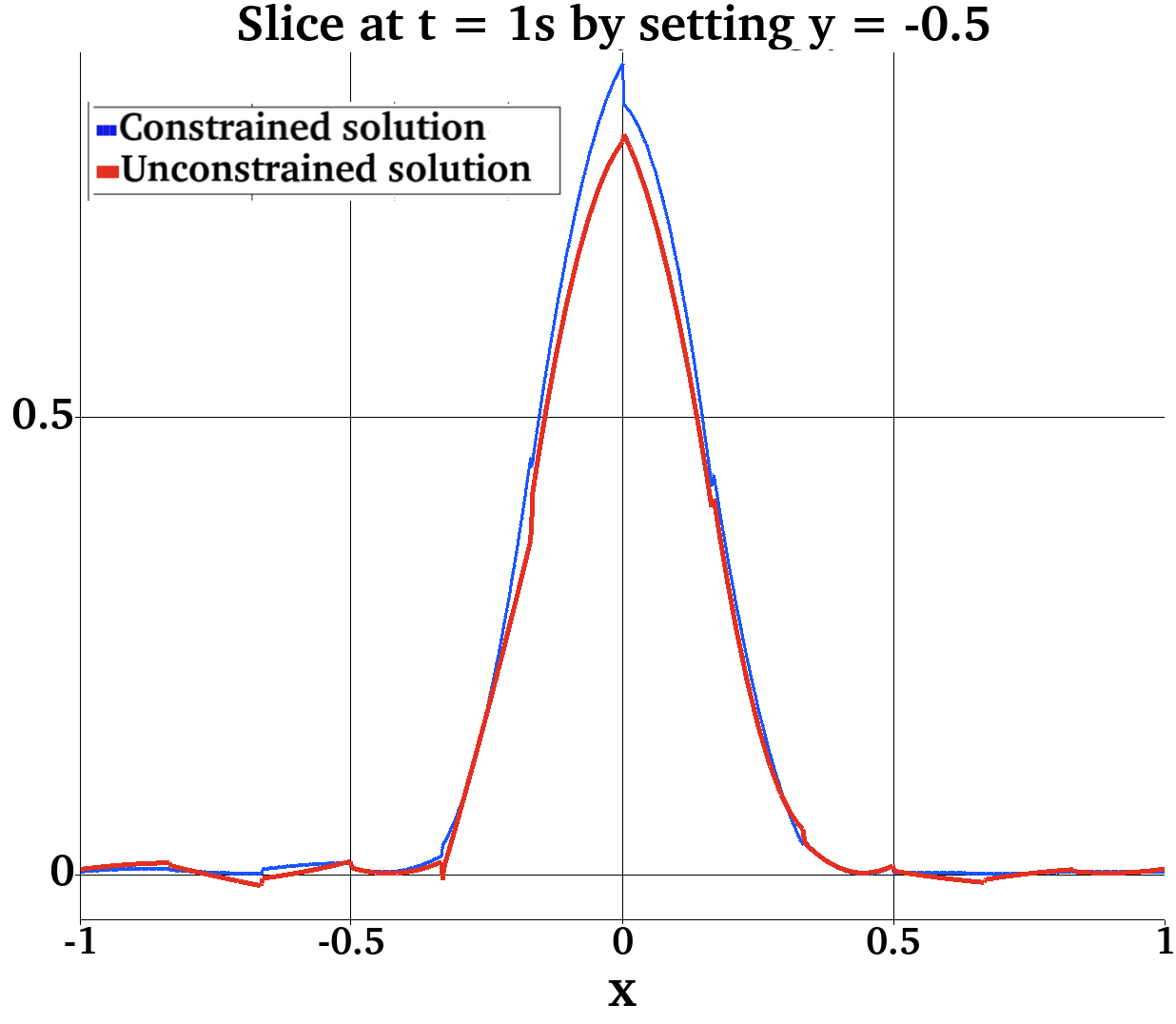}\\
\includegraphics[width=0.36\textwidth]{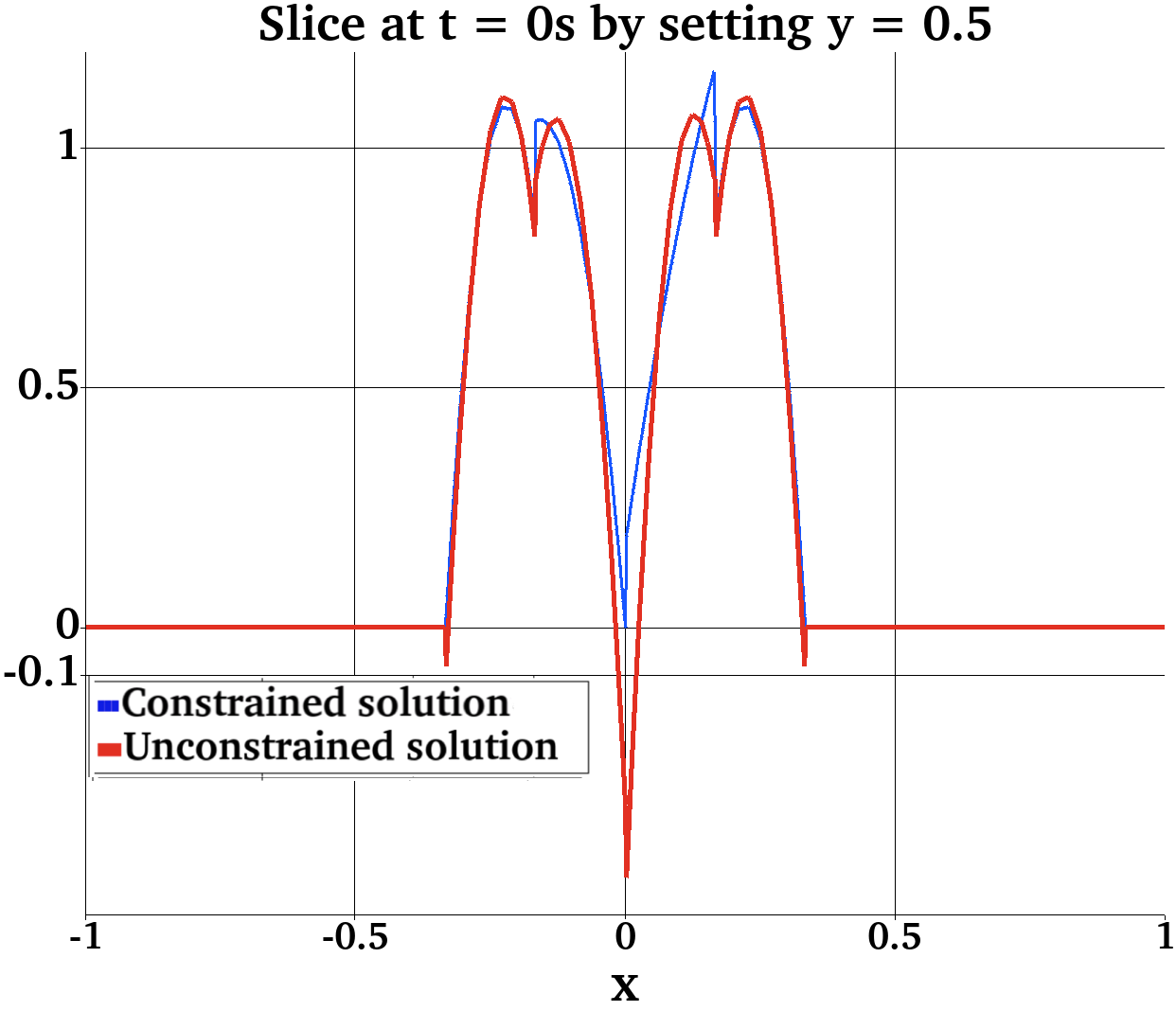}
\includegraphics[width=0.36\textwidth]{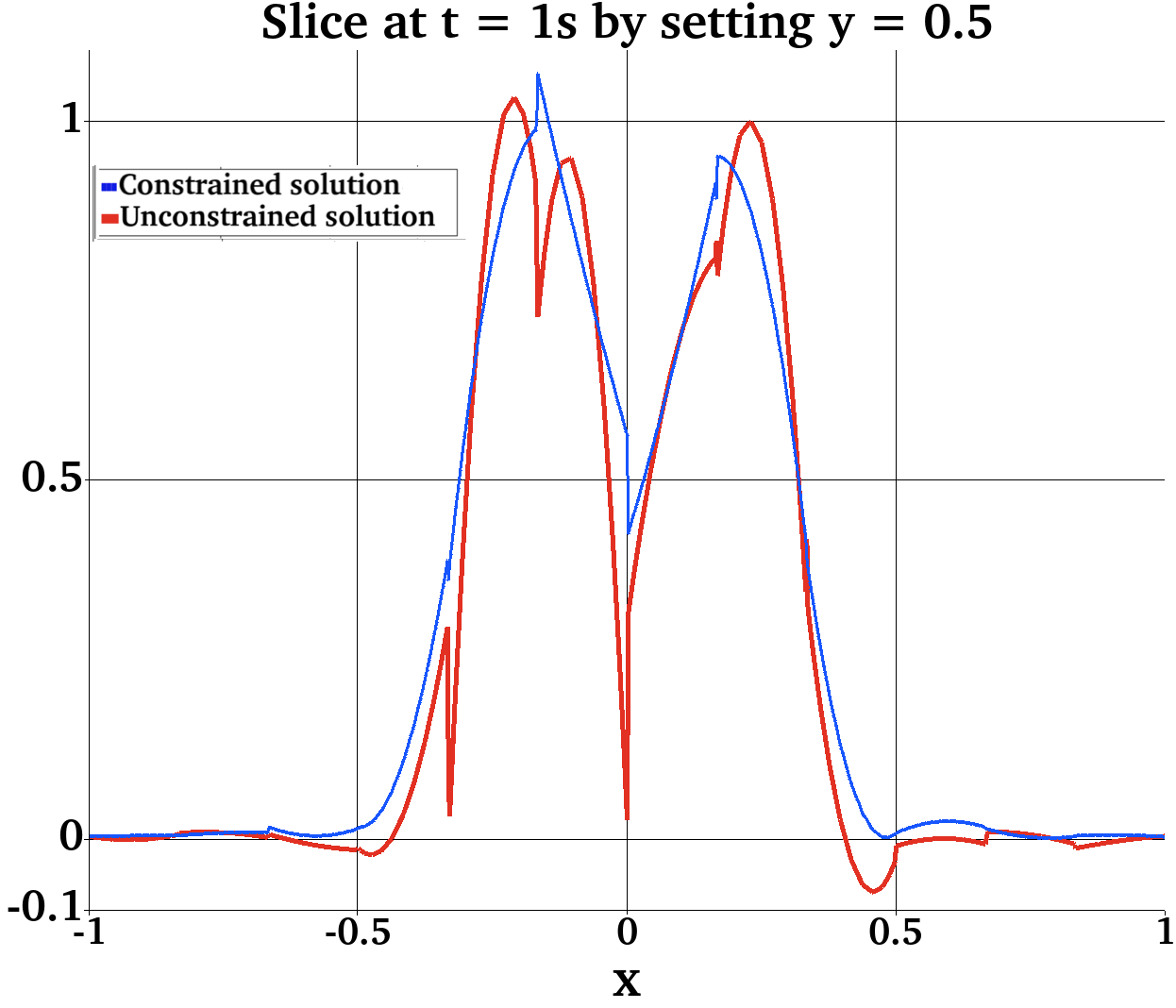}
\caption{{Snapshots at various points in time during advection. $\Delta t = 5e-4$, time period = 1 and polynomial order = 3.  Some subfigures show insets with the zoomed-in area of interest.
\textit{Column 1:} State of the system at t = 0. \textit{Column 2:} State of the system at t = 1s.  \textit{Row 1:} Slice at x = -0.5. \textit{Row 2:} Slice at x = 0. \textit{Row 3 :}Slice at y = -0.5 and \textit{Row 4: }Slice at y = 0.5}}
  \label{fig:adv3dlevequeslice1_apb}
\end{figure}

\begin{figure}[b]
\centering
\includegraphics[width=0.44\textwidth]{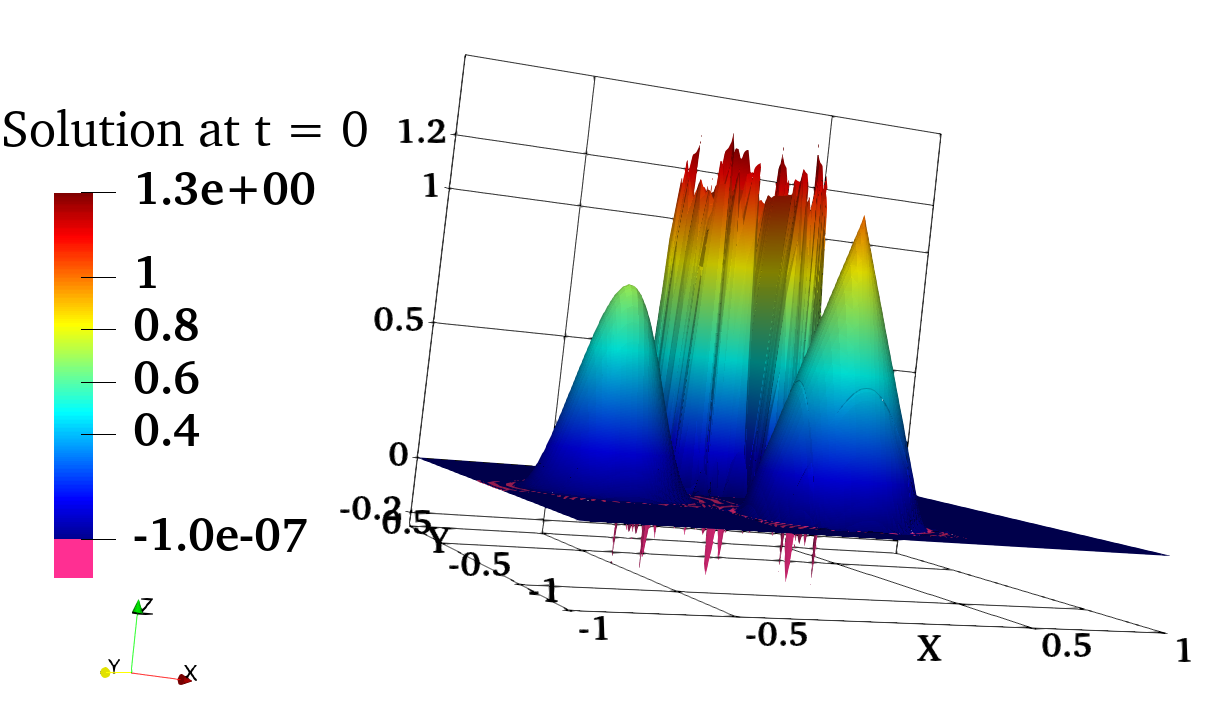}
\includegraphics[width=0.44\textwidth]{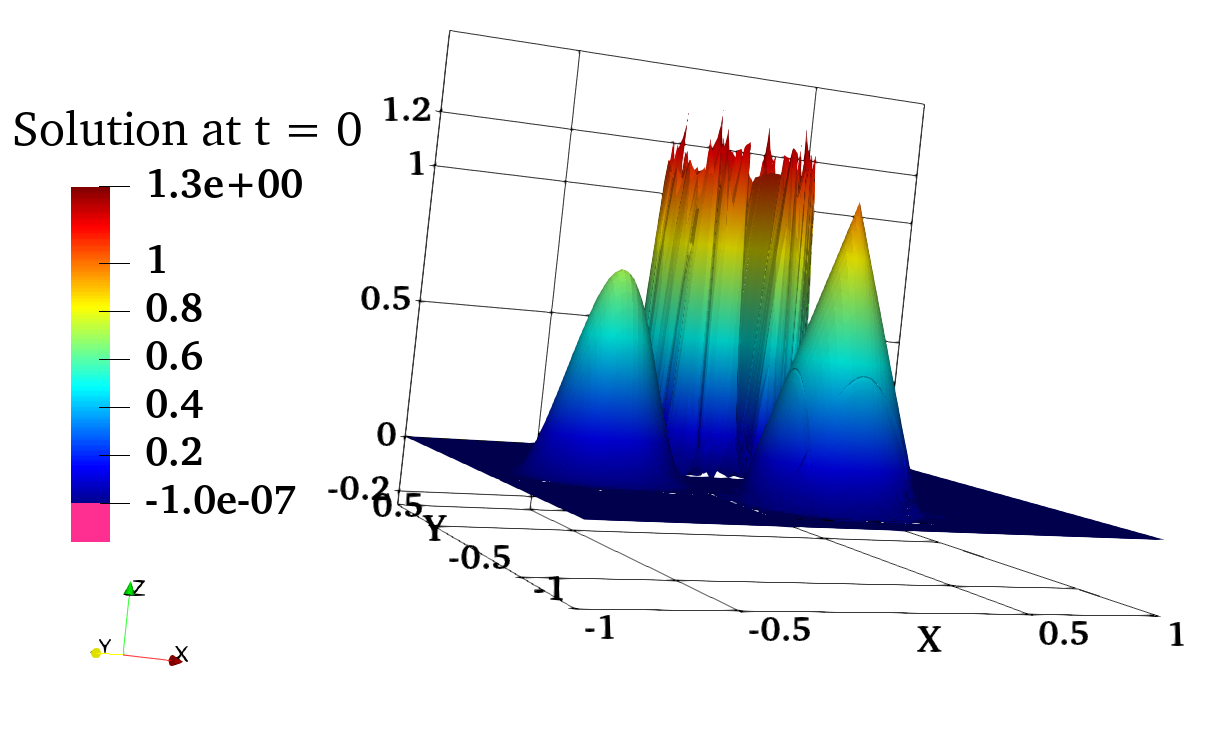}\\
\includegraphics[width=0.44\textwidth]{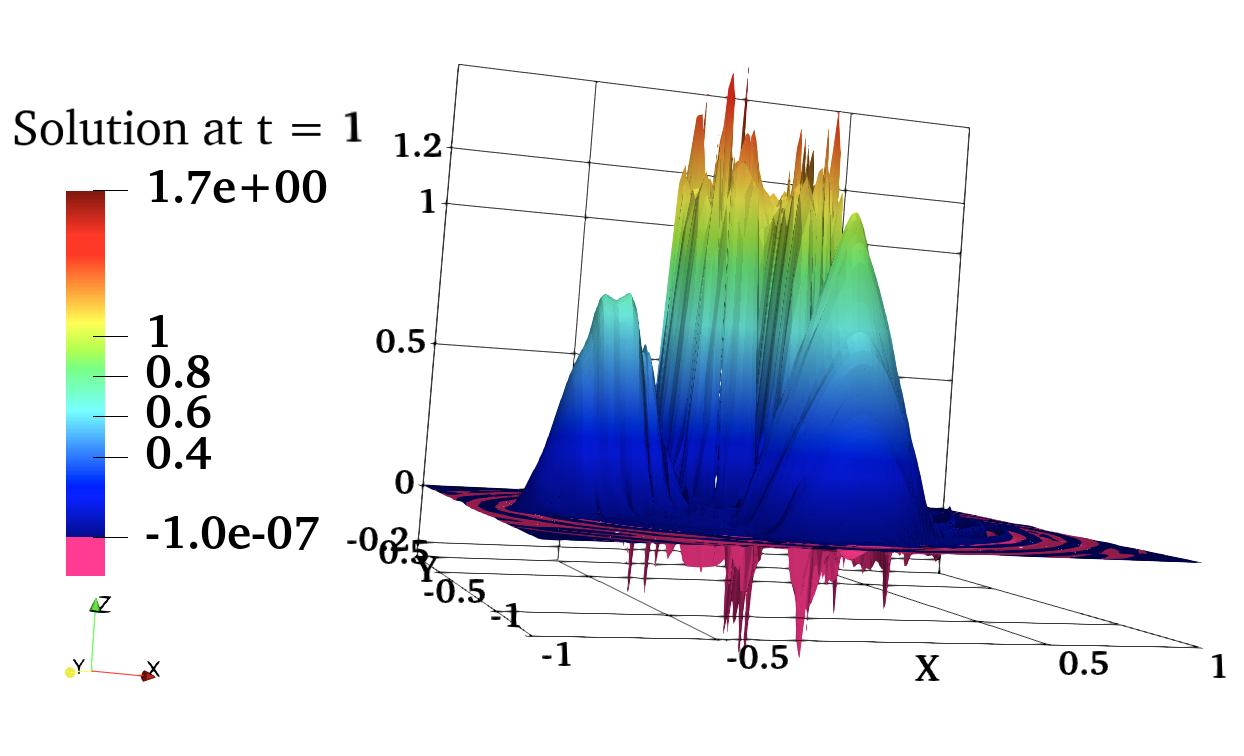}
\includegraphics[width=0.44\textwidth]{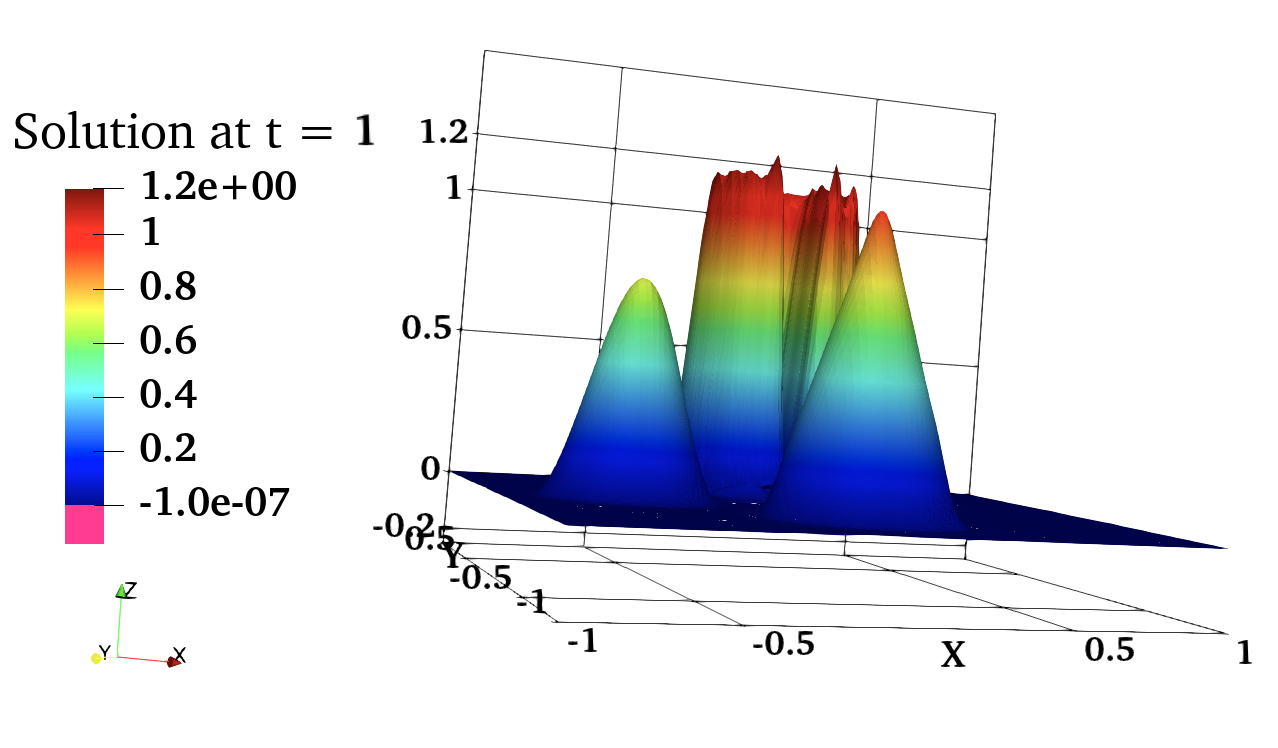}
\caption{{Snapshots at the beginning and end of advection. The highlighted region shows the values below tolerance for negativity $(10^{-7})$.  At t = 1s, the solid body finishes one rotation and returns to the original position.  Parameters for the test: Timestep = $5e-4$, polynomial order = 7. 
\textit{Row 1:} Time = 0s. \textit{Row 2:} Time = 1s. \textit{Column 1:} Unconstrained solution. \textit{Column 2:} Constrained solution.  }}
  \label{fig:newfn_apb2}
\end{figure}

\begin{figure}
\centering
\includegraphics[width=0.36\textwidth]{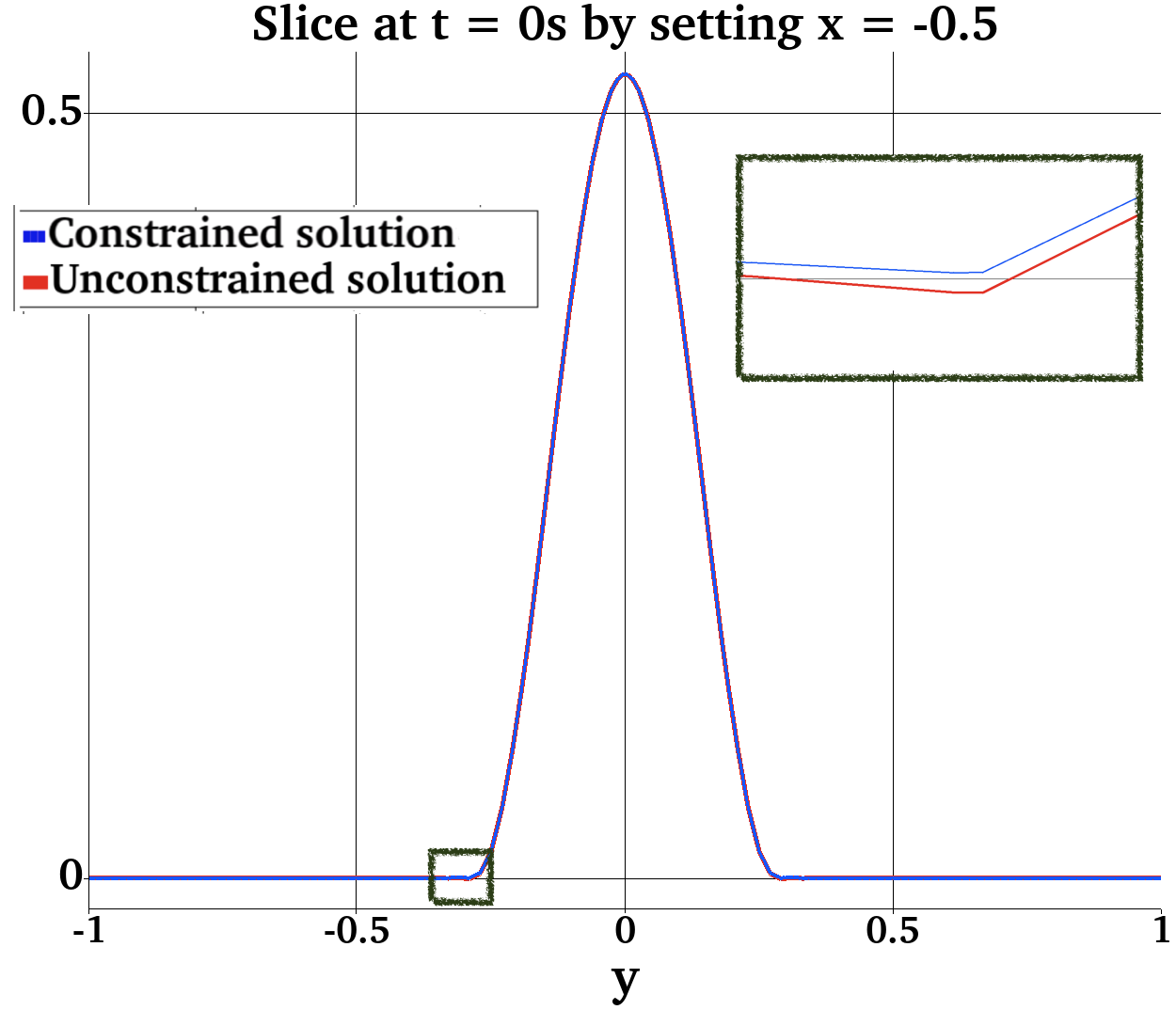}
\includegraphics[width=0.36\textwidth]{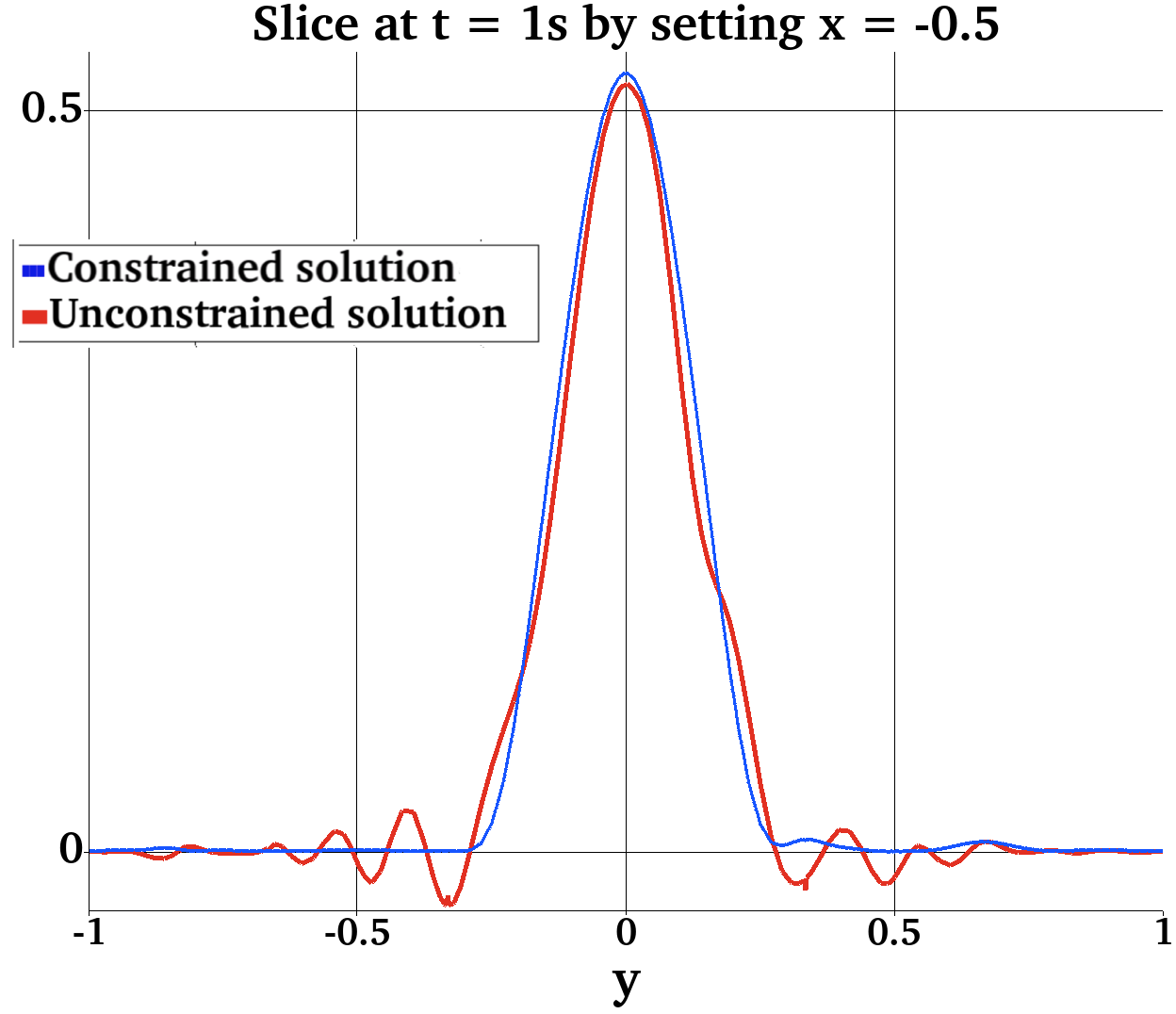}\\
\includegraphics[width=0.36\textwidth]{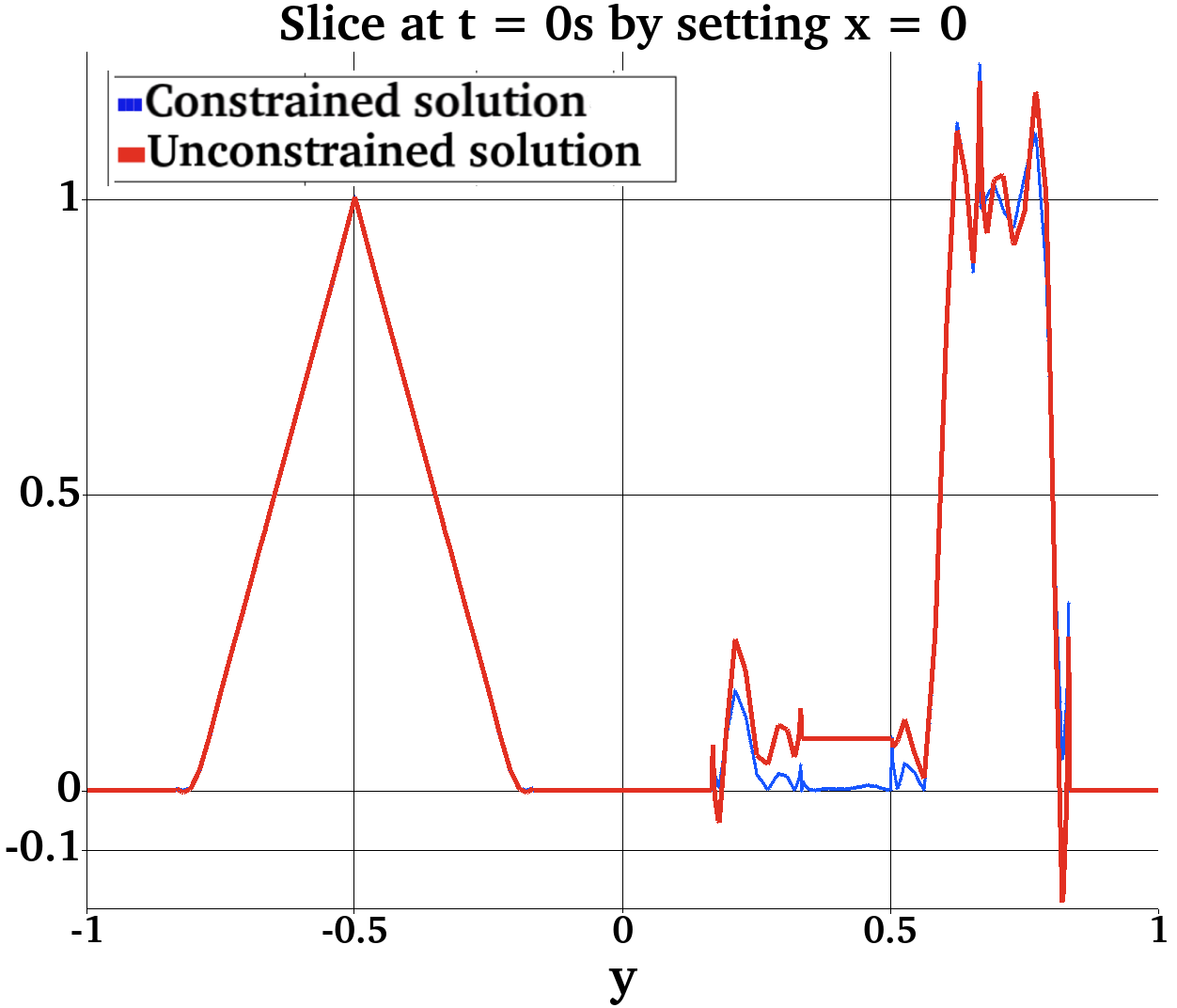}
\includegraphics[width=0.36\textwidth]{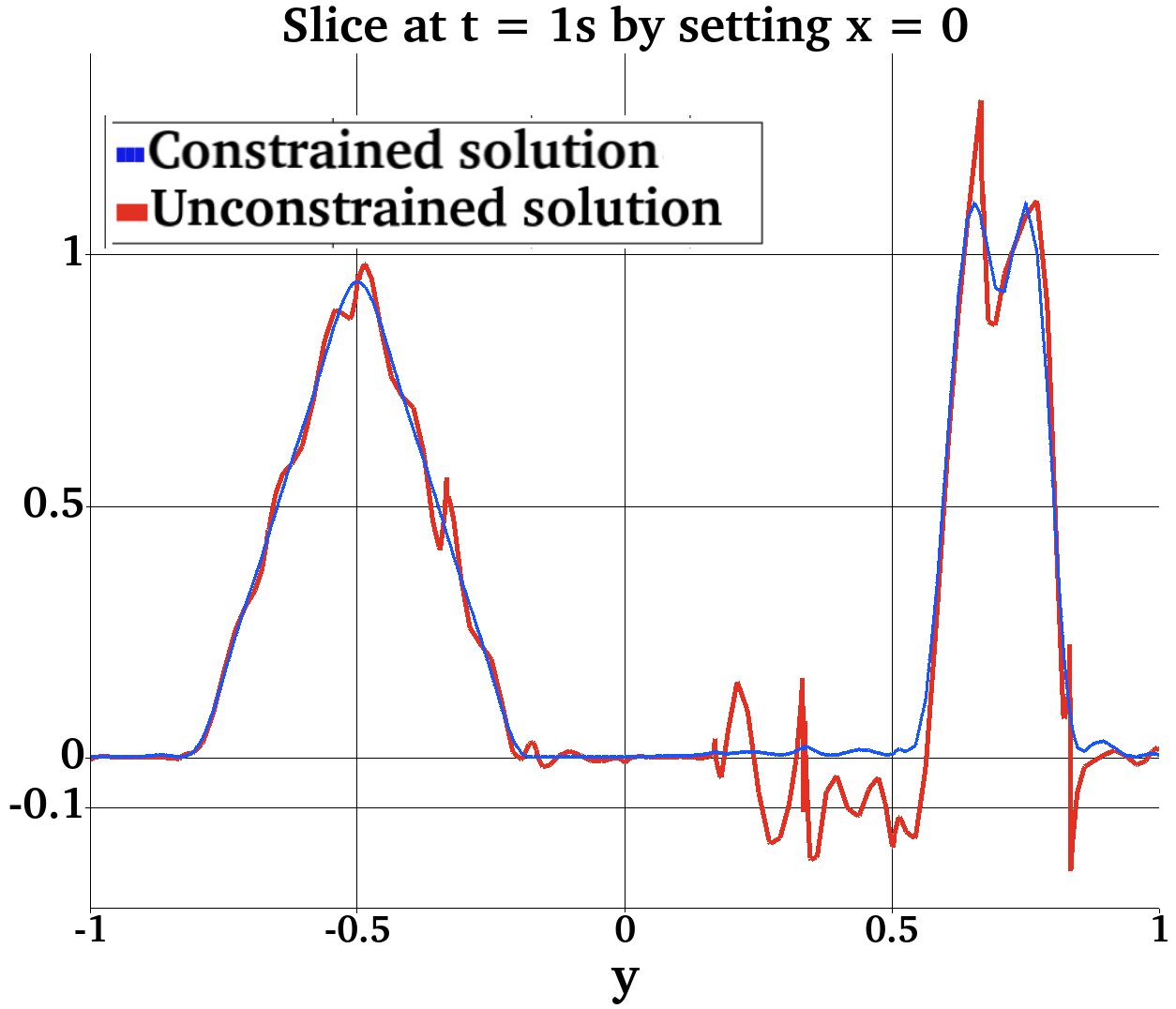}\\
\includegraphics[width=0.36\textwidth]{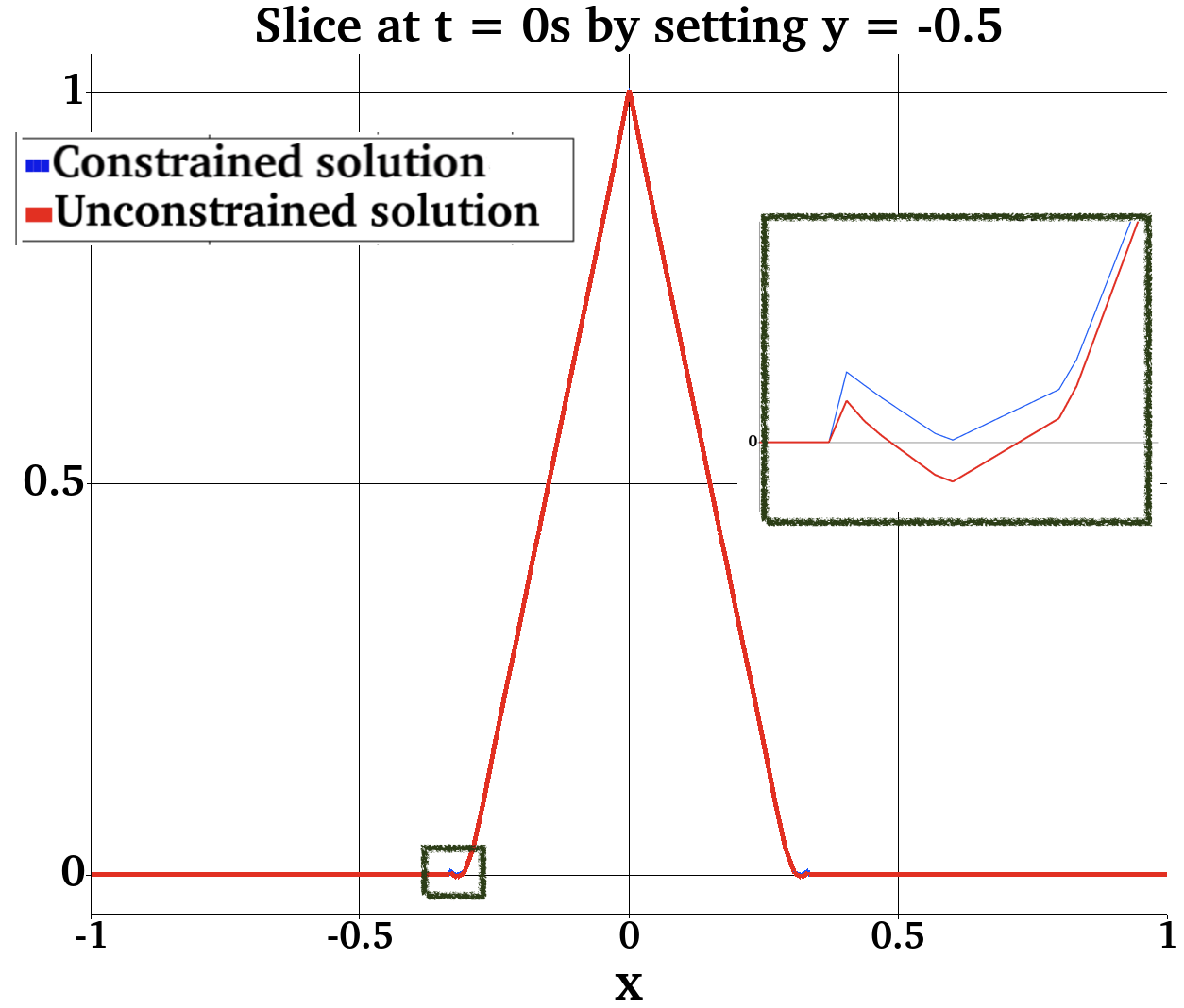}
\includegraphics[width=0.36\textwidth]{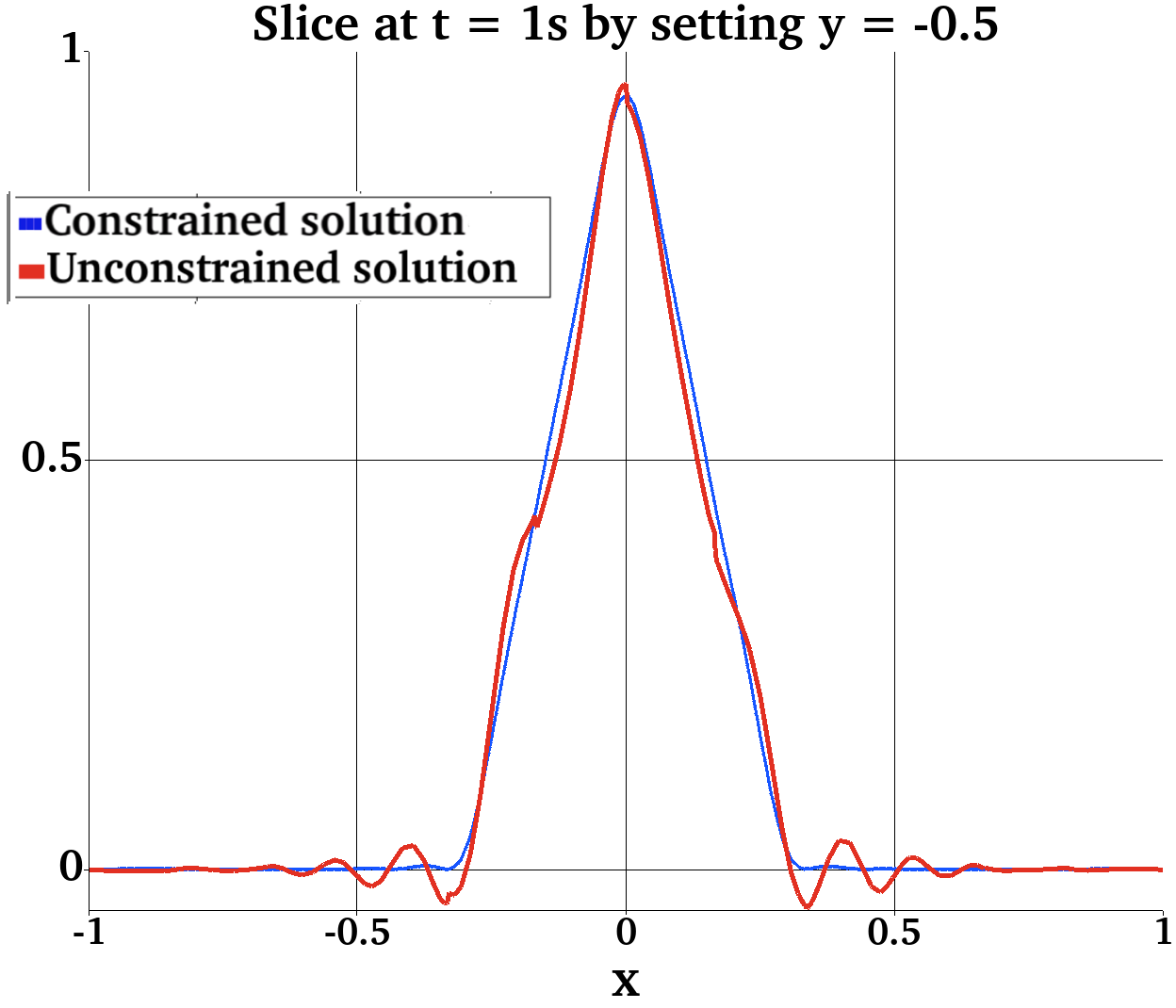}\\
\includegraphics[width=0.36\textwidth]{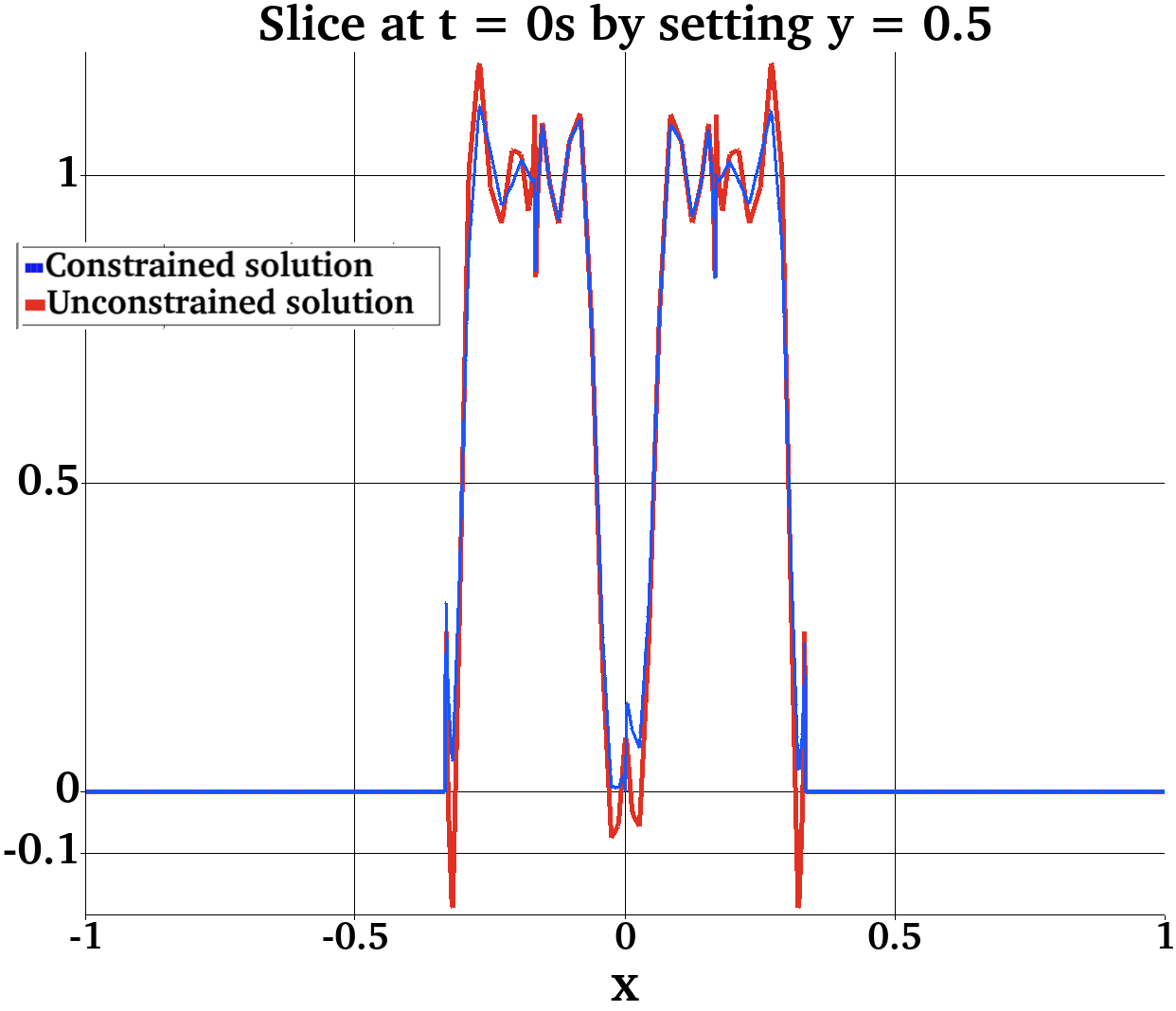}
\includegraphics[width=0.36\textwidth]{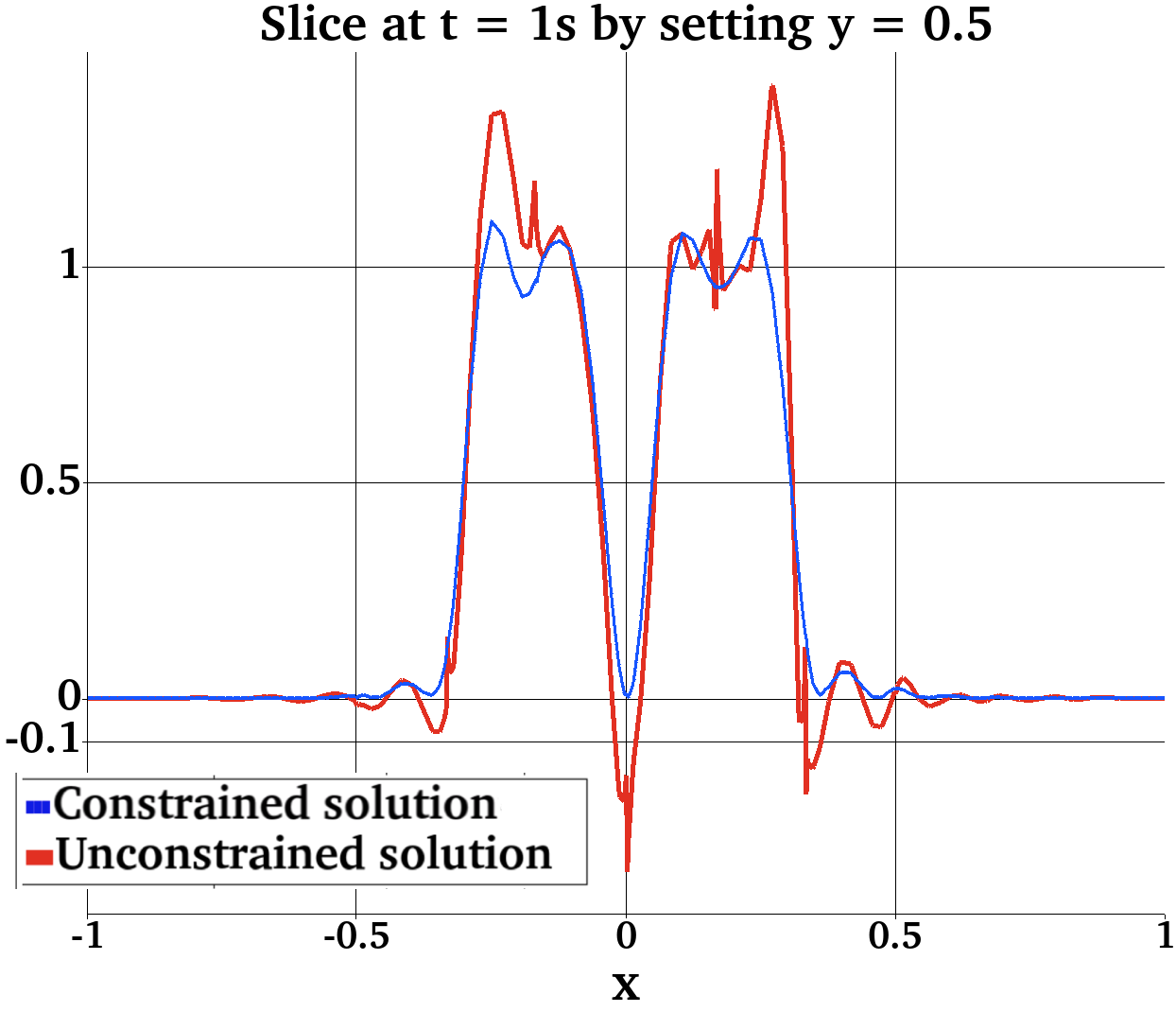}
\caption{{Snapshots at various points in time during advection. $\Delta t = 5e-4$, time period = 1 and polynomial order = 7.  Some subfigures show insets with the zoomed-in area of interest.
\textit{Column 1:} State of the system at t = 0. \textit{Column 2:} State of the system at t = 1s.  \textit{Row 1:} Slice at x = -0.5. \textit{Row 2:} Slice at x = 0. \textit{Row 3 :}Slice at y = -0.5 and \textit{Row 4: }Slice at y = 0.5}}
  \label{fig:adv3dlevequeslice1_apb2}
\end{figure}

\end{document}